\documentclass[leqno,openany,draft]{book}

\usepackage{makeidx}
\makeindex

\def\diam{\mathop{\rm diam}}
\def\dist{\mathop{\rm dist}}

\newtheorem{theorem}{Theorem}
\newtheorem{lemma}[theorem]{Lemma}
\newtheorem{proposition}[theorem]{Proposition}
\newtheorem{definition}[theorem]{Definition}
\newtheorem{corollary}[theorem]{Corollary}

\newcommand{\begintheorem}{\addtocounter{equation}{1}\begin{theorem}}
\newcommand{\beginlemma}{\addtocounter{equation}{1}\begin{lemma}}
\newcommand{\beginproposition}{\addtocounter{equation}{1}\begin{proposition}}
\newcommand{\begindefinition}{\addtocounter{equation}{1}\begin{definition}}
\newcommand{\begincorollary}{\addtocounter{equation}{1}\begin{corollary}}

\begin{document}

\frontmatter

\title{Some Basic Aspects of Analysis \\ on Metric and Ultrametric Spaces}

\author{Stephen Semmes \\
        Rice University}

\date{}

\maketitle

\chapter*{Preface}

        A number of topics involving metrics and measures are discussed,
including some of the special structure associated with ultrametrics.

\tableofcontents

\mainmatter

\chapter{Basic notions}
\label{basic notions}

\section{Metrics and ultrametrics}
\label{metrics, ultrametrics}

        Let $X$ be a set.  As usual, a \emph{metric}\index{metrics} on $X$ 
is a nonnegative real-valued function $d(x, y)$ defined for $x, y \in
X$ such that $d(x, y) = 0$ if and only if $x = y$,
\begin{equation}
\label{d(x, y) = d(y, x)}
        d(x, y) = d(y, x)
\end{equation}
for every $x, y \in X$, and
\begin{equation}
\label{d(x, z) le d(x, y) + d(y, z)}
        d(x, z) \le d(x, y) + d(y, z)
\end{equation}
for every $x, y, z \in X$.  If
\begin{equation}
\label{d(x, z) le max(d(x, y), d(y, z))}
        d(x, z) \le \max(d(x, y), d(y, z))
\end{equation}
for every $x, y, z \in X$, then $d(x, y)$ is said to be an
\emph{ultrametric}\index{ultrametrics} on $X$.

        Let $(X, d(x, y))$ be a metric space, and let $x \in X$ and a 
positive real number $r$ be given.  The corresponding 
\emph{open ball}\index{open balls} in $X$ is defined by
\begin{equation}
\label{B(x, r) = {y in X : d(x, y) < r}}
        B(x, r) = \{y \in X : d(x, y) < r\}.
\end{equation}
If $y \in B(x, r)$, then $t = r - d(x, y) > 0$, and
\begin{equation}
\label{B(y, t) subseteq B(x, r)}
        B(y, t) \subseteq B(x, r),
\end{equation}
by the triangle inequality.  However, if $d(\cdot, \cdot)$ is an
ultrametric on $X$, then one can check that
\begin{equation}
\label{B(y, r) subseteq B(x, r)}
        B(y, r) \subseteq B(x, r)
\end{equation}
for every $y \in B(x, r)$.  In fact,
\begin{equation}
\label{B(x, r) = B(y, r)}
        B(x, r) = B(y, r)
\end{equation}
for every $x, y \in X$ with $d(x, y) < r$, since we can also apply the
previous argument with the roles of $x$ and $y$ reversed.

        Similarly, the \emph{closed ball}\index{closed balls} in a metric
space $X$ centered at $x \in X$ and with radius $r \ge 0$ is defined by
\begin{equation}
\label{overline{B}(x, r) = {y in X  : d(x, y) le r}}
        \overline{B}(x, r) = \{y \in X  : d(x, y) \le r\}.
\end{equation}
If $d(\cdot, \cdot)$ is an ultrametric on $X$, and if $y \in
\overline{B}(x, r)$, then
\begin{equation}
\label{overline{B}(y, r) subseteq overline{B}(x, r)}
        \overline{B}(y, r) \subseteq \overline{B}(x, r),
\end{equation}
as before.  It follows that
\begin{equation}
\label{overline{B}(x, r) = overline{B}(y, r)}
        \overline{B}(x, r) = \overline{B}(y, r)
\end{equation}
when $d(x, y) \le r$, by reversing the roles of $x$ and $y$.

        Let us continue to ask for the moment that $d(\cdot, \cdot)$
be an ultrametric on $X$.  If $x, y, z \in X$ and $d(y, z) \le d(x, y)$,
then
\begin{equation}
\label{d(x, z) le d(x, y)}
        d(x, z) \le d(x, y),
\end{equation}
by (\ref{d(x, z) le max(d(x, y), d(y, z))}).  Of course, we also have
that
\begin{equation}
\label{d(x, y) le max(d(x, z), d(y, z))}
        d(x, y) \le \max(d(x, z), d(y, z)),
\end{equation}
by (\ref{d(x, z) le max(d(x, y), d(y, z))}) with the roles of $y$ and
$z$ exchanged.  This implies that
\begin{equation}
\label{d(x, y) le d(x, z)}
        d(x, y) \le d(x, z)
\end{equation}
when $d(y, z) < d(x, y)$, and hence that
\begin{equation}
\label{d(x, y) = d(x, z)}
        d(x, y) = d(x, z).
\end{equation}

        Put
\begin{equation}
\label{V(x, r) = {y in X : d(x, y) > r}}
        V(x, r) = \{y \in X : d(x, y) > r\}
\end{equation}
for every $x \in X$ and $r \ge 0$, which is the same as the complement
of $\overline{B}(x, r)$ in $X$.  If $d(\cdot, \cdot)$ is an ordinary
metric on $X$ and $y \in V(x, r)$, then $t = d(x, y) - r > 0$, and
one can check that
\begin{equation}
\label{B(y, t) subseteq V(x, r)}
        B(y, t) \subseteq V(x, r),
\end{equation}
using the triangle inequality.  If $d(\cdot, \cdot)$ is an ultrametric
on $X$ and $y \in V(x, r)$, then we get that
\begin{equation}
\label{B(y, d(x, y)) subseteq V(x, r)}
        B(y, d(x, y)) \subseteq V(x, r),
\end{equation}
by (\ref{d(x, y) le d(x, z)}).  Similarly,
\begin{equation}
\label{W(x, r) = {y in X : d(x, y) ge r}}
        W(x, r) = \{y \in X : d(x, y) \ge r\}
\end{equation}
is the same as the complement of $B(x, r)$ in $X$ for each $x \in X$
and $r > 0$.  If $d(\cdot, \cdot)$ is an ultrametric on $X$ and $y \in
W(x, r)$, then we also have that
\begin{equation}
\label{B(y, d(x, y)) subseteq W(x, r)}
        B(y, d(x, y)) \subseteq W(x, r),
\end{equation}
by (\ref{d(x, y) le d(x, z)}).

        If $X$ is any metric space, then every open ball in $X$ is an
open set in $X$ with respect to the topology determined by the metric.
Closed balls in $X$ are closed sets too, which is the same as saying
that $V(x, r)$ is an open set in $X$ for every $x \in X$ and $r \ge
0$.  If $d(\cdot, \cdot)$ is an ultrametric on $X$, then
(\ref{overline{B}(y, r) subseteq overline{B}(x, r)}) implies that
$\overline{B}(x, r)$ is an open set in $X$ for every $x \in X$ and $r > 0$.
In this case, $W(x, r)$ is an open set in $X$ for every $x \in X$ and $r > 0$,
by (\ref{B(y, d(x, y)) subseteq W(x, r)}), which implies that $B(x, r)$
is a closed set in $X$.

        Let $|x|$ be the absolute value of a real number $x$, which is 
equal to $x$ when $x \ge 0$ and to $-x$ when $x \le 0$.  Thus the
standard metric on the real line ${\bf R}$ is given by $|x - y|$.  Of
course, this is far from being an ultrametric.  By contrast, the
$p$-adic metric on the set ${\bf Q}$ of rational numbers is an
ultrametric for each prime number $p$.  This will be discussed in 
Section \ref{p-adic absolute value}.

\section{Abstract Cantor sets}
\label{abstract cantor sets}

        Let $X_1, X_2, X_3, \ldots$ be a sequence of finite sets, each of
which has at least two elements.  Also let $X = \prod_{j = 1}^\infty
X_j$ be their Cartesian product, which is the set of sequences $x =
\{x_j\}_{j = 1}^\infty$ such that $x_j \in X_j$ for each $j$.  Thus
$X$ is a compact Hausdorff space with respect to the product topology
corresponding to the discrete topology on each factor.  If $x, y \in
X$ and $x \ne y$, then let $l(x, y)$ be the largest nonnegative
integer such that $x_j = y_j$ when $1 \le j \le l(x, y)$.
Equivalently, $l(x, y) + 1$ is the smallest positive integer $j$ such
that $x_j \ne y_j$.  If $x = y$, then one can take $l(x, y) =
+\infty$.  Note that
\begin{equation}
\label{l(x, y) = l(y, x)}
        l(x, y) = l(y, x)
\end{equation}
for every $x, y \in X$, and that
\begin{equation}
\label{l(x, z) ge min(l(x, y), l(y, z))}
        l(x, z) \ge \min(l(x, y), l(y, z))
\end{equation}
for every $x, y, z \in X$.

        Let $\{t_l\}_{l = 0}^\infty$ be a strictly decreasing sequence
of positive real numbers that converges to $0$.  Put
\begin{equation}
\label{d(x, y) = t_{l(x, y)}}
        d(x, y) = t_{l(x, y)}
\end{equation}
when $x \ne y$, and $d(x, y) = 0$ when $x = y$, which corresponds to
(\ref{d(x, y) = t_{l(x, y)}}) with $t_\infty = 0$.  It is easy to see
that this defines an ultrametric on $X$, because of (\ref{l(x, y) =
  l(y, x)}) and (\ref{l(x, z) ge min(l(x, y), l(y, z))}), and that the
topology on $X$ determined by $d(x, y)$ is the same as the product
topology on $X$ corresponding to the discrete topology on each factor.
If $x \in X$ and $k$ is a nonnegative integer, then put
\begin{equation}
\label{B_k(x) = {y in X : y_j = x_j for each j le k}}
        B_k(x) = \{y \in X : y_j = x_j \hbox{ for each } j \le k\}.
\end{equation}
Equivalently, $B_k(x)$ is the closed ball in $X$ centered at $x$ with
radius $t_k$ with respect to (\ref{d(x, y) = t_{l(x, y)}}).

        Suppose now that $\mu_j$ is a probability measure on $X_j$
for each $j$, where all subsets of $X_j$ are measurable.  Thus $\mu_j$ 
assigns a weight to each element of $X_j$, and the sum of the weights
is equal to $1$.  This leads to a product probability measure
$\mu$ on $X$, where
\begin{equation}
\label{mu(B_k(x)) = prod_{j = 1}^k mu_j({x_j})}
        \mu(B_k(x)) = \prod_{j = 1}^k \mu_j(\{x_j\})
\end{equation}
for each $x \in X$ and $k \ge 1$.  Alternatively, one can first use the
$\mu_j$'s to define a nonnegative linear functional on the space
of continuous real-valued functions on $X$, as a limit of Riemann
sums.  One can then apply the Riesz representation theorem, to get
a Borel probability measure on $X$.

        Let $n_j \ge 2$ be the number of elements of $X_j$ for each 
positive integer $j$.  Also let $\mu_j$ be the probability measure on
$X_j$ that corresponds to the uniform distribution on $X_j$, which
assigns to each element of $X_J$ has the same weight $1/n_j$.  In this
case, (\ref{mu(B_k(x)) = prod_{j = 1}^k mu_j({x_j})}) reduces to
\begin{equation}
\label{mu(B_k(x)) = 1/N_k}
        \mu(B_k(x)) = 1/N_k
\end{equation}
for each $x \in X$ and $k \ge 1$, where
\begin{equation}
\label{N_k = prod_{j = 1}^k n_j}
        N_k = \prod_{j = 1}^k n_j.
\end{equation}
If we put $N_0 = 1$, then (\ref{mu(B_k(x)) = 1/N_k}) holds for $k = 0$
as well.  Note that $t_l = 1/N_l$ defines a strictly decreasing
sequence of positive real numbers that converges to $0$, as before.

\section{The $p$-adic absolute value}
\label{p-adic absolute value}

        Let $p$ be a prime number, and let $x$ be a rational number.
The \emph{$p$-adic absolute value}\index{p-adic absolute
 value@$p$-adic absolute value} $|x|_p$ of $x$ is defined as follows.
If $x = 0$, then $|x|_p = 0$, and otherwise $x$ can be expressed as
$p^l \, a / b$, where $a$, $b$, and $l$ are integers, and neither $a$
nor $b$ is an integer multiple of $p$.  In this case, we put
\begin{equation}
\label{|x|_p = p^{-l}}
        |x|_p = p^{-l},
\end{equation}
which is not affected by any other common factors that $a$ and $b$
might have.  It is easy to see that
\begin{equation}
\label{|x + y|_p le max(|x|_p, |y|_p)}
        |x + y|_p \le \max(|x|_p, |y|_p)
\end{equation}
and
\begin{equation}
\label{|x y|_p = |x|_p |y|_p}
        |x \, y|_p = |x|_p \, |y|_p
\end{equation}
for every $x, y \in {\bf Q}$.  The \emph{$p$-adic 
metric}\index{p-adic metric@$p$-adic metric} is defined on ${\bf Q}$ by
\begin{equation}
\label{d_p(x, y) = |x - y|_p}
        d_p(x, y) = |x - y|_p.
\end{equation}
This is an ultrametric on ${\bf Q}$, because of (\ref{|x + y|_p le
  max(|x|_p, |y|_p)}).

        If $y \in {\bf Q}$ and $n$ is a nonnegative integer, then
\begin{equation}
\label{(1 - y) sum_{j = 0}^n y^j = 1 - y^{n + 1}}
        (1 - y) \, \sum_{j = 0}^n y^j = 1 - y^{n + 1},
\end{equation}
by a standard computation.  Here $y^j$ is interpreted as being equal
to $1$ for all $y$ when $j = 0$, as usual.  If $|y|_p < 1$, then $y^{n
  + 1} \to 0$ as $n \to \infty$ with respect to the $p$-adic metric.
This implies that
\begin{equation}
\label{sum_{j = 0}^n y^j = frac{1 - y^{n + 1}}{1 - y} to frac{1}{1 - y}}
        \sum_{j = 0}^n y^j = \frac{1 - y^{n + 1}}{1 - y} \to \frac{1}{1 - y}
\end{equation}
as $n \to \infty$ with respect to the $p$-adic metric.

        Of course, $|x|_p \le 1$ for every integer $x$.  Now let 
$x \in {\bf Q}$ with $|x|_p \le 1$ be given.  Thus $x$ can be expressed as
$a / b$, where $a$ and $b$ are integers, $b \ne 0$, and $b$ is not an
integer multiple of $p$.  It is well known that there is a nonzero
integer $c$ such that $b \, c \equiv 1$ modulo $p$ under these
conditions.  Put $y = 1 - b \, c$, so that $y$ is an integer which is 
divisible by $p$, and hence $|y|_p \le 1/p < 1$.  It follows that
\begin{equation}
\label{x = frac{a}{b} = frac{a c}{b c} = frac{a c}{1 - y}}
        x = \frac{a}{b} = \frac{a \, c}{b \, c} = \frac{a \, c}{1 - y}
\end{equation}
can be approximated by integers with respect to the $p$-adic metric,
by (\ref{sum_{j = 0}^n y^j = frac{1 - y^{n + 1}}{1 - y} to frac{1}{1 - y}}).

\section{$p$-Adic numbers}
\label{p-adic numbers}

        The set ${\bf Q}_p$\index{Q_p@${\bf Q}_p$} of \emph{$p$-adic 
numbers}\index{p-adic numbers@$p$-adic numbers} can be obtained by 
completing ${\bf Q}$ as a metric space with respect to the $p$-adic
metric, in the same way that the real line ${\bf R}$ is obtained by
completing ${\bf Q}$ with respect to the standard Euclidean metric.
Sums and product of rational numbers can be extended to $p$-adic
numbers in a natural way, so that ${\bf Q}_p$ becomes a field.  The
$p$-adic absolute value $|x|_p$ and $p$-adic metric $d_p(x, y)$ can
also be extended to $x, y \in {\bf Q}_p$, in such a way that (\ref{|x
  + y|_p le max(|x|_p, |y|_p)}), (\ref{|x y|_p = |x|_p |y|_p}), and
(\ref{d_p(x, y) = |x - y|_p}) still hold.  By construction, ${\bf Q}$
is dense in ${\bf Q}_p$ with respect to the $p$-adic metric, and
$|x|_p$ is an integer power of $p$ for every $x \in {\bf Q}_p$ with $x
\ne 0$.  One can show that addition and multiplication are continuous
on ${\bf Q}_p$ with respect to the $p$-adic metric, in essentially the
same way as for real numbers.

        The set ${\bf Z}_p$\index{Z_p@${\bf Z}_p$} of \emph{$p$-adic 
integers}\index{p-adic integers@$p$-adic integers} is defined by
\begin{equation}
\label{{bf Z}_p = {x in {bf Q}_p : |x|_p le 1}}
        {\bf Z}_p = \{x \in {\bf Q}_p : |x|_p \le 1\}.
\end{equation}
This is the same as the closed unit ball in ${\bf Q}_p$, which is a
closed set in ${\bf Q}_p$ in particular.  This is also an open set in
${\bf Q}_p$ with respect to the $p$-adic metric, because the $p$-adic
metric is an ultrametric, as in Section \ref{metrics, ultrametrics}.
Of course, ${\bf Z}_p$ contains the set ${\bf Z}$ of ordinary
integers.  It is easy to see that ${\bf Q} \cap {\bf Z}_p$ is dense in
${\bf Z}_p$ with respect to the $p$-adic metric, because ${\bf Q}$ is
dense in ${\bf Q}_p$, and using the ultrametric version of the
triangle inequality.  As in the previous section, elements of ${\bf Q}
\cap {\bf Z}_p$ can be approximated by integers with respect to the
$p$-adic metric.  Combining these statements, we get that elements of
${\bf Z}_p$ can be approximated by elements of ${\bf Z}$ with respect
to the $p$-adic metric, so that ${\bf Z}_p$ is the same as the closure
of ${\bf Z}$ in ${\bf Q}_p$ with respect to the $p$-adic metric.  Note
that ${\bf Z}_p$ is also closed under addition and multiplication, by
(\ref{|x + y|_p le max(|x|_p, |y|_p)}) and (\ref{|x y|_p = |x|_p
  |y|_p}).

        Put
\begin{equation}
\label{p^l {bf Z}_p = {p^l x : x in Z_p} = {y in Q_p : |y|_p le p^{-l}}}
        p^l \, {\bf Z}_p = \{p^l \, x : x \in {\bf Z}_p\} 
                         = \{y \in {\bf Q}_p : |y|_p \le p^{-l}\}
\end{equation}
for each integer $l$.  This is the same as the closed ball in ${\bf
  Q}_p$ centered at $0$ with radius $p^{-l}$ with respect to the
$p$-adic metric, which is also an open set in ${\bf Q}_p$, as in
Section \ref{metrics, ultrametrics}.  Observe that $p^l \, {\bf Z}_p$
is a subgroup of ${\bf Q}_p$ with respect to addition for each $l$,
because of (\ref{|x + y|_p le max(|x|_p, |y|_p)}).  If $l \ge 0$, then
$p^l \, {\bf Z}_p$ is an ideal in ${\bf Z}_p$ as a commutative ring,
and hence the quotient ${\bf Z}_p / p^l \, {\bf Z}_p$ can be defined
as a commutative ring.  The composition of the obvious inclusion of
${\bf Z}$ in ${\bf Z}_p$ with the standard quotient homomorphism from
${\bf Z}_p$ onto ${\bf Z}_p / p^l \, {\bf Z}_p$ leads to a ring
homomorphism from ${\bf Z}$ into ${\bf Z}_p / p^l \, {\bf Z}_p$.  The
kernel of this homomorphism is
\begin{equation}
\label{{bf Z} cap (p^l {bf Z}_p) = p^l {bf Z}}
        {\bf Z} \cap (p^l \, {\bf Z}_p) = p^l \, {\bf Z},
\end{equation}
which is an ideal in ${\bf Z}$.  This leads to a natural injective
ring homomorphism from ${\bf Z} / p^l \, {\bf Z}$ into ${\bf Z}_p /
p^l \, {\bf Z}_p$.  Every element of ${\bf Z}_p$ can be expressed as
the sum of elements of ${\bf Z}$ and $p^l \, {\bf Z}_p$, because ${\bf
  Z}$ is dense in ${\bf Z}_p$ with respect to the $p$-adic metric.
Thus we get a natural ring isomorphism from ${\bf Z} / p^l \, {\bf Z}$
onto ${\bf Z}_p / p^l \, {\bf Z}_p$ for each nonnegative integer $l$.

        In particular, ${\bf Z}_p / p^l \, {\bf Z}_p$ has exactly $p^l$
elements for each nonnegative integer $l$.  This implies that ${\bf Z}_p$
can be expressed as the union of $p^l$ pairwise-disjoint translates of
$p^l \, {\bf Z}_p$ for each $l \ge 0$.  It follows that ${\bf Z}_p$
is totally bounded with respect to the $p$-adic metric, in the sense
that ${\bf Z}_p$ can be covered by finitely many balls of arbitrarily
small radius.  A well-known theorem implies that ${\bf Z}_p$ is compact
with respect to the topology determined on ${\bf Q}_p$ by the $p$-adic
metric, because ${\bf Z}_p$ is also a closed set in ${\bf Q}_p$
and ${\bf Q}_p$ is complete.  Of course, $p^k \, {\bf Z}_p$ is a
compact set in ${\bf Q}_p$ for every integer $k$ too, by continuity
of multiplication.

\section{Haar measure on ${\bf Q}_p$}
\label{haar measure on Q_p}

        If $A$ is a locally compact commutative topological group,
then it is well known that there is a nonnegative
translation-invariant Borel measure on $A$ which is finite on compact
subsets of $A$, positive on nonempty open subsets of $A$, and which
satisfies certain other regularity properties.  This is known as
\emph{Haar measure}\index{Haar measure} on $A$, and it is unique up to
multiplication by a positive real number.  The real line is a
commutative topological group with respect to addition and the
standard topology, for instance, and Lebesgue measure on ${\bf R}$
satisfies the requirements of Haar measure.  Similarly, the discussion
in the previous section implies that ${\bf Q}_p$ is a locally compact
commutative topological group with respect to addition and the
topology determined by the $p$-adic metric.  Let $|E|$ be the
corresponding Haar measure of a Borel set $E \subseteq {\bf Q}_p$,
normalized so that $|{\bf Z}_p| = 1$.

        If $l$ is a positive integer, then it follows that
\begin{equation}
\label{|p^l {bf Z}_p| = p^{-l}}
        |p^l \, {\bf Z}_p| = p^{-l}.
\end{equation}
This uses the fact that ${\bf Z}_p$ can be expressed as the union of
$p^l$ pairwise-disjoint translates of $p^l \, {\bf Z}_p$, as in the
previous section.  If $l$ is a negative integer, then $p^l \, {\bf
  Z}_p$ can be expressed as the union of $p^{-l}$ pairwise-disjoint
translates of ${\bf Z}_p$, by applying the previous statement to $-l$.
This implies that (\ref{|p^l {bf Z}_p| = p^{-l}}) also holds when $l <
0$, and hence for all $l \in {\bf Z}$.

        If $a \in {\bf Q}_p$ and $E \subseteq {\bf Q}_p$ is a Borel set,
then
\begin{equation}
\label{a E = {a x : x in E}}
        a \, E = \{a \, x : x \in E\}
\end{equation}
is also a Borel set in ${\bf Q}_p$.  This is trivial when $a = 0$, and
it follows from the fact that $x \mapsto a \, x$ is a homeomorphism on
${\bf Q}_p$ when $a \ne 0$.  If $a \ne 0$, then $x \mapsto a \, x$ is
an isomorphism of ${\bf Q}_p$ onto itself as a commutative topological
group, which implies that $|a \, E|$ satisfies the requirements of a
Haar measure on ${\bf Q}_p$.  The uniqueness of Haar measure implies
that $|a \, E|$ is a constant multiple of $|E|$, where the constant
depends on $a$ but not $E$.  To determine the constant, one can
consider the case where $E = {\bf Z}_p$, using (\ref{|p^l {bf Z}_p| =
  p^{-l}}).  If $|a|_p = p^{-l}$ for some $l \in {\bf Z}$, then it is
easy to see that $a \, {\bf Z}_p = p^l \, {\bf Z}_p$, so that
\begin{equation}
\label{|a {bf Z}_p| = |p^l {bf Z}_p| = p^{-l} = |a|_p}
       |a \, {\bf Z}_p| = |p^l \, {\bf Z}_p| = p^{-l} = |a|_p.
\end{equation}
It follows that
\begin{equation}
\label{|a E| = |a|_p |E|}
        |a \, E| = |a|_p \, |E|
\end{equation}
for every $a \in {\bf Q}_p$ and Borel set $E \subseteq {\bf Q}_p$,
which is trivial when $a = 0$.

        Let $A$ be a locally compact commutative topological group
again, and let $C_{com}(A)$ be the vector space of real-valued continuous
functions on $A$ with compact support.  Nonnegative linear functionals
on $C_{com}(A)$ correspond to nonnegative Borel measures on $A$ which are
finite on compact sets and have certain other regularity properties,
by the Riesz representation theorem.  The existence and uniqueness
of Haar measure on $A$ can also be considered in terms of Haar integrals,
which are nonnegative linear functionals on $C_{com}(A)$ that are invariant
under translations and positive on nonnegative elements of $C_{com}(A)$
that are positive somewhere on $A$.  The ordinary Riemann integral can
be used to define a Haar integral on the real line, for instance.
Similarly, one can get a Haar integral on ${\bf Q}_p$ as a limit
of suitable Riemann sums.

\section{Snowflake metrics and quasi-metrics}
\label{snowflake metrics, quasi-metrics}
\index{snowflake metrics}

        It is well known that
\begin{equation}
\label{(r + t)^a le r^a + t^a}
        (r + t)^a \le r^a + t^a
\end{equation}
for all nonnegative real numbers $r$, $t$ when $a \in {\bf R}$
satisfies $0 < a \le 1$.  Indeed,
\begin{equation}
\label{max(r, t) le (r^a + t^a)^{1/a}}
        \max(r, t) \le (r^a + t^a)^{1/a}
\end{equation}
for every $a > 0$, which implies that
\begin{eqnarray}
\label{ r + t le ... = (r^a + t^a)^{1/a}}
 r + t & \le & (r^a + t^a) \, \max(r, t)^{1 - a}   \\
       & \le & (r^a + t^a)^{1 + (1 - a)/a}
                               = (r^a + t^a)^{1/a} \nonumber
\end{eqnarray}
when $a \le 1$.  If $d(x, y)$ is a metric on a set $X$, then it
follows that $d(x, y)^a$ is also a metric on $X$ when $0 < a \le 1$.
Similarly, if $d(x, y)$ is an ultrametric on $X$, then $d(x, y)^a$ is
an ultrametric on $X$ for every $a > 0$.  In both cases, $d(x, y)^a$
determines the same topology on $X$ as $d(x, y)$.

        A \emph{quasi-metric}\index{quasi-metrics} on a set $X$ is a
nonnegative real-valued function $d(x, y)$ on $X \times X$ such that
$d(x, y) = 0$ if and only if $x = y$, $d(x, y) = d(y, x)$ for every
$x, y \in X$, and
\begin{equation}
\label{d(x, z) le C (d(x, y) + d(y, z))}
        d(x, z) \le C \, (d(x, y) + d(y, z))
\end{equation}
for some $C \ge 1$ and every $x, y, z \in X$.  Thus a quasi-metric
$d(x, y)$ on $X$ is a metric on $X$ if and only if one can take $C =
1$ in (\ref{d(x, z) le C (d(x, y) + d(y, z))}).  If $d(x, y)$ is a
quasi-metric on $X$, then the open ball $B(x, r)$ centered at a point
$x \in X$ with radius $r > 0$ with respect to $d(\cdot, \cdot)$ can
still be defined as in (\ref{B(x, r) = {y in X : d(x, y) < r}}).  One
can also define a topology on $X$ corresponding to $d(\cdot, \cdot)$
in the usual way, by saying that a set $U \subseteq X$ is an open set
if for each $x \in U$ there is an $r > 0$ such that $B(x, r) \subseteq
X$.  It is easy to see that this satisfies the requirements of a
topology on $X$, but the weaker version of the triangle inequality is
not sufficient to show that open balls are open sets in $X$ with
respect to this topology.  

        If $a \in {\bf R}$ and $a > 1$, then it is well known that 
$r^a$ is a convex function on the set of nonnegative real numbers $r$.
This implies that
\begin{equation}
\label{((r + t)/2)^a le (1/2) (r^a + t^a)}
        ((r + t)/2)^a \le (1/2) \, (r^a + t^a)
\end{equation}
for every $r, t \ge 0$, and hence that
\begin{equation}
\label{(r + t)^a le 2^{a - 1} (r + t)}
        (r + t)^a \le 2^{a - 1} \, (r + t).
\end{equation}
If $d(x, y)$ is a metric on a set $X$, then it follows that $d(x,
y)^a$ is a quasi-metric on $X$ for every $a > 1$.  Similarly, $d(x,
y)^a$ is a quasi-metric on $X$ for every $a > 0$ when $d(x, y)$ is a
quasi-metric on $X$.  Of course, the topology on $X$ determined
by $d(x, y)^a$ is the same as the topology on $X$ corresponding to
$d(x, y)$, and in fact the open ball in $X$ centered at a point $x \in
X$ and with radius $r > 0$ with respect to $d(\cdot, \cdot)$ is the
same as the open ball in $X$ centered at $x$ and with radius $r^a$
with respect to $d(\cdot, \cdot)^a$.

        If $d(x, y)$ is a quasi-metric on a set $X$, then one can
define a uniform structure on $X$ in the usual way, by considering the
subsets
\begin{equation}
\label{U_r = {(x, y) in X times X : d(x, y) < r}}
        U_r = \{(x, y) \in X \times X : d(x, y) < r\}
\end{equation}
of $X \times X$ for each $r > 0$.  The topology on $X$ determined by
this uniform structure is the same as the topology on $X$ defined in
terms of open balls associated to $d(x, y)$, as before.  Standard
results about uniform structures imply that for each $x \in X$ and $r
> 0$, $x$ is in the interior of the corresponding open ball $B(x, r)$
with respect to this topology.  This uniform structure on $X$
obviously has a countable base, corresponding to any sequence of
positive real numbers that converges to $0$.  This implies that there
is a metric on $X$ that determines the same uniform structure on $X$,
as in \cite{jk}.  In particular, there is a metric on $X$ that
determines the same topology on $X$ as $d(x, y)$.  In \cite{m-s-1}, it
is shown that there is a metric $\widetilde{d}(x, y)$ on $X$ and a
positive real number $a$ such that $d(x, y)$ is comparable to
$\widetilde{d}(x, y)^a$ on $X$, in the sense that each is bounded by a
constant times the other.

\section{Sequences and series}
\label{sequences, series}

        Let $d(x, y)$ be a quasi-metric on a set $X$.  As usual, a
sequence $\{x_j\}_{j = 1}^\infty$ of elements of $X$ is said to 
\emph{converge}\index{convergent sequences} to an element $x$ of $X$
if for every $\epsilon > 0$ there is an $L \ge 1$ such that
\begin{equation}
\label{d(x_j, x) < epsilon}
        d(x_j, x) < \epsilon
\end{equation}
for every $j \ge L$.  This is equivalent to saying that $\{x_j\}_{j =
  1}^\infty$ converges to $x$ with respect to the topology on $X$
determined by $d(\cdot, \cdot)$.  More precisely, this uses the fact
that every open ball in $X$ centered at $x$ contains an open set that
contains $x$ as an element, as in the previous section.  Similarly, a
sequence $\{x_j\}_{j = 1}^\infty$ of elements of $X$ is said to be a
\emph{Cauchy sequence}\index{Cauchy sequences} in $X$ if for every
$\epsilon > 0$ there is an $L \ge 1$ such that
\begin{equation}
\label{d(x_j, x_l) < epsilon}
        d(x_j, x_l) < \epsilon
\end{equation}
for every $j, l \ge L$.  In particular, it is easy to see that
convergent sequences are Cauchy sequences.  If $\{x_j\}_{j =
  1}^\infty$ is a Cauchy sequence in $X$, then
\begin{equation}
\label{lim_{j to infty} d(x_j, x_{j + 1}) = 0}
        \lim_{j \to \infty} d(x_j, x_{j + 1}) = 0,
\end{equation}
by taking $l = j + 1$ in (\ref{d(x_j, x_l) < epsilon}).  If $d(\cdot,
\cdot)$ is an ultrametric on $X$, and if $\{x_j\}_{j = 1}^\infty$ is a
sequence of elements of $X$ that satisfies (\ref{lim_{j to infty}
  d(x_j, x_{j + 1}) = 0}), then one can check that $\{x_j\}_{j =
  1}^\infty$ is a Cauchy sequence in $X$.

        Let $\sum_{j = 1}^\infty a_j$ be an infinite series whose terms
are real numbers, complex numbers, or $p$-adic numbers for some prime
number $p$.  If the corresponding sequence of partial sums
\begin{equation}
\label{s_n = sum_{j = 1}^n a_j}
        s_n = \sum_{j = 1}^n a_j
\end{equation}
converges in ${\bf R}$, ${\bf C}$, or ${\bf Q}_p$, as appropriate,
then $\sum_{j = 1}^\infty a_j$ is said to
\emph{converge},\index{convergent series} and the value of the sum is
defined to be the limit of $\{s_n\}_{n = 1}^\infty$.  Because ${\bf
  R}$, ${\bf C}$, and ${\bf Q}_p$ are complete with respect to their
standard metrics, convergence of $\sum_{j = 1}^\infty a_j$ is
equivalent to asking that $\{s_n\}_{n = 1}^\infty$ be a Cauchy
sequence.  In particular, a necessary condition for the convergence of
$\sum_{j = 1}^\infty a_j$ is that $\{a_j\}_{j = 1}^\infty$ converge as
a sequence to $0$ in ${\bf R}$, ${\bf C}$, or ${\bf Q}_p$, as
appropriate.  This is also a sufficient condition in the $p$-adic
case, because the $p$-adic metric is an ultrametric.

        An infinite series $\sum_{j = 1}^\infty a_j$ of real or complex
numbers is said to converge \emph{absolutely}\index{absolute
  convergence} if $\sum_{j = 1}^\infty |a_j|$ converges, where $|a_j|$
is the absolute value of $a_j$ in the real case, and the modulus of
$a_j$ in the complex case.  It is well known that absolute convergence
implies convergence, using the triangle inequality to show that the
partial sums of $\sum_{j = 1}^\infty a_j$ form a Cauchy sequence when
the partial sums of $\sum_{j = 1}^\infty |a_j|$ form a Cauchy
sequence.  One can also check that
\begin{equation}
        \biggl|\sum_{j = 1}^\infty a_j \biggr| \le \sum_{j = 1}^\infty |a_j|
\end{equation}
when $\sum_{j = 1}^\infty a_j$ converges absolutely.  Similarly,
if $\{a_j\}_{j = 1}^\infty$ is a sequence of $p$-adic numbers that
converges to $0$, then
\begin{equation}
\label{|sum_{j = 1}^infty a_j|_p le max_{j ge 1} |a_j|_p}
 \biggl|\sum_{j = 1}^\infty a_j\biggr|_p \le \max_{j \ge 1} |a_j|_p.
\end{equation}
Note that the maximum of $|a_j|_p$ over $j \in {\bf Z}_+$ exists in
this situation, because $|a_j|_p \to 0$ as $j \to \infty$.

        The \emph{Cauchy product}\index{Cauchy products} of two
infinite series $\sum_{j = 0}^\infty a_j$, $\sum_{k = 0}^\infty b_k$
of real, complex, or $p$-adic numbers is the infinite series 
$\sum_{l = 0}^\infty c_l$, where
\begin{equation}
\label{c_l = sum_{j = 0}^l a_j b_{l - j}}
        c_l = \sum_{j = 0}^l a_j \, b_{l - j}.
\end{equation}
If $\sum_{j = 0}^\infty a_j$ and $\sum_{k = 0}^\infty b_k$ are
absolutely convergent series of real or complex numbers, then it is
well known that $\sum_{l = 0}^\infty c_l$ converges absolutely too,
and that
\begin{equation}
\label{sum_{l = 0}^infty c_l = (sum_{j = 0}^infty a_j) (sum_{k = 0}^infty b_k)}
        \sum_{l = 0}^\infty c_l = \Big(\sum_{j = 0}^\infty a_j\Big)
                                 \, \Big(\sum_{k = 0}^\infty b_k\Big).
\end{equation}
Similarly, if $\{a_j\}_{j = 0}^\infty$ and $\{b_k\}_{k = 0}^\infty$
are sequences of $p$-adic numbers converging to $0$, then one can
check that $\{c_l\}_{l = 0}^\infty$ also converges to $0$ in ${\bf
  Q}_p$, using the fact that the $p$-adic metric is an ultrametric.
It is not too difficult to verify that (\ref{sum_{l = 0}^infty c_l =
  (sum_{j = 0}^infty a_j) (sum_{k = 0}^infty b_k)}) holds under these
conditions as well.

\chapter{Hausdorff measures}
\label{hausdorff measures}

\section{Diameters}
\label{diameters}

        Let $(M, d(x, y))$ be a metric space.  As usual, the 
\emph{diameter}\index{diameter of a set} of a nonempty bounded set
$A \subseteq M$ is defined by
\begin{equation}
\label{diam A = sup {d(x, y) : x, y in A}}
        \diam A = \sup \{d(x, y) : x, y \in A\}.
\end{equation}
It is sometimes convenient to define the diameter of the empty set
to be $0$, and to put $\diam A = \infty$ when $A$ is not bounded.
Note that
\begin{equation}
\label{diam overline{A} = diam A}
        \diam \overline{A} = \diam A
\end{equation}
for any set $A \subseteq M$, where $\overline{A}$ is the closure of
$A$ in $M$.

        Let $A \subseteq M$ and $r > 0$ be given, and put
\begin{equation}
\label{A_r = bigcup_{x in A} B(x, r)}
        A_r = \bigcup_{x \in A} B(x, r).
\end{equation}
Thus $A_r$ is an open set in $M$, since it is a union of open sets,
and $A \subseteq A_r$.  If $w, z \in A_r$, then there are $x, y \in A$
such that $d(x, w), d(y, z) < r$, and hence
\begin{equation}
\label{d(w, z) < d(x, y) + 2 r}
        d(w, z) < d(x, y) + 2 \, r.
\end{equation}
This implies that
\begin{equation}
\label{diam A_r le diam A + 2 r}
        \diam A_r \le \diam A + 2 \, r
\end{equation}
for each $r > 0$.  

        If $d(x, y)$ is an ultrametric on $M$, then we can replace
(\ref{d(w, z) < d(x, y) + 2 r}) with
\begin{equation}
\label{d(w, z) le max(d(x, y), r)}
        d(w, z) \le \max(d(x, y), r),
\end{equation}
so that
\begin{equation}
\label{diam A_r le max(diam A, r)}
        \diam A_r \le \max(\diam A, r)
\end{equation}
for each $r > 0$.  Of course,
\begin{equation}
\label{diam A le diam A_r}
        \diam A \le \diam A_r
\end{equation}
for every $r > 0$, since $A \subseteq A_r$.  Thus (\ref{diam A_r le
  max(diam A, r)}) implies that
\begin{equation}
\label{diam A_r = diam A}
        \diam A_r = \diam A
\end{equation}
when $r \le \diam A$.

        Similarly,
\begin{equation}
\label{diam overline{B}(x, r) le 2 r}
        \diam \overline{B}(x, r) \le 2 \, r
\end{equation}
for every $x \in M$ and $r \ge 0$, and for any metric $d(x, y)$ on
$M$.  If $d(x, y)$ is an ultrametric on $M$, then
\begin{equation}
\label{diam overline{B}(x, r) le r}
        \diam \overline{B}(x, r) \le r
\end{equation}
for every $x \in M$ and $r \ge 0$.  If $d(x, y)$ is any metric on $M$
and $A$ is a nonempty bounded set in $M$, then
\begin{equation}
\label{A subseteq overline{B}(x, diam A)}
        A \subseteq \overline{B}(x, \diam A)
\end{equation}
for every $x \in A$.  If $M$ is the real line with the standard
metric, and if $A$ is a nonempty bounded subset of ${\bf R}$, then
\begin{equation}
\label{I_A = [inf A, sup A]}
        I_A = [\inf A, \sup A]
\end{equation}
contains $A$ and has the same diameter as $A$.

\section{Hausdorff content}
\label{hausdorff content}

        Let $(M, d(x, y))$ be a metric space again, and let $\alpha$
be a positive real number.  The $\alpha$-dimensional 
\emph{Hausdorff content}\index{Hausdorff content} of $E \subseteq M$ 
is defined by
\begin{equation}
\label{H^alpha_{con}(E) = ...}
        H^\alpha_{con}(E) = \inf \bigg\{\sum_j (\diam A_j)^\alpha : 
                                       E \subseteq \bigcup_j A_j\bigg\},
\end{equation}
where more precisely the infimum is taken over all collections
$\{A_j\}_j$ of finitely or countably many subsets of $M$ such that $E
\subseteq \bigcup_j A_j$.  The sum
\begin{equation}
\label{sum_j (diam A_j)^alpha}
        \sum_j (\diam A_j)^\alpha
\end{equation}
is defined as usual as the supremum over all finite subsums when there
are infinitely many $A_j$'s, which may be infinite.  If $A_j$ is
unbounded for any $j$, then $\diam A_j = \infty$, and (\ref{sum_j
  (diam A_j)^alpha}) is infinite.  This definition can also be used
when $\alpha = 0$, with the conventions that $(\diam A)^0$ is equal to
$0$ when $A = \emptyset$, is equal to $1$ when $A$ is nonempty and
bounded, and is equal to $\infty$ when $A$ is unbounded.

        Note that $H^\alpha_{con}(\emptyset) = 0$ for every $\alpha \ge 0$,
and that
\begin{equation}
\label{H^alpha_{con}(E) le (diam E)^alpha}
        H^\alpha_{con}(E) \le (\diam E)^\alpha
\end{equation}
for every $E \subseteq M$ and $\alpha \ge 0$, by covering $E$ by itself.
If $E \subseteq \widetilde{E} \subseteq M$, then
\begin{equation}
\label{H^alpha_{con}(E) le H^alpha_{con}(widetilde{E})}
        H^\alpha_{con}(E) \le H^\alpha_{con}(\widetilde{E})
\end{equation}
for every $\alpha \ge 0$, because every covering of $\widetilde{E}$ in
$M$ is also a covering of $E$.  If $E_1, E_2, E_2, \ldots$ is any
sequence of subsets of $M$, then one can show that
\begin{equation}
\label{H^alpha_{con}(bigcup_{k = 1}^infty E_k) le ...}
        H^\alpha_{con}\Big(\bigcup_{k = 1}^\infty E_k\Big)
                      \le \sum_{k = 1}^\infty H^\alpha_{con}(E_k)
\end{equation}
for every $\alpha \ge 0$, by combining coverings of the $E_k$'s to get
coverings of $\bigcup_{k = 1}^\infty E_k$.  Of course, if
$H^\alpha_{con}(E_k) = \infty$ for some $k$, then the sum on the right
side of (\ref{H^alpha_{con}(bigcup_{k = 1}^infty E_k) le ...}) is
equal to $\infty$ too, in which case the inequality is trivial.
Otherwise, one can choose coverings of the $E_k$'s for which the
corresponding sums (\ref{sum_j (diam A_j)^alpha}) are as close as one
wants to $H^\alpha_{con}(E_k)$.  The main point is to do this in such
a way that the sum of the errors is arbitrarily small too.

        In the definition of the Hausdorff content, one might as well
restrict one's attention to coverings of $E$ by collections of closed
subsets of $M$, because of (\ref{diam overline{A} = diam A}).  One can
also restrict one's attention to coverings by collection of open
subsets of $M$, using (\ref{diam A_r le diam A + 2 r}).  If $E
\subseteq M$ is compact, then it follows that one can restrict one's
attention to coverings of $E$ by finitely many subsets of $M$.

        Remember that an \emph{outer measure}\index{outer measures}
on a $\sigma$-algebra $\mathcal{A}$ of subsets of $M$ is a nonnegative
extended real-valued function $\mu$ on $\mathcal{A}$ such that
$\mu(\emptyset) = 0$,
\begin{equation}
\label{mu(A) le mu(B)}
        \mu(A) \le \mu(B)
\end{equation}
for every $A, B \in \mathcal{A}$ with $A \subseteq B$, and $\mu$ is
countably-subadditive on $\mathcal{A}$.  Thus $H^\alpha_{con}$ is an
outer measure on the $\sigma$-algebra of all subsets of $M$ for each
$\alpha \ge 0$, for instance.  Let $\mu$ be an outer measure
defined on a $\sigma$-algebra $\mathcal{A}$ of subsets of $M$ that
contains the Borel sets, and suppose that
\begin{equation}
\label{mu(A) le C (diam A)^alpha}
        \mu(A) \le C \, (\diam A)^\alpha
\end{equation}
for some nonnegative real numbers $C$, $\alpha$ and every $A \in
\mathcal{A}$.  If $E \in \mathcal{A}$, and if $\{A_j\}_j$ are finitely
or countably many elements of $\mathcal{A}$ such that $E \subseteq
\bigcup_j A_j$, then
\begin{equation}
\label{mu(E) le sum_j mu(A_j) le C sum_j (diam A_j)^alpha}
        \mu(E) \le \sum_j \mu(A_j) \le C \, \sum_j (\diam A_j)^\alpha.
\end{equation}
This implies that
\begin{equation}
\label{mu(E) le C H^alpha_{con}(E)}
        \mu(E) \le C \, H^\alpha_{con}(E),
\end{equation}
since we can restrict our attention to coverings of $E$ by open or
closed subsets of $M$ in the definition of Hausdorff content, as in
the previous paragraph.

\section{Restricting the diameters}
\label{restricting the diameters}

        Let $(M, d(x, y))$ be a metric space, and let $0 \le \alpha < \infty$
and $0 < \delta \le \infty$ be given.  Put
\begin{equation}
\label{H^alpha_delta(E) = ...}
  \quad H^\alpha_\delta(E) = \inf \bigg\{\sum_j (\diam A_j)^\alpha : 
 E \subseteq \bigcup_j A_j, \, \diam A_j < \delta \hbox{ for each } j\bigg\}
\end{equation}
for each $E \subseteq M$, where more precisely the infimum is taken
over all collections $\{A_j\}_j$ of finitely or countably many subsets
of $M$ such that $E \subseteq \bigcup_j A_j$ and $\diam A_j < \delta$
for each $j$, if there are any.  If not, then put $H^\alpha_\delta(E)
= \infty$.  Of course, if $M$ is separable, then $M$ is contained in
the union of finitely or countably many balls of radius $r$ for every
$r > 0$, and this is not a problem.  This is also not a problem when
$\delta = \infty$, because every $E \subseteq M$ is covered by a
sequence of bounded subsets of $M$.

        By construction,
\begin{equation}
\label{H^alpha_{con}(E) le H^alpha_delta(E) le H^alpha_{eta}(E)}
        H^\alpha_{con}(E) \le H^\alpha_\delta(E) \le H^\alpha_{\eta}(E)
\end{equation}
for every $\alpha \ge 0$ and $E \subseteq M$ when $0 < \eta < \delta
\le \infty$, since one is restricting the class of admissible
coverings of $E$ as $\delta$ decreases.  It is easy to see that
\begin{equation}
\label{H^alpha_{con}(E) = H^alpha_infty(E)}
        H^\alpha_{con}(E) = H^\alpha_\infty(E)
\end{equation}
for every $\alpha \ge 0$ and $E \subseteq M$, because (\ref{sum_j
  (diam A_j)^alpha}) is infinite when $A_j$ is unbounded for any $j$.
As before, $H^\alpha_\delta(\emptyset) = 0$ for every $\alpha \ge 0$
and $\delta > 0$, and
\begin{equation}
\label{H^alpha_delta(E) le H^alpha_delta(widetilde{E})}
        H^\alpha_\delta(E) \le H^\alpha_\delta(\widetilde{E})
\end{equation}
when $E \subseteq \widetilde{E} \subseteq M$.  If $E_1, E_2, E_3,
\ldots$ is any sequence of subsets of $M$, then
\begin{equation}
\label{H^alpha_delta(bigcup_{k = 1}^infty E_k) le ...}
        H^\alpha_\delta\Big(\bigcup_{k = 1}^\infty E_k\Big)
                       \le \sum_{k = 1}^\infty H^\alpha_\delta(E_k)
\end{equation}
for every $\alpha \ge 0$ and $\delta > 0$, as in the previous section.
Thus $H^\alpha_\delta$ is an outer measure on the $\sigma$-algebra of all
subsets of $M$ for each $\alpha \ge 0$ and $\delta > 0$.

        One might as well restrict one's attention to coverings of $E$
by open or closed subsets of $M$ in (\ref{H^alpha_delta(E) = ...}),
for the same reasons as before.  In particular, if $E$ is compact,
then one can restrict one's attention to coverings of $E$ by finitely
many subsets of $M$.

        Suppose that $E_1, E_2 \subseteq M$ have the property that 
\begin{equation}
\label{d(x, y) ge delta}
        d(x, y) \ge \delta
\end{equation}
for some $\delta > 0$ and every $x \in E_1$ and $y \in E_2$.  Let
$\{A_j\}_{j \in I}$ be any collection of finitely or countably many
subsets of $M$ such that $\diam A_j < \delta$ for each $j$ and
\begin{equation}
\label{E_1 cup E_2 subseteq bigcup_{j in I} A_j}
        E_1 \cup E_2 \subseteq \bigcup_{j \in I} A_j.
\end{equation}
Let $I_1$, $I_2$ be the set of $j \in I$ such that $A_j$ intersects $E_1$,
$E_2$, respectively.  The separation condition (\ref{d(x, y) ge delta})
implies that $I_1$ and $I_2$ disjoint subsets of $I$, so that
\begin{eqnarray}
\label{H^alpha_delta(E_1) + H^alpha_delta(E_2) le ...}
        H^\alpha_\delta(E_1) + H^\alpha_\delta(E_2) & \le &
   \sum_{j \in I_1} (\diam A_j)^\alpha + \sum_{j \in I_2} (\diam A_j)^\alpha \\
                        & \le & \sum_{j \in I} (\diam A_j)^\alpha. \nonumber
\end{eqnarray}
for every $\alpha \ge 0$.  This implies that
\begin{equation}
\label{H^alpha_delta(E_1) + H^alpha_delta(E_2) le H^alpha_delta(E_1 cup E_2)}
 H^\alpha_\delta(E_1) + H^\alpha_\delta(E_2) \le H^\alpha_\delta(E_1 \cup E_2)
\end{equation}
for every $\alpha \ge 0$, by taking the infimum over all such coverings
$\{A_j\}_{j \in I}$ of $E_1 \cup E_2$.  The opposite inequality holds
automatically, as in (\ref{H^alpha_delta(bigcup_{k = 1}^infty E_k) le ...}).
Thus
\begin{equation}
\label{H^alpha_delta(E_1) + H^alpha_delta(E_2) = H^alpha_delta(E_1 cup E_2)}
        H^\alpha_\delta(E_1) + H^\alpha_\delta(E_2) = H^\alpha_\delta(E_1 \cup E_2)
\end{equation}
for all $\alpha \ge 0$ under these conditions.

\section{Hausdorff measures}
\label{hausdorff measures, section}

        Let $(M, d(x, y))$ be a metric space, and let $\alpha \ge 0$
be given.  The $\alpha$-dimensional \emph{Hausdorff 
measure}\index{Hausdorff measure} of $E \subseteq M$ is defined by
\begin{equation}
\label{H^alpha(E) = sup_{delta > 0} H^alpha_delta(E)}
        H^\alpha(E) = \sup_{\delta > 0} H^\alpha_\delta(E),
\end{equation}
where $H^\alpha_\delta(E)$ is as in the previous section.  This can
also be considered as the limit of $H^\alpha_\delta(E)$ as $\delta \to
0$, since $H^\alpha_\delta(E)$ increases monotonically as $\delta$
decreases.  As usual, $H^\alpha(\emptyset) = 0$ for every $\alpha \ge
0$, and
\begin{equation}
\label{H^alpha(E) le H^alpha(widetilde{E})}
        H^\alpha(E) \le H^\alpha(\widetilde{E})
\end{equation}
for every $\alpha \ge 0$ when $E \subseteq \widetilde{E} \subseteq M$.
If $E_1, E_2, E_3, \ldots$ is any sequence of subsets of $M$, then
\begin{equation}
\label{H^alpha(bigcup_{k = 1}^infty E_k) le sum_{k = 1}^infty H^alpha(E_k)}
        H^\alpha\Big(\bigcup_{k = 1}^\infty E_k\Big)
                 \le \sum_{k = 1}^\infty H^\alpha(E_k)
\end{equation}
for every $\alpha \ge 0$, by (\ref{H^alpha_delta(bigcup_{k = 1}^infty
  E_k) le ...}), so that $H^\alpha$ is an outer measure on the
$\sigma$-algebra of all subsets of $M$ for each $\alpha \ge 0$.

        Let $E \subseteq M$, $0 \le \alpha < \beta$, and $\delta > 0$
be given.  If $\{A_j\}_j$ is a collection of finitely or countably many
subsets of $M$ such that $E \subseteq \bigcup_j A_j$ and $\diam A_j < \delta$
for each $j$, then
\begin{equation}
\label{sum_j (diam A_j)^beta le delta^{beta - alpha} sum_j (diam A_j)^alpha}
 \sum_j (\diam A_j)^\beta \le \delta^{\beta - \alpha} \, \sum_j (\diam A_j)^\alpha.
\end{equation}
This implies that
\begin{equation}
\label{H^beta_delta(E) le delta^{beta - alpha} H^alpha_delta(E)}
        H^\beta_\delta(E) \le \delta^{\beta - \alpha} \, H^\alpha_\delta(E).
\end{equation}
If $H^\alpha(E) < \infty$, then one can pass to the limit as $\delta
\to 0$, to get that $H^\beta(E) = 0$.  The \emph{Hausdorff
  dimension}\index{Hausdorff dimension} of $E \subseteq M$ may be
defined as the infimum of the $\alpha \ge 0$ such that $H^\alpha(E) <
\infty$, if there is such an $\alpha$, and otherwise the Hausdorff
dimension of $E$ is $\infty$.  If $H^\alpha_{con}(E) = 0$ for some $\alpha \ge 0$
and $E \subseteq M$, then it is easy to see that $H^\alpha_\delta(E) =
0$ for every $\delta > 0$, and hence that $H^\alpha(E) = 0$.  The main
point is that if $\{A_j\}_j$ is a collection of finitely or countable
many subsets of $M$ such that $E \subseteq \bigcup_j A_j$ and the
corresponding sum (\ref{sum_j (diam A_j)^alpha}) is small, then $\diam
A_j$ has to be small for each $j$.

        Suppose that $H^\alpha(E) < \infty$ for some $\alpha \ge 0$ again.
Thus for each positive integer $n$ there is a collection $\{A_{j,
  n}\}_{j \in I_n}$ of finitely or countably many open subsets of $M$
such that $E \subseteq \bigcup_{j \in I_n} A_{j, n}$, $\diam A_{j, n}
< 1/n$ for every $j \in I_n$, and
\begin{equation}
\label{sum_{j in I_n} (diam A_{j, n})^alpha < H^alpha(E) + 1/n}
        \sum_{j \in I_n} (\diam A_{j, n})^\alpha < H^\alpha(E) + 1/n.
\end{equation}
Put
\begin{equation}
\label{widetilde{E} = bigcap_{n = 1}^infty (bigcup_{j in I_n} A_{j, n})}
 \widetilde{E} = \bigcap_{n = 1}^\infty \Big(\bigcup_{j \in I_n} A_{j, n}\Big),
\end{equation}
so that $E \subseteq \widetilde{E}$, $\widetilde{E}$ is the
intersection of a sequence of open subsets of $M$, and $\widetilde{E}
\subseteq \bigcup_{j \in I_n} A_j$ for each $n$.  The latter implies
that $H^\alpha(\widetilde{E}) \le H^\alpha(E)$, and hence that
$H^\alpha(\widetilde{E}) = H^\alpha(E)$.

        If $E_1, E_2 \subseteq M$ have the property that
\begin{equation}
\label{d(x, y) ge eta}
        d(x, y) \ge \eta
\end{equation}
for some $\eta > 0$ and every $x \in E_1$ and $y \in E_2$, then
(\ref{H^alpha_delta(E_1) + H^alpha_delta(E_2) = H^alpha_delta(E_1 cup
  E_2)}) holds when $0 < \delta \le \eta$.  This implies that
\begin{equation}
\label{H^alpha(E_1) + H^alpha(E_2) = H^alpha(E_1 cup E_2)}
        H^\alpha(E_1) + H^\alpha(E_2) = H^\alpha(E_1 \cup E_2),
\end{equation}
by taking the limit as $\delta \to 0$ in (\ref{H^alpha_delta(E_1) +
  H^alpha_delta(E_2) = H^alpha_delta(E_1 cup E_2)}).  This shows that
$H^\alpha$ satisfies a well-known criterion of Carath\'eodory, and
hence that $H^\alpha$ is countably additive on a suitable
$\sigma$-algebra of measurable subsets of $M$ that includes the Borel
sets.  If $\alpha = 0$, then Hausdorff measure reduces to counting
measure on $M$.

\section{Some special cases}
\label{some special cases}

        Suppose that $M$ is the real line, with the standard metric.
As in Section \ref{diameters}, every nonempty bounded subset of ${\bf R}$
is contained in a closed interval with the same diameter.  This implies
that one may as well restrict one's attention to coverings of 
$E \subseteq {\bf R}$ by closed intervals in the definition of
$H^\alpha_\delta(E)$ for every $\alpha, \delta > 0$.  If $\alpha = 0$,
then one should consider the empty set as a closed interval too.
One might also consider the real line itself as a closed interval,
for the analogous statement for $H^\alpha_{con}(E)$, although this
does not really matter.

        Let us restrict our attention now to $\alpha = 1$.  It is easy
to see that
\begin{equation}
\label{H^1_delta(E) = H^1_{con}(E)}
        H^1_\delta(E) = H^1_{con}(E)
\end{equation}
for every $\delta > 0$ and $E \subseteq {\bf R}$, by subdividing intervals
in ${\bf R}$ into finitely many arbitrarily small subintervals.  It follows
that
\begin{equation}
\label{H^1(E) = H^1_{con}(E)}
        H^1(E) = H^1_{con}(E)
\end{equation}
for every $E \subseteq {\bf R}$, which is of course the same as the
Lebesgue outer measure of $E$.  If $E$ is a closed interval in ${\bf
  R}$, then this is less than or equal to the diameter of $E$, which
is the same as the length of $E$ as an interval, as in
(\ref{H^alpha_{con}(E) le (diam E)^alpha}).  As usual, one can show
that $H^1(E)$ is equal to the diameter of $E$ in this case, by
considering coverings of $E$ by finitely many intervals in ${\bf R}$.

        Suppose now that $(M, d(x, y))$ is an ultrametric space.
Every nonempty bounded subset of $M$ is contained in a closed ball in
$M$ with the same diameter, as in Section \ref{diameters} again.  Thus
one may as well restrict one's attention to coverings of $E \subseteq
M$ by closed balls in $M$ in the definition of $H^\alpha_\delta(E)$
for every $\alpha, \delta > 0$.  As before, one should consider the
empty set as a closed ball in $M$ when $\alpha = 0$, and one might
also consider $M$ as a closed ball even when $M$ is unbounded, in the
context of Hausdorff content.

        In particular, these remarks can be applied to ${\bf Q}_p$,
with the $p$-adic metric.  Let us restrict our attention to $\alpha =
1$ again.  Remember that ${\bf Z}_p$ can be expressed as the union of
$p^l$ pairwise-disjoint translates of $p^l \, {\bf Z}_p$ for every
positive integer $l$, as in Section \ref{p-adic numbers}.  This
implies that any closed ball in ${\bf Q}_p$ of radius $p^k$ for some
$k \in {\bf Z}$ can be expressed as the pairwise-disjoint union of
$p^l$ closed balls of radius $p^{k - l}$ for every positive integer
$l$.  Note that every closed ball in ${\bf Q}_p$ of radius $p^j$ for
some $j \in {\bf Z}$ has diameter equal to $p^j$ too.  Using this, one
can check that (\ref{H^1_delta(E) = H^1_{con}(E)}) also holds in this
case for every $E \subseteq {\bf Q}_p$.  This implies that
(\ref{H^1(E) = H^1_{con}(E)}) holds for every $E \subseteq {\bf Q}_p$
as well.

        If $B$ is a closed ball in ${\bf Q}_p$ with radius $p^j$ for some 
$j \in {\bf Z}$, then $H^1(B) \le p^j$, by the previous discussion.
Let us verify that $H^1(B) \ge p^j$ under these conditions, and hence that
\begin{equation}
\label{H^1(B) = p^j}
        H^1(B) = p^j.
\end{equation}
To do this, it suffices to consider coverings of $B$ by finitely many
closed balls in ${\bf Q}_p$, because $B$ is compact, and closed balls
in ${\bf Q}_p$ are open sets.  More precisely, it suffices to consider
coverings of $B$ by finitely many closed balls of the same radius
$p^{j - l}$ for some nonnegative integer $l$, by subdividing balls of
different radii to get balls of the same radius.  To show that
$H^1(B) \ge p^j$, it is enough to check that $B$ cannot be covered
by fewer than $p^l$ closed balls of radius $p^{j - l}$, for any
nonnegative integer $l$.  If $B = {\bf Z}_p$, which is the closed
unit ball in ${\bf Q}_p$, then this follows from the discussion
in Section \ref{p-adic numbers}.  If $B$ is any other closed ball
in ${\bf Q}_p$, then one can reduce to the case of ${\bf Z}_p$, using
translations and dilations.

        Now let $X_1, X_2, X_3, \ldots$ be a sequence of finite sets,
where $X_j$ has exactly $n_j \ge 2$ elements for each $j$.  Also let
$X = \prod_{j = 1}^\infty X_j$ be their Cartesian product, as in
Section \ref{abstract cantor sets}.  Put $t_0 = 1$, and let $t_l > 0$
be defined by
\begin{equation}
\label{1/t_l = prod_{j = 1}^l n_j}
        1/t_l = \prod_{j = 1}^l n_j
\end{equation}
when $l \ge 1$.  Thus $\{t_l\}_{l = 0}^\infty$ is a strictly
decreasing sequence of positive real numbers that converges to $0$,
which leads to an ultrametric $d(x, y)$ on $X$, as in (\ref{d(x, y) =
  t_{l(x, y)}}).  Remember that the closed ball in $X$ centered at a
point $x \in X$ and with radius $t_k$ for some nonnegative integer $k$
is of the form $B_k(x)$ as in (\ref{B_k(x) = {y in X : y_j = x_j for
    each j le k}}).  The diameter of this ball is also equal to $t_k$.
The radius of any closed ball in $X$ with respect to $d(x, y)$ in
(\ref{d(x, y) = t_{l(x, y)}}) can be taken to be $t_k$ for some
nonnegative integer $k$, since these are the only positive values of
$d(x, y)$.

        By construction, every closed ball $B$ in $X$ of radius $t_k$
for some nonnegative integer $k$ is the union of
\begin{equation}
\label{prod_{j = k + 1}^l n_j}
        \prod_{j = k + 1}^l n_j
\end{equation}
pairwise-disjoint closed balls in $X$ of radius $t_l$, for every
integer $l > k$.  This implies that (\ref{H^1_delta(E) =
  H^1_{con}(E)}) also holds for every $E \subseteq X$ in this
situation, and hence that (\ref{H^1(E) = H^1_{con}(E)}) holds for
every $E \subseteq X$ too.  In particular,
\begin{equation}
\label{H^1(B) le diam B}
        H^1(B) \le \diam B.
\end{equation}
As before, one can show that $H^1(B) \ge \diam B$, by considering
coverings of $B$ by finitely many closed balls, and subdividing the
balls to get finitely many smaller balls of the same radius.  This
implies that
\begin{equation}
\label{H^1(B) = diam B}
        H^1(B) = \diam B,
\end{equation}
which corresponds exactly to (\ref{mu(B_k(x)) = prod_{j = 1}^k
  mu_j({x_j})}), in the case where $\mu_j$ is uniformly distributed on
$X_j$ for each $j$.

\section{Carath\'eodory's construction}
\label{caratheodory's construction}

        Let $(M, d(x, y))$ be a metric space, let $\mathcal{F}$
be a collection of subsets of $M$, and let $\zeta$ be a nonnegative
extended real-valued function on $\mathcal{F}$.  Also let 
$0 < \delta \le \infty$ be given, and put
\begin{eqnarray}
\label{H_delta(E) = ...}
 H_\delta(E) & = & \inf\bigg\{\sum_j \zeta(A_j) : E \subseteq \bigcup_j A_j, \, 
                A_j \in \mathcal{F} \hbox{ for each } j, \\
       & &   \qquad\qquad\qquad\qquad\hbox{and } \diam A_j < \delta 
                                       \hbox{ for each } j\bigg\} \nonumber
\end{eqnarray}
for each $E \subseteq M$.  More precisely, the infimum is taken over
all collections $\{A_j\}_j$ of finitely or countably many elements of
$\mathcal{F}$ such that $E \subseteq \bigcup_j A_j$ and $\diam A_j <
\delta$ for each $j$, if there are any.  If there are no such
coverings of $E$, then we put $H_\delta(E) = +\infty$.  If $E =
\emptyset$, then we interpret $H_\delta(E)$ as being equal to $0$,
using the empty covering of $E$, and interpreting an empty sum as
being $0$.

        As in \cite{mat}, one can avoid these problems with two very
mild additional hypotheses.  The first is that for each $\delta > 0$,
there be a collection $\{A_j\}_j$ of finitely or countably many
elements of $\mathcal{F}$ such that $\bigcup_j A_j = M$ and $\diam A_j
< \delta$ for every $j$.  If $\mathcal{F}$ is the collection of all
subsets of $M$, then this is equivalent to asking that $M$ be
separable.  This condition ensures that the coverings used in the
definition of $H_\delta(E)$ always exist.  The second additional
hypothesis is that for each $\delta > 0$, there be an $A \in
\mathcal{F}$ such that $\diam A < \delta$ and $\zeta(A) < \delta$.
This implies that $H_\delta(\emptyset) = 0$, without using the empty
covering.  In particular, this holds when $\emptyset \in \mathcal{F}$
and $\zeta(\emptyset) = 0$, so that one can cover the empty set by itself.

        Observe that
\begin{equation}
\label{H_delta(E) le H_delta(widetilde{E})}
        H_\delta(E) \le H_\delta(\widetilde{E})
\end{equation}
for every $\delta > 0$ when $E \subseteq \widetilde{E} \subseteq M$.
This simply uses the fact that every covering of $\widetilde{E}$ as in
(\ref{H_delta(E) = ...}) is also a covering of $E$, so that
$H_\delta(E)$ is the infimum of a larger collection of sums than for
$H_\delta(\widetilde{E})$.  In many situations, $\zeta$ may enjoy the
monotonicity property
\begin{equation}
\label{zeta(A) le zeta(B)}
        \zeta(A) \le \zeta(B)
\end{equation}
for every $A, B \in \mathcal{F}$ with $A \subseteq B$, but this is not
needed to get (\ref{H_delta(E) le H_delta(widetilde{E})}).  One can
also show that $H_\delta$ is countably subadditive for each $\delta >
0$, by standard arguments, so that $H_\delta$ is an outer measure on
the $\sigma$-algebra of all subsets of $M$.  As before, if $0 < \delta
< \eta \le \infty$, then
\begin{equation}
\label{H_eta(E) le H_delta(E)}
        H_\eta(E) \le H_\delta(E)
\end{equation}
for every $E \subseteq M$, because $H_\eta(E)$ is the infimum of a
larger class of sums than for $H_\delta(E)$.  If $E_1, E_2 \subseteq M$
satisfy $d(x, y) \ge \delta$ for every $x \in E_1$ and $y \in E_2$,
then it is easy to see that
\begin{equation}
\label{H_delta(E_1) + H_delta(E_2) le H_delta(E_1 cup E_2)}
        H_\delta(E_1) + H_\delta(E_2) \le H_\delta(E_1 \cup E_2),
\end{equation}
for the same reasons as in Section \ref{restricting the diameters}.
The opposite inequality holds for any $E_1, E_2 \subseteq M$, so that
equality holds in (\ref{H_delta(E_1) + H_delta(E_2) le H_delta(E_1 cup
  E_2)}) under these conditions.

        Put
\begin{equation}
\label{H(E) = sup_{delta > 0} H_delta(E)}
        H(E) = \sup_{\delta > 0} H_\delta(E)
\end{equation}
for each $E \subseteq M$, which can also be interpreted as a limit as
$\delta \to 0$, because of (\ref{H_eta(E) le H_delta(E)}).  As usual,
$H(\emptyset) = 0$, and
\begin{equation}
\label{H(E) le H(widetilde{E})}
        H(E) \le H(\widetilde{E})
\end{equation}
when $E \subseteq \widetilde{E} \subseteq M$, by (\ref{H_delta(E) le
  H_delta(widetilde{E})}).  Similarly, the countable subadditivity of
$H_\delta$ for each $\delta > 0$ implies the same property for $H$,
and hence that $H$ is an outer measure on the $\sigma$-algebra of all
subsets of $M$.  If $E_1, E_2 \subseteq M$ satisfy $d(x, y) \ge \eta$
for some $\eta > 0$ and every $x \in E_1$ and $y \in E_2$, then
(\ref{H_delta(E_1) + H_delta(E_2) le H_delta(E_1 cup E_2)}) holds when
$0 < \delta \le \eta$, and hence
\begin{equation}
\label{H(E_1) + H(E_2) le H(E_1 cup E_2)}
        H(E_1) + H(E_2) \le H(E_1 \cup E_2).
\end{equation}
The opposite inequality holds automatically, and it follows that $H$
is countably additive on a suitable $\sigma$-algebra of measurable
sets that includes the Borel sets, by Carath\'eodory's criterion.

        Suppose that $E \subseteq M$ satisfies $H(E) < \infty$, and let
$n \in {\bf Z}_+$ be given.  As in Section \ref{hausdorff measures, section},
there is a collection $\{A_{j, n}\}_{j \in I_n}$ of finitely or
countably many elements of $\mathcal{F}$ such that $E \subseteq
\bigcup_{j \in I_n} A_{j, n}$, $\diam A_{j, n} \le 1/n$ for every $j
\in I_n$, and
\begin{equation}
\label{sum_{j in I_n} zeta(A_{j, n}) < H(E) + 1/n}
        \sum_{j \in I_n} \zeta(A_{j, n}) < H(E) + 1/n.
\end{equation}
If we put
\begin{equation}
\label{widetilde{E} = bigcap_{n = 1}^infty (bigcup_{j in I_n} A_{j, n}), 2}
 \widetilde{E} = \bigcap_{n = 1}^\infty \Big(\bigcup_{j \in I_n} A_{j, n}\Big),
\end{equation}
then $E \subseteq \widetilde{E}$ and $\widetilde{E} \subseteq
\bigcup_{j \in I_n} A_{j, n}$ for each $n$.  This implies that $H(E) =
H(\widetilde{E})$, since the first inclusion implies that (\ref{H(E)
  le H(widetilde{E})}) holds, and the opposite inequality can be
derived from the second inclusion and the definition of
$H(\widetilde{E})$.  If every element of $\mathcal{F}$ is a Borel set,
then $\widetilde{E}$ is a Borel set too.

        Alternatively, one might define $H_\delta'(E)$ in the same
way as $H_\delta(E)$, except for replacing the requirement that $\diam
A_j < \delta$ for each $j$ in (\ref{H_delta(E) = ...})  with the
weaker condition that $\diam A_j \le \delta$ for each $j$.  If $\delta
= \infty$, then this condition on $\diam A_j$ is vacuous, so that
$H'_\infty$ is analogous to Hausdorff content.  It is easy to see that
$H_\delta'(E)$ satisfies the analogues of (\ref{H_delta(E) le
  H_delta(widetilde{E})}) and (\ref{H_eta(E) le H_delta(E)}) for each
$\delta > 0$, and that $H_\delta'(E)$ is countably subadditive, for
the same reasons as before.  Thus $H_\delta'$ is also an outer measure
on the $\sigma$-algebra of all subsets of $M$ for each $\delta > 0$.
If $E_1, E_2 \subseteq M$ satisfy $d(x, y) > \delta$ for some $\delta
> 0$ and every $x \in E_1$ and $y \in E_2$, then one can check that
(\ref{H_delta(E_1) + H_delta(E_2) le H_delta(E_1 cup E_2)}) holds, as
before.  Of course,
\begin{equation}
\label{H_delta'(E) le H_delta(E)}
        H_\delta'(E) \le H_\delta(E)
\end{equation}
for every $\delta > 0$ and $E \subseteq M$, because $H_\delta'(E)$
is the infimum of a larger collection of sums than for $H_\delta(E)$.
Similarly,
\begin{equation}
\label{H_eta(E) le H_delta'(E)}
        H_\eta(E) \le H_\delta'(E)
\end{equation}
for every $E \subseteq M$ when $0 < \delta < \eta \le +\infty$,
because $H_\eta(E)$ is the infimum of a larger collection of sums than
for $H_\delta'(E)$.  It follows that the supremum of $H_\delta'(E)$
over $\delta > 0$ is the same as the supremum of $H_\delta(E)$ over
$\delta > 0$, which is equal to $H(E)$.

        If $\mathcal{F}$ is the collection of all subsets of $M$ and 
\begin{equation}
\label{zeta(A) = (diam A)^alpha}
        \zeta(A) = (\diam A)^\alpha
\end{equation}
for some $\alpha \ge 0$ and every $A \subseteq M$, then $H_\delta(E)$
is the same as $H^\alpha_\delta(E)$ in Section \ref{restricting the
  diameters}, and $H(E)$ is the same as the $\alpha$-dimensional
Hausdorff measure of $E$.  We have also seen that we can take
$\mathcal{F}$ to be the collection of all closed subsets of $M$, or
the collection of all open subsets of $M$, when $\zeta(A)$ is as in
(\ref{zeta(A) = (diam A)^alpha}), and get the same results for
$H_\delta(E)$ and $H(E)$.  Similarly, if $\zeta(A)$ is as in
(\ref{zeta(A) = (diam A)^alpha}), then we can take $\mathcal{F}$ to be
the collection of all closed subsets of $M$ and get the same result
for $H_\delta'(E)$ as when $\mathcal{F}$ is the collection of all
subsets of $M$, for each $\delta > 0$.  However, the analogous
argument for open sets does not work for $H_\delta'(E)$ when $0 <
\delta < \infty$, because approximations of a set $A \subseteq M$ by
open sets that contain $A$ may have diameter greater than $\delta$
when $\diam A = \delta$.

        Let $\mathcal{F}$ and $\zeta$ be given as before, and let
$\mathcal{A}$ be a $\sigma$-algebra of subsets of $M$ that contains
$\mathcal{F}$.  Suppose that $\mu$ is an outer measure on $\mathcal{A}$
such that
\begin{equation}
\label{mu(A) le C zeta(A)}
        \mu(A) \le C \, \zeta(A)
\end{equation}
for some nonnegative real number $C$ and every $A \in \mathcal{F}$.
If $E \in \mathcal{A}$, and if $\{A_j\}_j$ are finitely or countably
many elements of $\mathcal{F}$ such that $E \subseteq \bigcup_j A_j$,
then
\begin{equation}
\label{mu(E) le sum_j mu(A_j) le C sum_j zeta(A_j)}
        \mu(E) \le \sum_j \mu(A_j) \le C \, \sum_j \zeta(A_j).
\end{equation}
This implies that
\begin{equation}
\label{mu(E) le C H_infty'(E)}
        \mu(E) \le C \, H_\infty'(E)
\end{equation}
for every $E \in \mathcal{A}$, where $H_\infty'$ is the outer measure
on $M$ corresponding to $\delta = \infty$ discussed earlier.

        Let $\mathcal{F}$ and $\zeta$ be given again, and let 
$\widetilde{d}(x, y)$ be another metric on $M$.  Also let
$\widetilde{H}_\delta(E)$, $\widetilde{H}_\delta'(E)$, and $\widetilde{H}(E)$ 
be the analogues of $H_\delta(E)$, $H_\delta'(E)$, and $H(E)$, using
$\widetilde{d}(x, y)$ to define diameters of subsets of $M$ instead
of $d(x, y)$.  If the identity mapping on $M$ is uniformly continuous
as a mapping from $M$ equipped with $d(x, y)$ to $M$ equipped with
$\widetilde{d}(x, y)$, then for each $\epsilon > 0$ there is a
$\delta > 0$ such that
\begin{equation}
\label{widetilde{H}_epsilon(E) le H_delta(E)}
        \widetilde{H}_\epsilon(E) \le H_\delta(E)
\end{equation}
for every $E \subseteq M$, and similarly for
$\widetilde{H}_\epsilon'(E)$ and $H_\delta'(E)$.  In the limit as
$\epsilon \to 0$, we get that
\begin{equation}
\label{widetilde{H}(E) le H(E)}
        \widetilde{H}(E) \le H(E)
\end{equation}
for every $E \subseteq M$ under these conditions.  If the identity mapping
on $M$ is uniformly continuous as a mapping from $M$ equipped with
$\widetilde{d}(x, y)$ to $M$ equipped with $d(x, y)$, then
\begin{equation}
\label{H(E) le widetilde{H}(E)}
        H(E) \le \widetilde{H}(E)
\end{equation}
for every $E \subseteq M$, for the same reasons.  This implies that
\begin{equation}
\label{widetilde{H}(E) = H(E)}
        \widetilde{H}(E) = H(E)
\end{equation}
for every $E \subseteq M$ when $d(x, y)$ and $\widetilde{d}(x, y)$
determine the same uniform structure on $M$.  Of course, it is
important here that we are using the same function $\zeta(A)$ for both
metrics.

\section{Snowflakes and quasi-metrics}
\label{snowflakes, quasi-metrics}

        Let $(M, d(x, y))$ be a metric space, and suppose that
$d(x, y)^a$ is also a metric on $M$ for some $a > 0$.  As in Section
\ref{snowflake metrics, quasi-metrics}, this holds when $0 < a \le 1$
and $d(x, y)$ is any metric on $M$, and for all $a > 0$ when $d(x, y)$
is an ultrametric on $M$.  It is easy to see that the diameter of $A
\subseteq M$ with respect to $d(x, y)^a$ is equal to $(\diam A)^a$,
where $\diam A$ is the diameter of $A$ with respect to $d(x, y)$.
This implies that the $\alpha$-dimensional Hausdorff content of $E
\subseteq M$ with respect to $d(x, y)^a$ is equal to the $(\alpha \,
a)$-dimensional Hausdorff content of $E$ with respect to $d(x, y)$,
for each $\alpha \ge 0$.  Similarly, the analogue of
$H^\alpha_\delta(E)$ with respect to $d(x, y)^a$ corresponds to
$H^{\alpha'}_{\delta'}(E)$ with respect to $d(x, y)$, where $\alpha' =
\alpha \, a$ and $\delta' = \delta^{1/a}$.  It follows that the
$\alpha$-dimensional Hausdorff measure of $E$ with respect to $d(x,
y)^a$ is the same as the $(\alpha \, a)$-dimensional Hausdorff measure
of $E$ with respect to $d(x, y)$.  In particular, the Hausdorff
dimension of $E$ with respect to $d(x, y)^a$ is equal to the Hausdorff
dimension of $E$ with respect to $d(x, y)$ divided by $a$.

        As in Section \ref{snowflake metrics, quasi-metrics}, $d(x, y)^a$ 
is a quasi-metric on $M$ for every $a > 0$ when $d(x, y)$ is a metric
on $M$, or even a quasi-metric on $M$.  One could define diameters,
Hausdorff measures, and so on with respect to quasi-metrics, in which
case the remarks in the previous paragraph would hold for all $a > 0$.
However, there are some technical problems with this, related to the
continuity properties of $d(x, y)$.  If $d(x, y)$ is a metric on $M$,
then the diameter of a set $A \subseteq M$ is the same as the diameter
of the closure of $A$, and $A$ is contained open subsets of $M$ with
approximately the same diameter, as in Section \ref{diameters}.  Of
course, this also works for quasi-metrics on $M$ of the form $d_0(x,
y)^a$ for some metric $d_0(x, y)$ on $M$ and $a > 0$, by reducing to
the corresponding statements for $d_0(x, y)$.

        If $d(x, y)$ is a quasi-metric on $M$ of the form $d_0(x, y)^a$
for some metric $d_0(x, y)$ on $M$ and $a > 0$, then one might as well
use Hausdorff measures with respect to $d_0(x, y)$ on $M$ to get
Hausdorff measures with respect to $d(x, y)$, with suitable
adjustments to the dimensions, as before.  Alternatively, let $d(x,
y)$ be a quasi-metric on $M$, and suppose that $d_1(x, y)$ is a metric
on $M$ that defines the same uniform structure on $M$.  This is
equivalent to saying that the identity mapping on $M$ is uniformly
continuous as a mapping from $M$ equipped with $d(x, y)$ to $M$
equipped with $d_1(x, y)$, and as a mapping from $M$ equipped with
$d_1(x, y)$ to $M$ equipped with $d(x, y)$, where uniform continuity
can be characterized in the usual way in terms of $\epsilon$'s and
$\delta$'s.  One can then define Hausdorff measures on $M$ with
respect to $d(x, y)$ using the construction described in the previous
section, where the metric $d(x, y)$ in the previous section is taken
to be $d_1(x, y)$, and where $\zeta(A)$ is defined in terms of the
diameter of $A$ with respect to $d(x, y)$.  If $d_2(x, y)$ is another
metric on $M$ that defines the same uniform structure on $M$ as $d(x, y)$,
then $d_1(x, y)$ and $d_2(x, y)$ also determine the same uniform structure
on $M$, and they lead to the same measures on $M$ as before.

\section{Other Hausdorff measures}
\label{other hausdorff measures}

        Let $(M, d(x, y))$ be a metric space, and let $\mathcal{F}$
be the collection of all subsets of $M$.  Also let $h$ be a nonnegative
real-valued function on the set $[0, +\infty)$ of nonnegative real numbers,
and put
\begin{equation}
\label{zeta(A) = h(diam A)}
        \zeta(A) = h(\diam A)
\end{equation}
for every bounded set $A \subseteq M$.  Let us interpret this as being
equal to $0$ when $A = \emptyset$, which is automatic when $h(0) = 0$.
One can include unbounded sets $A \subseteq M$ as well, with the
convention that
\begin{equation}
\label{h(+infty) = sup_{t ge 0} h(t)}
        h(+\infty) = \sup_{t \ge 0} h(t),
\end{equation}
which may be infinite.  This leads to outer measures $H_\delta$ and
$H_\delta'$ on $M$ for each $\delta > 0$ as in Section
\ref{caratheodory's construction}, and to an outer measure $H$ on $M$,
which is the Hausdorff measure\index{Hausdorff measure} associated to
$h$. 

        Of course, this reduces to the previous situation when $h(t) =
t^\alpha$ for some $\alpha \ge 0$.  As usual, one can get the same results
for $H_\delta$, $H_\delta'$, and $H$ by taking $\mathcal{F}$ to be the
collection of all closed subsets of $M$, because of (\ref{diam
  overline{A} = diam A}).  If $h(t)$ is continuous from the right at
each $t \ge 0$, then one can also get the same results for $H_\delta$
and hence $H$ by taking $\mathcal{F}$ to be the collection of all open
subsets of $M$.  If $M$ is the real line with the standard metric,
then one can get the same results for $H_\delta$, $H_\delta'$, and $H$
using the collection of all closed intervals in ${\bf R}$, as in
Section \ref{some special cases}.  Similarly, if $d(x, y)$ is an
ultrametric on any set $M$, then one can get the same results for
$H_\delta$, $H_\delta'$, and $H$ using the collection of all closed
balls in $M$, as in Section \ref{some special cases}.

        Let $X_1, X_2, X_3, \ldots$ be a sequence of finite sets,
where $X_j$ has exactly $n_j \ge 2$ elements for each $j$, and let $X$
be their Cartesian product, as in Section \ref{abstract cantor sets}.
Also let $\{t_l\}_{l = 0}^\infty$ be a strictly decreasing sequence of
positive real numbers that converges to $0$, and let $d(x, y)$ be the
corresponding ultrametric on $X$, as in (\ref{d(x, y) = t_{l(x, y)}}).
Put $\widetilde{t}_0 = 1$, and let $\widetilde{t}_l$ be defined for $l
\ge 1$ by
\begin{equation}
\label{1/widetilde{t}_l = prod_{j = 1}^l n_j}
        1/\widetilde{t}_l = \prod_{j = 1}^l n_j,
\end{equation}
so that $\{\widetilde{t}_l\}_{l = 0}^\infty$ is also a strictly
decreasing sequence of positive real numbers that converges to $0$.
If $\widetilde{d}(x, y)$ is the ultrametric on $X$ that corresponds to
$\{\widetilde{t}_l\}_{l = 0}^\infty$ as in (\ref{d(x, y) = t_{l(x,
    y)}}), then $\widetilde{d}(x, y)$ is the same as the ultrametric
considered in Section \ref{some special cases}.  Note that $d(x, y)$
and $\widetilde{d}(x, y)$ determine the same uniform structure on $X$.

        Let $h$ be a nonnegative real-valued function on $[0, +\infty)$
such that $h(0) = 0$ and
\begin{equation}
\label{h(t_l) = widetilde{t}_l}
        h(t_l) = \widetilde{t}_l
\end{equation}
for each $l \ge 0$.  If $\diam A$ is the diameter of $A \subseteq M$
with respect to $d(x, y)$, then $h(\diam A)$ is equal to the diameter
of $A$ with respect to $\widetilde{d}(x, y)$.  Let $H(E)$ be the outer
measure on $X$ corresponding to (\ref{zeta(A) = h(diam A)}) and the
collection $\mathcal{F}$ of all subsets of $X$ as in Section
\ref{caratheodory's construction}, and let $\widetilde{H}^1(E)$ be
one-dimensional Hausdorff measure on $X$ with respect to
$\widetilde{d}(x, y)$.  It is easy to see that $H(E) =
\widetilde{H}^1(E)$ for every $E \subseteq M$ under these conditions,
using the remarks at the end of Section \ref{caratheodory's
  construction}.  Remember that $\widetilde{H}^1$ can be analyzed as
in Section \ref{some special cases}.

\section{Product spaces}
\label{product spaces}

        Let $(M_1, d_1(x_1, y_1))$ and $(M_2, d_2(x_2, y_2))$ be metric
spaces, and let $M = M_1 \times M_2$ be their Cartesian product.  It
is easy to see that
\begin{equation}
\label{d(x, y) = max(d_1(x_1, y_1), d_2(x_2, y_2))}
        d(x, y) = \max(d_1(x_1, y_1), d_2(x_2, y_2))
\end{equation}
defines a metric on $M$, where $x = (x_1, x_2)$, $y = (y_1, y_2)$.
This metric has the nice property that the open ball in $M$ centered
at a point $x = (x_1, x_2)$ and with radius $r > 0$ is equal to the
Cartesian product of the open balls in $M_1$, $M_2$ centered at $x_1$,
$x_2$ with radii equal to $r$.  In particular, the topology on $M$
determined by (\ref{d(x, y) = max(d_1(x_1, y_1), d_2(x_2, y_2))}) is
the same as the product topology associated to the topologies on $M_1$
and $M_2$ determined by the metrics $d_1(x_1, y_1)$ and $d_2(x_2,
y_2)$, respectively.  Alternatively,
\begin{equation}
\label{D_p(x, y) = (d_1(x_1, y_1)^p + d_2(x_2, y_2)^p)^{1/p}}
        D_p(x, y) = (d_1(x_1, y_1)^p + d_2(x_2, y_2)^p)^{1/p}
\end{equation}
defines a metric on $M$ when $1 \le p < \infty$, because of the
triangle inequality for $\ell^p$ norms.  This is especially simple
when $p = 1$, and the $p = 2$ case is very natural in the context of
Euclidean geometry.  Observe that
\begin{equation}
\label{d(x, y) le D_p(x, y) le 2^{1/p} d(x, y)}
        d(x, y) \le D_p(x, y) \le 2^{1/p} \, d(x, y)
\end{equation}
for every $x, y \in M$ and $1 \le p < \infty$, which implies that
$D_p(x, y)$ determines the same topology on $M$ as $d(x, y)$.  This
also implies analogous relations between diameters of subsets of $M$
with respect to these metrics, and permits one to compare Hausdorff
measures on $M$ with respect to these metrics.  Another nice property
of (\ref{d(x, y) = max(d_1(x_1, y_1), d_2(x_2, y_2))}) is that it is
an ultrametric on $M$ when $d_1(x_1, y_1)$ and $d_2(x_2, y_2)$ are
ultrametrics on $M_1$ and $M_2$, respectively.

        Let $p_1 : M \to M_1$ and $p_2 : M \to M_2$ be the obvious
coordinate projections, so that $p_1(x) = x_1$ and $p_2(x) = x_2$
for every $x = (x_1, x_2) \in M$.  If $A \subseteq M$, then
\begin{equation}
\label{diam A = max(diam p_1(A), diam p_2(A))}
        \diam A = \max(\diam p_1(A), \diam p_2(A)),
\end{equation}
where $\diam A$ is the diameter of $A$ with respect to (\ref{d(x, y) =
  max(d_1(x_1, y_1), d_2(x_2, y_2))}), and $\diam p_1(A)$, $\diam
p_2(A)$ are the diameters of $p_1(A)$, $p_2(A)$ in $M_1$, $M_2$,
respectively.  It follows that the diameters of $A$ and $p_1(A) \times
p_2(A)$ with respect to (\ref{d(x, y) = max(d_1(x_1, y_1), d_2(x_2,
  y_2))}) on $M$ are the same.  This implies that Hausdorff measures
of a set $E \subseteq M$ with respect to (\ref{d(x, y) = max(d_1(x_1,
  y_1), d_2(x_2, y_2))}) can be defined equivalently in terms of
coverings of $E$ by products of subsets of $M_1$ and $M_2$.  More
precisely, one can restrict one's attention to coverings of $E$ by
products of closed subsets of $M_1$ and $M_2$, because of (\ref{diam
  overline{A} = diam A}).

        Let $h_1$, $h_2$ be monotone increasing nonnegative 
real-valued functions on $[0, +\infty)$, and put
\begin{equation}
\label{h_j(+infty) = sup_{t ge 0} h_j(t)}
        h_j(+\infty) = \sup_{t \ge 0} h_j(t)
\end{equation}
for $j = 1, 2$, which may be infinite.  Suppose that $\mu_1$, $\mu_2$
are nonnegative Borel measures on $M_1$ and $M_2$ such that
\begin{equation}
\label{mu_1(A_1) le C_1 h_1(diam A_1)}
        \mu_1(A_1) \le C_1 \, h_1(\diam A_1)
\end{equation}
and
\begin{equation}
\label{mu_2(A_2) le C_2 h_2(diam A_2)}
        \mu_2(A_2) \le C_2 \, h_2(\diam A_2)
\end{equation}
for some nonnegative real numbers $C_1$, $C_2$ and all Borel sets $A_1
\subseteq M_1$ and $A_2 \subseteq M_2$.  In particular, this ensures
that $M_1$, $M_2$ are $\sigma$-finite with respect to $\mu_1$,
$\mu_2$, so that the product measure $\mu = \mu_1 \times \mu_2$ can be
defined on a suitable $\sigma$-algebra of subsets $M$.  If $M_1$ and
$M_2$ are separable, then $M$ is separable, which implies that open
subsets of $M$ can be expressed as unions of finitely or countably
many products of open subsets of $M_1$ and $M_2$.  In this case, open
subsets of $M$ are measurable with respect to the product measure
construction, and hence Borel subsets of $M$ are measurable too.

        Put
\begin{equation}
\label{h(t) = h_1(t) h_2(t)}
        h(t) = h_1(t) \, h_2(t),
\end{equation}
when $0 \le t < \infty$, which is also a monotone increasing
nonnegative real-valued function on $[0, +\infty)$.  Note that
\begin{equation}
\label{sup_{t ge 0} h(t) = h_1(+infty) h_2(+infty)}
        \sup_{t \ge 0} h(t) = h_1(+\infty) \, h_2(+\infty),
\end{equation}
with the convention that $r \cdot (+\infty) = (+\infty) \cdot r$ is
equal to $+\infty$ when $r > 0$, and to $0$ when $r = 0$.  Thus we
take $h(+\infty)$ to be (\ref{sup_{t ge 0} h(t) = h_1(+infty)
  h_2(+infty)}).  If $A_1 \subseteq M_1$, $A_2 \subseteq M_2$ are
Borel sets, $A \subseteq M$ is measurable with respect to the product
measure construction, and $A \subseteq A_1 \times A_2$, then
\begin{eqnarray}
\label{mu(A) le mu(A_1 times A_2) = mu_1(A_1) mu_2(A_2) le ...}
 \mu(A) \le \mu(A_1 \times A_2) & = & \mu_1(A_1) \, \mu_2(A_2) \\
      & \le & C_1 \, C_2 \, h_1(\diam A_1) \, h_2(\diam A_2) \nonumber \\
      & \le & C_1 \, C_2 \, h(\max(\diam A_1, \diam A_2))     \nonumber
\end{eqnarray}
by (\ref{mu_1(A_1) le C_1 h_1(diam A_1)}) and (\ref{mu_2(A_2) le C_2
  h_2(diam A_2)}).  It follows that
\begin{equation}
\label{mu(A) le C_1 C_2 h(diam A)}
        \mu(A) \le C_1 \, C_2 \, h(\diam A),
\end{equation}
by taking $A_1$, $A_2$ to be the closures of $p_1(A)$, $p_2(A)$ in
$M_1$, $M_2$, respectively, and using (\ref{diam A = max(diam p_1(A),
  diam p_2(A))}).

\chapter{Lipschitz mappings}
\label{lipschitz mappings}

\section{Basic properties}
\label{basic properties}

        Let $(M, d(x, y))$ and $(N, \rho(w, z))$ be metric spaces.
A mapping $f : M \to N$ is said to be 
\emph{Lipschitz}\index{Lipschitz mappings} if there is a nonnegative 
real number $C$ such that
\begin{equation}
\label{rho(f(x), f(y)) le C d(x, y)}
        \rho(f(x), f(y)) \le C \, d(x, y)
\end{equation}
for every $x, y \in M$.  In this case, one might also say that $f$ is
$C$-Lipschitz, or Lipschitz with constant $C$, to indicate the
constant $C$.  Of course, $f$ is Lipschitz with constant $C = 0$ if
and only if $f$ is constant.  Note that the composition of two
Lipshitz mappings with constants $C_1$, $C_2$ is Lipschitz with
constant $C_1 \, C_2$.

        Suppose that $f : M \to N$ is Lipschitz with constant $C$,
and that $A$ is a nonempty bounded subset of $M$.  Under these conditions,
$f(A)$ is a nonempty bounded set in $N$, and
\begin{equation}
\label{diam f(A) le C diam A}
        \diam f(A) \le C \, \diam A,
\end{equation}
where more precisely $\diam A = \diam_M A$ is defined using the metric
on $M$, and $\diam f(A) = \diam_N f(A)$ uses the metric on $N$.  This
also works when $A$ is unbounded, with the convention that the right
side of (\ref{diam f(A) le C diam A}) is infinite when $C > 0$ and
equal to $0$ when $C = 0$.  It follows that
\begin{equation}
\label{H^alpha_{con}(f(E)) le C^alpha H^alpha_{con}(E)}
        H^\alpha_{con}(f(E)) \le C^\alpha \, H^\alpha_{con}(E)
\end{equation}
for every $E \subseteq M$ and $\alpha \ge 0$, where
$H^\alpha_{con}(E)$ is defined using the metric on $M$, and
$H^\alpha_{con}(f(E))$ is defined using the metric on $N$, as before.
If $\alpha = 0$, then $C^\alpha$ should be interpreted as being equal
to $1$ for every $C \ge 0$.

        Similarly,
\begin{equation}
\label{H^alpha_{C delta}(f(E)) le C^alpha H^alpha_delta(E)}
        H^\alpha_{C \, \delta}(f(E)) \le C^\alpha \, H^\alpha_\delta(E)
\end{equation}
for every $E \subseteq M$, $\alpha \ge 0$, and $\delta > 0$, at least
when $C > 0$, so that $C \, \delta > 0$.  This implies that
\begin{equation}
\label{H^alpha(f(E)) le C^alpha H^alpha(E)}
        H^\alpha(f(E)) \le C^\alpha \, H^\alpha(E)
\end{equation}
for every $E \subseteq M$ and $\alpha \ge 0$ when $C > 0$, which also
holds trivially when $C = 0$.  Indeed, if $C = 0$ and $\alpha > 0$,
then $H^\alpha(f(E)) = 0$ automatically.  If $\alpha = 0$, then
$H^\alpha$ reduces to counting measure, and the counting measure of
$f(E)$ is less than or equal to the counting measure of $E$ for any
mapping $f : M \to N$ and $E \subseteq M$.

        A mapping $f : M \to N$ is said to be 
\emph{bilipschitz}\index{bilipschitz mappings} if there is a $C \ge 1$
such that
\begin{equation}
\label{C^{-1} d(x, y) le rho(f(x), f(y)) le C d(x, y)}
        C^{-1} \, d(x, y) \le \rho(f(x), f(y)) \le C \, d(x, y)
\end{equation}
for every $x, y \in M$.  As before, one might say that $f$ is
$C$-bilipschitz, or bilipschitz with constant $C$, to indicate the
constant $C$.  If $f$ is bilipschitz with constant $C$, then
\begin{equation}
\label{C^{-1} diam A le diam f(A) le C diam A}
        C^{-1} \, \diam A \le \diam f(A) \le C \, \diam A
\end{equation}
for every nonempty bounded set $A \subseteq M$.  This implies that
\begin{equation}
\label{... le H^alpha_{con}(f(E)) le C^alpha H^alpha_{con}(E)}
        C^{-\alpha} \, H^\alpha_{con}(E) \le H^\alpha_{con}(f(E)) 
                                       \le C^\alpha \, H^\alpha_{con}(E)
\end{equation}
and
\begin{equation}
\label{C^{-alpha} H^alpha(E) le H^alpha(f(E)) le C^alpha H^alpha(E)}
        C^{-\alpha} \, H^\alpha(E) \le H^\alpha(f(E)) \le C^\alpha \, H^\alpha(E)
\end{equation}
for every $E \subseteq M$ and $\alpha \ge 0$.  Of course, the counting
measure of $E$ is equal to the counting measure of $f(E)$ for every $E
\subseteq M$ when $f : M \to N$ is one-to-one.

\section{Real-valued functions}
\label{real-valued functions}

        Let $(M, d(x, y))$ be a metric space, and let $f$ be a real-valued
function on $M$.  Thus $f$ is Lipschitz with constant $C \ge 0$ with respect
to the standard metric on ${\bf R}$ if and only if
\begin{equation}
\label{|f(x) - f(y)| le C d(x, y)}
        |f(x) - f(y)| \le C \, d(x, y)
\end{equation}
for every $x, y \in X$.  Of course, this implies that
\begin{equation}
\label{f(x) le f(y) + C d(x, y)}
        f(x) \le f(y) + C \, d(x, y)
\end{equation}
for every $x, y \in M$.  Conversely, if $f$ satisfies (\ref{f(x) le
  f(y) + C d(x, y)}) for every $x, y \in M$, then we also have that
\begin{equation}
\label{f(y) le f(x) + C d(x, y)}
        f(y) \le f(x) + C \, d(x, y)
\end{equation}
for every $x, y \in M$, by interchanging the roles of $x$ and $y$.  It
is easy to see that (\ref{|f(x) - f(y)| le C d(x, y)}) is implied by
(\ref{f(x) le f(y) + C d(x, y)}) and (\ref{f(y) le f(x) + C d(x, y)}),
so that (\ref{|f(x) - f(y)| le C d(x, y)}) and (\ref{f(x) le f(y) + C
  d(x, y)}) are equivalent to each other.

        In particular,
\begin{equation}
\label{f_p(x) = d(p, x)}
        f_p(x) = d(p, x)
\end{equation}
satisfies (\ref{f(x) le f(y) + C d(x, y)}) for every $p, x, y \in M$
with $C = 1$, by the triangle inequality.  This shows that
(\ref{f_p(x) = d(p, x)}) is a Lipschitz function on $M$ with constant
$C = 1$ for every $p \in M$.  Now let $A$ be a nonempty subset of $M$,
and put
\begin{equation}
\label{dist(x, A) = inf {d(x, z) : z in A}}
        \dist(x, A) = \inf \{d(x, z) : z \in A\}
\end{equation}
for every $x \in M$.  Observe that
\begin{equation}
\label{dist(x, A) le d(x, z) le d(x, y) + d(y, z)}
        \dist(x, A) \le d(x, z) \le d(x, y) + d(y, z)
\end{equation}
for every $x, y \in M$ and $z \in A$, which implies that
\begin{equation}
\label{dist(x, A) le d(x, y) + dist(y, A)}
        \dist(x, A) \le d(x, y) + \dist(y, A)
\end{equation}
for every $x, y \in M$.  Thus (\ref{dist(x, A) = inf {d(x, z) : z in
    A}}) is also a Lipschitz function on $M$ with constant $C = 1$,
for each nonempty set $A \subseteq M$.

        Suppose now that $d(\cdot, \cdot)$ is an ultrametric on $M$.
In this case, we have that
\begin{equation}
\label{dist(x, A) le d(x, z) le max(d(x, y), d(y, z))}
        \dist(x, A) \le d(x, z) \le \max(d(x, y), d(y, z))
\end{equation}
for every $x, y \in M$ and $z \in A$, which is stronger than
(\ref{dist(x, A) le d(x, z) le d(x, y) + d(y, z)}).  If
\begin{equation}
\label{d(x, y) < dist(x, A)}
        d(x, y) < \dist(x, A),
\end{equation}
then it follows that
\begin{equation}
\label{dist(x, A) le d(y, z)}
        \dist(x, A) \le d(y, z)
\end{equation}
for every $z \in A$, and hence
\begin{equation}
\label{dist(x, A) le dist(y, A)}
        \dist(x, A) \le \dist(y, A).
\end{equation}
Combining (\ref{d(x, y) < dist(x, A)}) and (\ref{dist(x, A) le dist(y,
  A)}), we get that
\begin{equation}
\label{d(x, y) < dist(y, A)}
        d(x, y) < \dist(y, A),
\end{equation}
so that
\begin{equation}
\label{dist(y, A) le dist(x, A)}
        \dist(y, A) \le \dist(x, A),
\end{equation}
by the same argument.  This shows that
\begin{equation}
\label{dist(x, A) = dist(y, A)}
        \dist(x, A) = \dist(y, A)
\end{equation}
when $x, y \in M$ satisfy (\ref{d(x, y) < dist(x, A)}).

        Let $d(x, y)$ be any metric on $M$ again, and let $E$ be a 
connected subset of $M$.  If $p, q \in E$ and $f_p(x)$ is as in
(\ref{f_p(x) = d(p, x)}), then $f_p(E)$ is a connected subset of ${\bf
  R}$ that contains $0$ and $d(p, q)$, and hence contains $[0, d(p,
  q)]$.  This implies that
\begin{equation}
\label{d(p, q) le H^1(f_p(E)) le H^1(E)}
        d(p, q) \le H^1(f_p(E)) \le H^1(E)
\end{equation}
for every $p, q \in E$, so that
\begin{equation}
\label{diam E le H^1(E)}
        \diam E \le H^1(E).
\end{equation}

\section{Some examples}
\label{some examples}

        Let $n_1, n_2, n_3, \ldots$ be a sequence of integers with
$n_j \ge 2$ for each $j$, and put
\begin{equation}
\label{X_j = {0, 1, ldots, n_j - 1}}
        X_j = \{0, 1, \ldots, n_j - 1\}
\end{equation}
for each $j \in {\bf Z}_+$.  Thus $X_j$ has exactly $n_j$ elements for
each $j$, and we let $X = \prod_{j = 1}^\infty X_j$ be their Cartesian
product, as in Section \ref{abstract cantor sets}.  Also put $N_k =
\prod_{j = 1}^k n_j$ for each $k \in {\bf Z}_+$ and $N_0 = 1$, and
$t_l = 1/N_l$ for every $l \ge 0$.  This leads to an ultrametric $d(x,
y)$ on $X$ as in (\ref{d(x, y) = t_{l(x, y)}}), for which the
corresponding one-dimensional Hausdorff measure was discussed in
Section \ref{some special cases}.

        Observe that
\begin{equation}
\label{N_{j - 1}^{-1} - N_j^{-1} = n_j N_j^{-1} - N_j^{-1} = (n_j - 1) N_j^{-1}}
 N_{j - 1}^{-1} - N_j^{-1} = n_j \, N_j^{-1} - N_j^{-1} = (n_j - 1) \, N_j^{-1}
\end{equation}
for each $j \in {\bf Z}_+$, and hence
\begin{equation}
\label{sum_{j = k}^l (n_j - 1) N_j^{-1} = N_{k - 1}^{-1} - N_l^{-1}}
        \sum_{j = k}^l (n_j - 1) \, N_j^{-1} = N_{k - 1}^{-1} - N_l^{-1}
\end{equation}
when $1 \le k \le l$.  Put
\begin{equation}
\label{f_k(x) = sum_{j = 1}^k x_j N_j^{-1}}
        f_k(x) = \sum_{j = 1}^k x_j \, N_j^{-1}
\end{equation}
for each $x \in X$ and $k \in {\bf Z}_+$, and $f_0(x) = 0$.  Thus
$f_k(x)$ is an integer multiple of $N_k^{-1}$ for each $x \in X$ and
$k \ge 0$, and
\begin{equation}
\label{0 le f_k(x) le 1 - N_k^{-1} < 1}
        0 \le f_k(x) \le 1 - N_k^{-1} < 1,
\end{equation}
by (\ref{sum_{j = k}^l (n_j - 1) N_j^{-1} = N_{k - 1}^{-1} -
  N_l^{-1}}).  One can check that every nonnegative integer multiple
of $N_k^{-1}$ strictly less than $1$ can be expressed as $f_k(x)$ for
some $x \in X$, using induction on $k$.

        If $x, y \in X$ satisfy $x_j \le y_j$ for $j = 1, \ldots, k$, then
\begin{equation}
\label{f_k(x) le f_k(y)}
        f_k(x) \le f_k(y).
\end{equation}
If $x_j = y_j$ when $j \le k$ and $k < l$, then
\begin{equation}
\label{f_l(y) le ... le f_k(x) + N_k^{-1} - N_l^{-1}}
        f_l(y) \le f_k(x) + \sum_{j = k + 1}^l (n_j - 1) \, N_j^{-1}
                \le f_k(x) + N_k^{-1} - N_l^{-1},
\end{equation}
by (\ref{sum_{j = k}^l (n_j - 1) N_j^{-1} = N_{k - 1}^{-1} - N_l^{-1}}).
Applying this to $y = x$, we get that
\begin{equation}
\label{f_l(x) le f_k(x) + N_k^{-1} - N_l^{-1}}
        f_l(x) \le f_k(x) + N_k^{-1} - N_l^{-1}
\end{equation}
when $k < l$.  If $x_j = y_j$ when $j \le k$ and $x_{k + 1} < y_{k + 1}$, then
\begin{equation}
\label{f_{k + 1}(x) + N_{k + 1}^{-1} le f_{k + 1}(y)}
        f_{k + 1}(x) + N_{k + 1}^{-1} \le f_{k + 1}(y).
\end{equation}
This implies that
\begin{equation}
\label{f_l(x) + N_l^{-1} le f_{k + 1}(y) le f_l(y)}
        f_l(x) + N_l^{-1} \le f_{k + 1}(y) \le f_l(y)
\end{equation}
for every $l \ge k + 1$, because of (\ref{f_l(x) le f_k(x) + N_k^{-1}
  - N_l^{-1}}) applied to $k + 1$ instead of $k$.  In particular,
\begin{equation}
\label{f_l(x) < f_l(y)}
        f_l(x) < f_l(y)
\end{equation}
for each $l \ge k + 1$ under these conditions.  It follows that
\begin{equation}
\label{f_l(x) ne f_l(y)}
        f_l(x) \ne f_l(y)
\end{equation}
 when $x_j \ne y_j$ for some $j \le l$, by considering the smallest
 such $j$.

        Taking the limit as $l \to \infty$ in 
(\ref{sum_{j = k}^l (n_j - 1) N_j^{-1} = N_{k - 1}^{-1} - N_l^{-1}}),
we get that
\begin{equation}
\label{sum_{j = k}^infty (n_j - 1) N_j^{-1} = N_{k - 1}^{-1}}
        \sum_{j = k}^\infty (n_j - 1) \, N_j^{-1} = N_{k - 1}^{-1}
\end{equation}
for each $k \in {\bf Z}_+$, which is equal to $1$ when $k = 1$.  Put
\begin{equation}
\label{f(x) = sum_{j = 1}^infty x_j N_j^{-1}}
        f(x) = \sum_{j = 1}^\infty x_j \, N_j^{-1}
\end{equation}
for each $x \in X$, where the series converges by comparison with
(\ref{sum_{j = k}^infty (n_j - 1) N_j^{-1} = N_{k - 1}^{-1}}).  Thus
\begin{equation}
\label{0 le f(x) le 1}
        0 \le f(x) \le 1
\end{equation}
for every $x \in X$ and $k \in {\bf Z}_+$, and
\begin{equation}
\label{f(x) le f(y)}
        f(x) \le f(y)
\end{equation}
when $x, y \in X$ satisfy $x_j \le y_j$ for each $j$.  If $x_j = y_j$
when $j \le k$ for some $k \ge 0$, and $x_{k + 1} < y_{k + 1}$, then
we also have (\ref{f(x) le f(y)}), by taking the limit as $l \to
\infty$ in (\ref{f_l(x) + N_l^{-1} le f_{k + 1}(y) le f_l(y)}).  In
this case, the only way that equality can hold in (\ref{f(x) le f(y)}) is if
\begin{equation}
\label{y_{k + 1} = x_{k + 1} + 1, and x_l = n_l - 1, y_l = 0 for l ge k + 2}
        y_{k + 1} = x_{k + 1} + 1, \hbox{ and } x_l = n_l - 1, \, y_l = 0
                                  \hbox{ for each } l \ge k + 2.
\end{equation}
If $x \ne y$, then $x_j \ne y_j$ for some $j$, and we can choose $k
\ge 0$ as small as possible so that $x_{k + 1} < y_{k + 1}$.  It
follows that $f(x) = f(y)$ only when $x = y$, or when there is a $k
\ge 0$ such that $x_j = y_j$ for $j \le k$, and (\ref{y_{k + 1} = x_{k
    + 1} + 1, and x_l = n_l - 1, y_l = 0 for l ge k + 2}) holds.

        Suppose again that $x, y \in X$ satisfy $x_j = y_j$ when $j \le k$ 
for some $k \ge 0$.  Of course, $f_k(x) = f_k(y) \le f(y)$, and hence
\begin{equation}
\label{f_k(x) le f(y) le f_k(x) + N_k^{-1}}
        f_k(x) \le f(y) \le f_k(x) + N_k^{-1},
\end{equation}
by taking the limit as $l \to \infty$ in (\ref{f_l(y) le ... le f_k(x)
  + N_k^{-1} - N_l^{-1}}).  In particular,
\begin{equation}
\label{f_k(x) le f(x) le f_k(x) + N_k^{-1}}
        f_k(x) \le f(x) \le f_k(x) + N_k^{-1},
\end{equation}
which implies that
\begin{equation}
\label{|f(x) - f(y)| le N_k^{-1}}
        |f(x) - f(y)| \le N_k^{-1}
\end{equation}
under these conditions.  This shows that $f$ is Lipschitz with
constant $C = 1$ as a mapping from $X$ into ${\bf R}$, with respect
to the ultrametric $d(x, y)$ on $X$ described at the beginning of
the section, and the standard metric on ${\bf R}$.

        Let $x \in X$ and $k \ge 0$ be given, and let $B_k(x)$ be the
set of $y \in X$ such that $x_j = y_j$ when $j \le k$, as in Section
\ref{abstract cantor sets}.  Thus
\begin{equation}
\label{f(B_k(x)) subseteq [f_k(x), f_k(x) + N_k^{-1}]}
        f(B_k(x)) \subseteq [f_k(x), f_k(x) + N_k^{-1}],
\end{equation}
by (\ref{f_k(x) le f(y) le f_k(x) + N_k^{-1}}).  One can check that
\begin{equation}
\label{f(B_k(x)) = [f_k(x), f_k(x) + N_k^{-1}]}
        f(B_k(x)) = [f_k(x), f_k(x) + N_k^{-1}]
\end{equation}
for every $x \in X$ and $k \ge 0$, by standard arguments.  In particular,
\begin{equation}
\label{f(X) = [0, 1]}
        f(X) = [0, 1],
\end{equation}
which is the same as (\ref{f(B_k(x)) = [f_k(x), f_k(x) + N_k^{-1}]})
when $k = 0$.

        Note that
\begin{equation}
\label{H^1(B_k(x)) = N_k^{-1}}
        H^1(B_k(x)) = N_k^{-1}
\end{equation}
for every $x \in X$ and $k \ge 0$, where $H^1(B_k(x))$ is the
one-dimensional Hausdorff measure of $B_k(x)$ with respect to the
ultrametric $d(x, y)$ on $X$ mentioned earlier.  This follows from the
discussion at the end of Section \ref{some special cases}.  In particular,
\begin{equation}
\label{H^1(X) = 1}
        H^1(X) = 1.
\end{equation}
This is also consistent with the discussion of Hausdorff measure and
Lipschitz mappings in Section \ref{basic properties}, since the
one-dimensional Hausdorff measure of an interval in the real line is
the same as the length of the interval.

\section{Other Lipschitz conditions}
\label{other lipschitz conditions}

        Let $(M, d(x, y))$ and $(N, \rho(w, z))$ be metric spaces,
and let $a$ be a positive real number.  A mapping $f : M \to N$ is
said to be \emph{Lipschitz of order $a$}\index{Lipschitz mappings}
if there is a nonnegative real number $C$ such that
\begin{equation}
\label{rho(f(x), f(y)) le C d(x, y)^a}
        \rho(f(x), f(y)) \le C \, d(x, y)^a
\end{equation}
for every $x, y \in M$.  As before, this condition holds with $C = 0$
if and only if $f$ is a constant mapping.  If $a = 1$, then this
condition is equivalent to the one discussed in Section \ref{basic
  properties}.

        If $a \le 1$, then $d(x, y)^a$ is also a metric on $M$, as in
Section \ref{snowflake metrics, quasi-metrics}.  In this case, the
condition described in the previous paragraph is equivalent to saying
that $f$ is Lipschitz of order $1$ with respect to the metric $d(x,
y)^a$ on $M$, and with the same constant $C$.  Similarly, if $d(x, y)$
is an ultrametric on $M$, then $d(x, y)^a$ is also an ultrametric on
$M$ for every $a > 0$, and the condition in the previous paragraph is
equivalent to saying that $f$ is Lipschitz of order $1$ with respect
to $d(x, y)^a$ on $M$.  Otherwise, $d(x, y)^a$ is a quasi-metric on
$M$ for every $a > 0$, as in Section \ref{snowflake metrics,
  quasi-metrics}.  One can define Lipschitz conditions with respect to
quasi-metrics in the same way as for metrics, so that a mapping $f : M
\to N$ is Lipschitz of order $a > 0$ with respect to $d(x, y)$ on $M$
if and only if it is Lipschitz of order $1$ with respect to $d(x,
y)^a$ on $M$.

        There are always a lot of real-valued Lipschitz functions of order
$1$ on any metric space $(M, d(x, y))$, as in Section \ref{real-valued
functions}.  If $0 < a \le 1$, then $d(x, y)^a$ is also a metric on $M$,
and the same discussion can be applied to get a lot of real-valued
Lipschitz functions of order $1$ on $M$ with respect to $d(x, y)^a$,
which are the same as real-valued Lipschitz functions of order $a$ on
$M$.  Of course, the property of being Lipschitz of order $a$ becomes
stronger on bounded sets as $a$ increases, and bounded Lipschitz
functions of order $1$ are also Lipschitz functions of order $a$ when
$0 < a \le 1$.  If $M$ is the real line with the standard metric, then
the only Lipschitz functions of order $a > 1$ are constant, because
the derivative of such a function must be equal to $0$ at every point.
Equivalently, the only Lipschitz functions of order $1$ on ${\bf R}$
with respect to the quasi-metric $|x - y|^a$ are the constant
functions when $a > 1$.  If $d(x, y)$ is any quasi-metric on a set
$M$, then there is a metric $\widetilde{d}(x, y)$ on $M$ and a
positive real number $a$ such that $d(x, y)$ is comparable to
$\widetilde{d}(x, y)^a$, as shown in \cite{m-s-1} and recalled in
Section \ref{snowflake metrics, quasi-metrics}.  This implies that
there are a lot of real-valued Lipschitz functions of order $1$ on $M$
with respect to $\widetilde{d}(x, y)$, which are Lipschitz of order
$1/a$ with respect to $d(x, y)$.

        Suppose that $f : M \to N$ is Lipschitz of some order $a > 0$
with constant $C$, as in (\ref{rho(f(x), f(y)) le C d(x, y)^a}).  If
$A$ is a nonempty bounded subset of $M$, then $f(A)$ is a nonempty
bounded set in $N$, and
\begin{equation}
\label{diam f(A) le C (diam A)^a}
        \diam f(A) \le C \, (\diam A)^a.
\end{equation}
More precisely, $\diam A$ is the diameter of $A$ with respect to the
metric on $M$, and $\diam f(A)$ is the diameter of $f(A)$ with respect
to the metric on $N$.  This implies that
\begin{equation}
\label{H^alpha(f(E)) le C^alpha H^{alpha a}(E)}
        H^\alpha(f(E)) \le C^\alpha \, H^{\alpha \, a}(E)
\end{equation}
for every $E \subseteq M$ and $\alpha \ge 0$, as in Section \ref{basic
  properties}.

\section{Subadditive functions}
\label{subadditive functions}

        Let $\sigma(t)$ be a monotone increasing real-valued function
on the set $[0, +\infty)$ of nonnegative real numbers such that
  $\sigma(0) = 0$, $\sigma(t) > 0$ when $t > 0$, and
\begin{equation}
\label{lim_{t to 0+} sigma(t) = 0}
            \lim_{t \to 0+} \sigma(t) = 0.
\end{equation}
If $d(x, y)$ is an ultrametric on a set $M$, then $\sigma(d(x, y))$ is
also an ultrametric on $M$, which determines the same topology on $M$
as $d(x, y)$.  Of course, this includes the case where $\sigma(t) =
t^a$ for some $a > 0$, as in Section \ref{snowflake metrics,
  quasi-metrics}.  This is also related to the examples discussed in
Section \ref{abstract cantor sets}.  If, in addition to the conditions
just mentioned, $\sigma(t)$ satisfies
\begin{equation}
\label{sigma(r + t) le sigma(r) + sigma(t)}
        \sigma(r + t) \le \sigma(r) + \sigma(t)
\end{equation}
for every $r, t \ge 0$, then $\sigma(t)$ is said to be
\emph{subadditive}.\index{subadditive functions} Remember that
$\sigma(t) = t^a$ is subadditive when $0 < a \le 1$, as in Section
\ref{snowflake metrics, quasi-metrics}.  If $\sigma(t)$ is subadditive
and $d(x, y)$ is a metric on $M$, then $\sigma(d(x, y))$ is also a
metric on $M$, which determines the same topology on $M$ as $d(x, y)$.

        In both cases, the identity mapping on $M$ is uniformly continuous
as a mapping from $M$ equipped with $d(x, y)$ to $M$ equipped with
$\sigma(d(x, y))$, because of (\ref{lim_{t to 0+} sigma(t) = 0}).
Similarly, the identity mapping on $M$ is uniformly continuous as a
mapping from $M$ equipped with $\sigma(d(x, y))$ to $M$ equipped with
$d(x, y)$.  More precisely, let $\epsilon > 0$ be given, and put
\begin{equation}
\label{delta = sigma(epsilon)}
        \delta = \sigma(\epsilon) > 0.
\end{equation}
Thus $\sigma(t) \ge \delta$ when $t \ge \epsilon$, because $\sigma(t)$
is monotone increasing.  Equivalently, this means that $t < \epsilon$
when $\sigma(t) < \delta$, which is exactly what we wanted.

        If $\sigma(t)$ is subadditive on $[0, +\infty)$, then
\begin{equation}
\label{0 le sigma(r + t) - sigma(r) le sigma(t)}
        0 \le \sigma(r + t) - \sigma(r) \le \sigma(t)
\end{equation}
for every $r, t \ge 0$, since $\sigma(\cdot)$ is also supposed to be
monotone increasing on $[0, +\infty)$.  This implies that
  $\sigma(\cdot)$ is uniformly continuous on $[0, +\infty)$, using
(\ref{lim_{t to 0+} sigma(t) = 0}).  Alternatively, it follows from
(\ref{0 le sigma(r + t) - sigma(r) le sigma(t)}) that $\sigma$ is Lipschitz
of order $1$ with constant $C = 1$ as a mapping from $[0, +\infty)$
equipped with the metric $\sigma(|x - y|)$ into the real line with the
standard metric.

        Put
\begin{equation}
\label{sigma(t-) = lim_{r to t-} sigma(r) = sup {sigma(r) : 0 le r < t}}
 \sigma(t-) = \lim_{r \to t-} \sigma(r) = \sup \{\sigma(r) : 0 \le r < t\}
\end{equation}
for each positive real number $t$, so that $\sigma(t-) \le \sigma(t)$
for each $t > 0$, and $\sigma(t-)$ is monotone increasing in $t$.  If
the diameter of a set $A \subseteq M$ with respect to $d(x, y)$ is
equal to $t$, $0 < t < \infty$, then the diameter $T$ of $A$ with
respect to $\sigma(d(x, y))$ satisfies
\begin{equation}
\label{sigma(t-) le T le sigma(t)}
         \sigma(t-) \le T \le \sigma(t).
\end{equation}
In particular,
\begin{equation}
\label{T = sigma(t)}
        T = \sigma(t)
\end{equation}
when $\sigma(t-) = \sigma(t)$, which holds automatically when $\sigma$
is subadditive, as in the previous paragraph.  Of course, if the
diameter of $A$ with respect to $d(x, y)$ is equal to $0$, then the
diameter of $A$ with respect to $\sigma(d(x, y))$ is equal to $0$ too.
If $A$ is unbounded with respect to $d(x, y)$, then the diameter of
$A$ with respect to $\sigma(d(x, y))$ is equal to
\begin{equation}
\label{sigma(+infty) = sup_{r ge 0} sigma(r)}
        \sigma(+\infty) = \sup_{r \ge 0} \sigma(r),
\end{equation}
which is either a positive real number or $+\infty$.

\section{Moduli of continuity}
\label{moduli of continuity}

        Let $(M, d(x, y))$ and $(N, \rho(w, z))$ be metric spaces,
and let $\sigma(t)$ be a monotone increasing nonnegative real-valued
function on $[0, +\infty)$ such that $\sigma(0) = 0$ and $\sigma(t)$
is continuous at $0$.  Suppose that $f : M \to N$ satisfies
\begin{equation}
\label{rho(f(x), f(y)) le sigma(d(x, y))}
        \rho(f(x), f(y)) \le \sigma(d(x, y))
\end{equation}
for every $x, y \in M$, which implies that $f$ is uniformly continuous
in particular.  This includes the Lipschitz condition (\ref{rho(f(x),
  f(y)) le C d(x, y)^a}) as a special case, with $\sigma(t) = C \,
t^\alpha$.  If $\sigma(d(x, y))$ is a metric on $M$, as in the
previous section, then (\ref{rho(f(x), f(y)) le sigma(d(x, y))}) is
the same as saying that $f$ is Lipschitz of order $1$ with constant $C
= 1$ as a mapping from $M$ equipped with the metric $\sigma(d(x, y))$
into $N$ equipped with the metric $\rho(w, z)$.

        If $A$ is a nonempty bounded set in $M$, and $f : M \to N$
satisfies (\ref{rho(f(x), f(y)) le sigma(d(x, y))}), then
\begin{equation}
\label{diam f(A) le sigma(diam A)}
        \diam f(A) \le \sigma(\diam A)
\end{equation}
for every nonempty bounded set $A \subseteq M$.  Here $\diam A$ is the
diameter of $A$ with respect to $d(x, y)$ on $M$, and $\diam f(A)$ is
the diameter of $f(A)$ with respect to $\rho(w, z)$ on $N$, as usual.
This also works when $A$ is unbounded, with $\sigma(+\infty)$ defined
as in (\ref{sigma(+infty) = sup_{r ge 0} sigma(r)}).  Using (\ref{diam
  f(A) le sigma(diam A)}), one can estimate Hausdorff measures of
$f(A)$ in terms of Hausdorff measures of $A$, where the Hausdorff
measures are defined in terms of functions of diameters of sets, as in
Section \ref{caratheodory's construction}.

        If $f$ is any mapping from $M$ into $N$, then put
\begin{equation}
\label{sigma_f(t) = sup {rho(f(x), f(y)) : x, y in M, d(x, y) le t}}
        \sigma_f(t) = \sup \{\rho(f(x), f(y)) : x, y \in M, \, d(x, y) \le t\}
\end{equation}
for each nonnegative real number $t$, where the supremum may be equal
to $+\infty$.  Thus $\sigma_f(t) \ge 0$ for every $t \ge 0$,
$\sigma_f(0) = 0$, $\sigma_f(t)$ is monotone increasing, and
(\ref{rho(f(x), f(y)) le sigma(d(x, y))}) holds with $\sigma(t) =
\sigma_f(t)$ for every $x, y \in M$, by construction.  Note that $f$
is uniformly continuous if and only if $\sigma_f(t) < +\infty$ when
$t$ is sufficiently small, and $\lim_{t \to 0+} \sigma_f(t) = 0$.  The
finiteness of $\sigma_f(t)$ for every $t > 0$ is another matter, and
is trivial when $f(M)$ is bounded in $N$.  However, in order to
estimate Hausdorff measures, it suffices to have a condition like
(\ref{diam f(A) le sigma(diam A)}) when the diameter of $A$ is small.

        Suppose that $f : {\bf R} \to {\bf R}$ satisfies
(\ref{rho(f(x), f(y)) le sigma(d(x, y))}), where $d(x, y)$ and
$\rho(w, z)$ are both equal to the standard metric on ${\bf R}$, and
$\lim_{t \to 0+} \sigma(t)/t = 0$.  This implies that the derivative
of $f$ is equal to $0$ everywhere on ${\bf R}$, and hence that $f$ is
constant on ${\bf R}$.  This includes the case where $f$ is Lipschitz
of order $\alpha > 1$, as in Section \ref{other lipschitz conditions}.
One can check that the analogous statement also holds when $\liminf_{t
  \to 0+} \sigma(t)/ t = 0$.  If $M$ and $N$ are arbitrary metric
spaces, $f : M \to N$ satisfies (\ref{rho(f(x), f(y)) le sigma(d(x,
  y))}), and $\sigma(t) = 0$ for some $t > 0$, then $f$ is locally
constant on $M$, and in particular $f$ is constant on $M$ when $M$ is 
connected.

\section{Isometries and similarities}
\label{isometries, similarities}

        Let $(M, d(x, y))$ and $(N, \rho(w, z))$ be metric spaces.
A mapping $f : M \to N$ is said to be an \emph{isometry}\index{isometries}
if
\begin{equation}
\label{rho(f(x), f(y)) = d(x, y)}
        \rho(f(x), f(y)) = d(x, y)
\end{equation}
for every $x, y \in M$.  Equivalently, $f$ is an isometry if it is a
bilipschitz mapping with constant $C = 1$.  Let us say that $f : M \to
N$ is a \emph{similarity}\index{similarities} if there is a positive
real number $\lambda$ such that
\begin{equation}
\label{rho(f(x), f(y)) = lambda d(x, y)}
        \rho(f(x), f(y)) = \lambda \, d(x, y)
\end{equation}
for every $x, y \in M$.  In this case, it is easy to see that
\begin{equation}
\label{diam f(A) = lambda diam A}
        \diam f(A) = \lambda \, \diam A
\end{equation}
for every nonempty bounded set $A \subseteq M$, and hence that
\begin{equation}
\label{H^alpha(f(E)) = lambda^alpha H^alpha(E)}
        H^\alpha(f(E)) = \lambda^\alpha \, H^\alpha(E)
\end{equation}
for every $E \subseteq M$ and $\alpha \ge 0$.

        Remember that a mapping $f : M \to N$ is said to be \emph{bounded}
if $f(M)$ is a bounded set in $N$.  The space of bounded continuous mappings
from $M$ into $N$ is denoted $C_b(M, N)$, and the supremum metric on
$C_b(M, N)$ is defined by
\begin{equation}
\label{sup {rho(f(x), g(x)) : x in M}}
        \sup \{\rho(f(x), g(x)) : x \in M\}.
\end{equation}
Note that the collection of $f \in C_b(M, N)$ such that $f(M)$ is
dense in $N$ is a closed set in $C_b(M, N)$ with respect to the
supremum metric.  If $M$ is compact, then it follows that the
collection of $f \in C_b(M, N)$ such that $f(M) = N$ is a closed
set in $C_b(M, N)$ with respect to the supremum metric.

        Let $\mathcal{I}(M, N)$ be the collection of isometric embeddings
of $M$ into $N$.  If $M$ is bounded, then $\mathcal{I}(M, N) \subseteq 
C_b(M, N)$, and $\mathcal{I}(M, N)$ is a closed set in $C_b(M, N)$
with respect to the supremum metric.  If $M$ is complete and $f : M \to N$
is an isometry, then $f(M)$ is a closed set in $N$.  In particular,
$f(M) = N$ when $M$ is complete, $f : M \to N$ is an isometry, and $f(M)$
is dense in $N$.  If $M$ and $N$ are compact, then $\mathcal{I}(M, N)$
is a compact set in $C_b(M, N)$ with respect to the supremum metric,
by standard Arzela--Ascoli arguments.

        Let $\mathcal{I}(M)$ be the collection of isometric mappings
of $M$ onto itself, which is a group with respect to composition.  If
$M$ is bounded, then the restriction of the supremum metric to
$\mathcal{I}(M)$ is invariant under left and right translations, and
one can check that $\mathcal{I}(M)$ is a topological group with
respect to the topology determined by the supremum metric.  If $M$ is
complete, then $\mathcal{I}(M)$ is the same as the collection of
isometric mappings $f$ from $M$ into itself such that $f(M)$ is dense
in $M$, as in the previous paragraph.  If $M$ is bounded and complete,
then it follows that $\mathcal{I}(M)$ is a closed subset of
$\mathcal{I}(M, M)$ with respect to the supremum metric, and hence is
a closed subset of $C_b(M, M)$.  If $M$ is compact, then
$\mathcal{I}(M)$ is also compact, with respect to the topology
determined by the supremum metric, because $\mathcal{I}(M, M)$
is compact.

        In fact, if $M$ is compact and $f$ is an isometry of $M$ into
itself, then $f(M) = M$.  To see this, suppose for the sake of a contradiction
that there is an element $x_1$ of $M$ not in $f(M)$.  Because $M$ is compact,
$f(M)$ is compact, and hence there is an $r > 0$ such that
\begin{equation}
\label{d(x_1, f(y)) ge r}
        d(x_1, f(y)) \ge r
\end{equation}
for every $y \in M$.  If $\{x_j\}_{j = 1}^\infty$ is the sequence of
elements of $M$ defined recursively by $x_{j + 1} = f(x_j)$ for each
$j \in {\bf Z}_+$, then one can check that
\begin{equation}
\label{d(x_j, x_k) ge r}
        d(x_j, x_k) \ge r
\end{equation}
when $j < k$, using (\ref{d(x_1, f(y)) ge r}) and the hypothesis that
$f$ be an isometry.  This implies that $\{x_j\}_{j = 1}^\infty$ has no
convergent subsequences, contradicting the compactness of $M$, as
desired.

        As a variant of this, suppose that $f$ and $g$ are similarities
from $M$ into $N$, with the same constant $\lambda$.  If $g(M) = N$,
then $f \circ g^{-1}$ is an isometry from $M$ into itself.  If $M$
is compact, then $f \circ g^{-1}$ maps $M$ onto itself, as in the
previous paragraph.  Thus $f(M) = N$ under these conditions.

        Suppose now that $d(x, y)$ is an ultrametric on $M$, and let 
$r > 0$ be given.  Also let $\sim_r$ be the relation on $M$ defined by 
$x \sim_r y$ when $d(x, y) \le r$.  This is an equivalence relation on
$M$, because $d(x, y)$ is an ultrametric on $M$.  The corresponding
equivalence classes are closed balls of radius $r$ in $M$.  If $f$ is
an isometry of $M$ into itself, then $f(x) \sim_r f(y)$ if and only if
$x \sim_r y$ for every $x, y \in M$.  This implies that $f$ maps each
equivalence class of $M$ with respect to $\sim_r$ into another
equivalence class, which is the same as saying that $f$ maps each
closed ball in $M$ with radius $r$ into another closed ball of radius
$r$.  More precisely, $f$ maps distinct equivalence classes in $M$
with respect to $\sim_r$ into distinct equivalence classes in $M$,
which is the same as saying that $f$ maps disjoint closed balls in
$M$ with radius $r$ into disjoint closed balls with radius $r$.

        If $M$ is totally bounded, then there are only finitely many
equivalence classes in $M$ with respect to $\sim_r$ for each $r > 0$.
In this case, it follows that every such equivalence class contains
an element of $f(M)$.  This means that $f(M)$ is dense in $M$, since
this holds for each $r > 0$.  If $M$ is complete, then we get that
$f(M) = M$.  Of course, $M$ is compact when $M$ is complete and totally
bounded.

        Let $d(x, y)$ be an arbitrary metric on $M$ again, and suppose
that $M$ is totally bounded.  Let $r$ be a positive real number, and
let $n(r)$ be the smallest number of subsets of $M$ with diameter less
than or equal to $r$ needed to cover $M$.  If $E$ is any subset of $M$
which is not dense in $M$, then $E$ can be covered by fewer than
$n(r)$ subsets of $M$ with diameter less than or equal to $r$ when
$r$ is sufficiently small, because at least one of the sets used to
cover $M$ will not intersect $E$.  If $f$ is an isometry of $M$ into
itself, then the minimal number of sets of diameter less than or equal
to $r$ needed to cover $f(M)$ is the same as $n(r)$.  This implies that
$f(M)$ is dense in $M$ when $M$ is totally bounded, and hence that
$f(M) = M$ when $M$ is also complete and thus compact.

        Suppose that $H$ is a Hausdorff measure on $M$, defined in terms
of some function of the diameter of subsets of $M$.  Thus $H(f(M)) =
H(M)$ when $f$ is an isometry of $M$ into itself.  If $H(M) < +\infty$
and nonempty open subsets of $M$ have positive measure with respect to
$H$, then it follows that $f(M)$ is dense in $M$, so that $f(M) = M$
when $M$ is compact.  One can also show that $f(M) = M$ when $f$ is an
isometry from $M$ into itself and $M$ is compact using compactness of
$\mathcal{I}(M, M)$.  The covering argument in the preceding paragraph
and the earlier approach using sequential compactness were suggested by
students in a class, and some instances of this type of situation will
be discussed in the next chapter.

\chapter{Functions on ${\bf Q}_p$}
\label{functions on Q_p}

\section{Polynomials on ${\bf Q}_p$}
\label{polynomials on Q_p}

        Let $p$ be a prime number, and let
\begin{equation}
\label{f(x) = a_n x^n + a_{n - 1} x^{n - 1} + cdots + a_1 x + a_0}
        f(x) = a_n \, x^n + a_{n - 1} \, x^{n - 1} + \cdots + a_1 \, x + a_0
\end{equation}
be a polynomial with coefficients in ${\bf Q}_p$.  Of course,
\begin{equation}
\label{(x + h)^k = sum_{j = 0}^k {k choose j} h^j x^{k - j}}
        (x + h)^k = \sum_{j = 0}^k {k \choose j} \, h^j \, x^{k - j}
\end{equation}
for every nonnegative integer $k$ and $x, h \in {\bf Q}_p$, where ${k
  \choose j}$ is the usual binomial coefficient.  Thus
\begin{equation}
\label{f(x + h) = sum_{k = 0}^n a_k (x + h)^k = ...}
        f(x + h) = \sum_{k = 0}^n a_k \, (x + h)^k 
    = \sum_{k = 0}^n \sum_{j = 0}^k a_k \, {k \choose j} \, h^j \, x^{k - j}
\end{equation}
for every $x, h \in {\bf Q}_p$.  This implies that
\begin{equation}
\label{f(x + h) - f(x) = ...}
        f(x + h) - f(x) = \sum_{k = 1}^n \sum_{j = 1}^k a_k \, {k \choose j} \, 
                                                            h^j \, x^{k - j},
\end{equation}
by subtracting the $j = 0$ terms from (\ref{f(x + h) = sum_{k = 0}^n
  a_k (x + h)^k = ...}), and using the simple fact that ${k \choose 0}
= 1$ for each $k$.

        The formal derivative of $f(x)$ is the polynomial defined by
\begin{equation}
\label{f'(x) = n a_n x^{n - 1} + (n - 1) a_{n - 1} x^{n - 2} + cdots + a_1}
 f'(x) = n \, a_n \, x^{n - 1} + (n - 1) \, a_{n - 1} \, x^{n - 2} + \cdots + a_1.
\end{equation}
Subtracting the $j = 1$ terms from (\ref{f(x + h) - f(x) = ...}), we get that
\begin{equation}
\label{f(x + h) - f(x) - f'(x) h = ...}
 f(x + h) - f(x) - f'(x) \, h = \sum_{k = 2}^n \sum_{j = 2}^k a_k \, 
                                          {k \choose j} \, h^j \, x^{k - j}
\end{equation}
for every $x, h \in {\bf Q}_p$, because ${k \choose 1} = k$.  In particular,
\begin{equation}
\label{lim_{h to 0} frac{f(x + h) - f(x)}{h} = f'(x)}
        \lim_{h \to 0} \frac{f(x + h) - f(x)}{h} = f'(x)
\end{equation}
for every $x \in {\bf Q}_p$, since each term on the right side of
(\ref{f(x + h) - f(x) - f'(x) h = ...}) is a multiple of $h^2$.

        It follows from (\ref{f(x + h) - f(x) = ...}) that $f(x)$
is Lipschitz of order $1$ on bounded subsets of ${\bf Q}_p$.  Suppose
now that $a_k \in {\bf Z}_p$ for each $k$, so that $f$ maps ${\bf
  Z}_p$ into itself.  In this case, (\ref{f(x + h) - f(x) = ...})
implies that
\begin{equation}
\label{|f(x + h) - f(x)|_p le |h|_p}
        |f(x + h) - f(x)|_p \le |h|_p
\end{equation}
for every $x, h \in {\bf Z}_p$, since the binomial coefficients ${k
  \choose j}$ are integers.  Of course, the coefficients of $f'(x)$
are elements of ${\bf Z}_p$ too, so that
\begin{equation}
\label{|f'(x + h) - f'(x)| le |h|_p}
        |f'(x + h) - f'(x)| \le |h|_p
\end{equation}
for every $x, h \in {\bf Z}_p$ as well.  Using (\ref{f(x + h) - f(x) -
  f'(x) h = ...}), we also get that
\begin{equation}
\label{|f(x + h) - f(x) - f'(x) h|_p le |h|_p^2}
        |f(x + h) - f(x) - f'(x) \, h|_p \le |h|_p^2
\end{equation}
for every $x, h \in {\bf Z}_p$ under these conditions.

\section{Hensel's lemma (first version)}
\label{hensel's lemma (first version)}

        Let $f(x)$ be a polynomial with coefficients in ${\bf Z}_p$,
so that $f(x)$ and $f'(x)$ are elements of ${\bf Z}_p$ for every $x
\in {\bf Z}_p$.  Suppose that $x_0 \in {\bf Z}_p$ satisfies $f(x_0)
\in p \, {\bf Z}_p$ and $|f'(x_0)|_p = 1$.  Under these conditions,
\emph{Hensel's lemma}\index{Hensel's lemma} states that there is an $x
\in {\bf Z}_p$ such that $x - x_0 \in p \, {\bf Z}_p$ and $f(x) = 0$.
The proof uses Newton's method, as follows.  If $x_1 \in {\bf Z}_p$ is
close to $x_0$, then $f(x_1)$ is approximately
\begin{equation}
\label{f(x_0) + f'(x_0) (x_1 - x_0)}
        f(x_0) + f'(x_0) \, (x_1 - x_0),
\end{equation}
as in (\ref{|f(x + h) - f(x) - f'(x) h|_p le |h|_p^2}).  In order to
make this approximation equal to $0$, we take
\begin{equation}
\label{x_1 = x_0 - f'(x_0)^{-1} f(x_0)}
        x_1 = x_0 - f'(x_0)^{-1} \, f(x_0).
\end{equation}
This satisfies $x_1 - x_0 \in p \, {\bf Z}_p$, since $f(x_0) \in p \,
{\bf Z}_p$ and $|f'(x_0)|_p = 1$.

        Repeating the process, we shall choose a sequence of elements
$x_1, x_2, x_3, \ldots$ of ${\bf Z}_p$ such that
\begin{equation}
\label{x_j - x_{j - 1} in p {bf Z}_p}
        x_j - x_{j - 1} \in p \, {\bf Z}_p
\end{equation}
for each $j \ge 1$.  In particular, this ensures that
\begin{equation}
\label{x_j - x_0 in p {bf Z}_p}
        x_j - x_0 \in p \, {\bf Z}_p
\end{equation}
for every $j \ge 1$, and hence that $f(x_j) - f(x_0) \in p \, {\bf
  Z}_p$ for every $j \ge 1$, by (\ref{|f(x + h) - f(x)|_p le |h|_p}).
Of course, this implies that
\begin{equation}
\label{f(x_j) in p {bf Z}_p}
        f(x_j) \in p \, {\bf Z}_p
\end{equation}
for every $j \ge 1$, since $f(x_0) \in p \, {\bf Z}_p$ by hypothesis.
Similarly,
\begin{equation}
\label{f'(x_j) - f'(x_0) in p {bf Z}_p}
        f'(x_j) - f'(x_0) \in p \, {\bf Z}_p
\end{equation}
for each $j \ge 1$, by (\ref{x_j - x_0 in p {bf Z}_p}) and (\ref{|f'(x
  + h) - f'(x)| le |h|_p}).  It follows that
\begin{equation}
\label{|f'(x_j)|_p = 1}
        |f'(x_j)|_p = 1
\end{equation}
for every $j \ge 1$, since $|f'(x_0)|_p = 1$ by hypothesis.

        If $x_{j - 1}$ has already been chosen, then we would like to
choose $x_j$ so that
\begin{equation}
\label{f(x_{j - 1}) + f'(x_{j - 1}) (x_j - x_{j - 1}) = 0}
        f(x_{j - 1}) + f'(x_{j - 1}) \, (x_j - x_{j - 1}) = 0,
\end{equation}
which is the same as saying that
\begin{equation}
\label{x_j = x_{j - 1} - f'(x_{j - 1}) f(x_{j - 1})}
        x_j = x_{j - 1} - f'(x_{j - 1}) \, f(x_{j - 1}).
\end{equation}
In particular, if $x_{j - 1} - x_0 \in p \, {\bf Z}_p$, then $f(x_{j -
  1}) \in p \, {\bf Z}_p$ and $|f'(x_{j - 1})|_p = 1$, as in the
previous paragraph.  This implies that (\ref{x_j - x_{j - 1} in p {bf
    Z}_p}) holds, so that the process can be repeated.  More precisely,
\begin{equation}
\label{|x_j - x_{j - 1}|_p = |f(x_{j - 1})|_p}
        |x_j - x_{j - 1}|_p = |f(x_{j - 1})|_p.
\end{equation}

        Under these conditions, we also have that
\begin{equation}
\label{|f(x_j)|_p le |x_j - x_{j - 1}|_p^2}
        |f(x_j)|_p \le |x_j - x_{j - 1}|_p^2,
\end{equation}
by applying (\ref{|f(x + h) - f(x) - f'(x) h|_p le |h|_p^2}) to $x =
x_{j - 1}$ and $h = x_j - x_{j - 1}$, and using (\ref{f(x_{j - 1}) +
  f'(x_{j - 1}) (x_j - x_{j - 1}) = 0}).  Combining this with
(\ref{|x_j - x_{j - 1}|_p = |f(x_{j - 1})|_p}), we get that
\begin{equation}
\label{|f(x_j)|_p le |f(x_{j - 1})|_p^2}
        |f(x_j)|_p \le |f(x_{j - 1})|_p^2
\end{equation}
for each $j \ge 1$.  This implies that $|f(x_j)|_p \to 0$ as $j \to
\infty$, because $|f(x_0)|_p < 1$, by hypothesis.  Thus $|x_j - x_{j -
  1}|_p \to 0$ as $j \to \infty$, by (\ref{|x_j - x_{j - 1}|_p =
  |f(x_{j - 1})|_p}) again.  It follows that $\{x_j\}_{j = 1}^\infty$
is a Cauchy sequence in ${\bf Z}_p$, as in Section \ref{sequences,
  series}, since the $p$-adic metric is an ultrametric.  By
completeness, $\{x_j\}_{j = 1}^\infty$ converges to an element $x$ of
${\bf Z}_p$, and in fact $x - x_0 \in {\bf Z}_p$, because of (\ref{x_j
  - x_0 in p {bf Z}_p}).  Of course, $f(x) = 0$, as desired, because
$f$ is continuous on ${\bf Q}_p$, and $|f(x_j)|_p \to 0$ as $j \to
\infty$.

\section{Hensel's lemma (second version)}
\label{hensel's lemma (second version)}

        Let $f(x)$ be a polynomial with coefficients in ${\bf Z}_p$ again, 
and suppose that $x_0 \in {\bf Z}_p$ satisfies
\begin{equation}
\label{|f(x_0)|_p < |f'(x_0)|_p^2}
        |f(x_0)|_p < |f'(x_0)|_p^2.
\end{equation}
We would like to find an $x \in {\bf Z}_p$ that is close to $x_0$ and
satisfies $f(x) = 0$.  Of course, $f(x_0), f'(x_0) \in {\bf Z}_p$, so
that $|f(x_0)|_p, |f'(x_0)|_p \le 1$.  If $|f'(x_0)|_p = 1$, then we
are back in the situation discussed in the previous section.
Otherwise, Newton's method is still applicable, but we should be a bit
more careful about some of the estimates.

        Let $j$ be a positive integer, and suppose that $x_{j - 1} \in 
{\bf Z}_p$ has been chosen in such a way that
\begin{equation}
\label{|x_{j - 1} - x_0| < |f'(x_0)|_p}
        |x_{j - 1} - x_0| < |f'(x_0)|_p
\end{equation}
and
\begin{equation}
\label{|f(x_{j - 1})|_p le |f(x_0)|_p}
        |f(x_{j - 1})|_p \le |f(x_0)|_p.
\end{equation}
Thus
\begin{equation}
\label{|f'(x_{j - 1}) - f'(x_0)|_p le |x_{j - 1} - x_0|_p < |f'(x_0)|_p}
        |f'(x_{j - 1}) - f'(x_0)|_p \le |x_{j - 1} - x_0|_p < |f'(x_0)|_p,
\end{equation}
by (\ref{|f'(x + h) - f'(x)| le |h|_p}), which implies that
\begin{equation}
\label{|f'(x_{j - 1})|_p = |f'(x_0)|_p}
        |f'(x_{j - 1})|_p = |f'(x_0)|_p,
\end{equation}
because of the ultrametric version of the triangle inequality.  Let us
choose $x_j \in {\bf Q}_p$ as in (\ref{x_j = x_{j - 1} - f'(x_{j - 1})
  f(x_{j - 1})}), so that
\begin{equation}
\label{|x_j - x_{j - 1}|_p = ... = |f'(x_0)|_p^{-1} |f(x_{j - 1})|_p}
        |x_j - x_{j - 1}|_p = |f'(x_{j - 1})|_p^{-1} \, |f(x_{j - 1})| 
                           = |f'(x_0)|_p^{-1} \, |f(x_{j - 1})|_p.
\end{equation}
Combining this with (\ref{|f(x_0)|_p < |f'(x_0)|_p^2}) and
(\ref{|f(x_{j - 1})|_p le |f(x_0)|_p}), we get that
\begin{equation}
\label{|x_j - x_{j - 1}|_p < |f'(x_0)|_p}
        |x_j - x_{j - 1}|_p < |f'(x_0)|_p.
\end{equation}
It follows that
\begin{equation}
\label{|x_j - x_0|_p < |f'(x_0)|_p}
        |x_j - x_0|_p < |f'(x_0)|_p,
\end{equation}
by (\ref{|x_{j - 1} - x_0| < |f'(x_0)|_p}), and in particular that
$x_j \in {\bf Z}_p$.  This permits us to apply (\ref{|f(x + h) - f(x)
  - f'(x) h|_p le |h|_p^2}) with $x = x_{j - 1}$ and $h = x_j - x_{j -
  1}$, to get that
\begin{equation}
\label{|f(x_j)|_p le |x_j - x_{j - 1}|_p^2, 2}
        |f(x_j)|_p \le |x_j - x_{j - 1}|_p^2,
\end{equation}
using also (\ref{f(x_{j - 1}) + f'(x_{j - 1}) (x_j - x_{j - 1}) = 0}).
This implies that
\begin{equation}
\label{|f(x_j)|_p le |f(x_{j - 1})|_p}
        |f(x_j)|_p \le |f(x_{j - 1})|_p,
\end{equation}
by (\ref{|x_j - x_{j - 1}|_p = ... = |f'(x_0)|_p^{-1} |f(x_{j -
    1})|_p}) and (\ref{|x_j - x_{j - 1}|_p < |f'(x_0)|_p}).  In particular,
\begin{equation}
\label{|f(x_j)|_p le |f(x_0)|_p}
        |f(x_j)|_p \le |f(x_0)|_p,
\end{equation}
by (\ref{|f(x_{j - 1})|_p le |f(x_0)|_p}).  This and (\ref{|x_j -
  x_0|_p < |f'(x_0)|_p}) show that $x_j$ satisfies the same conditions
as $x_{j - 1}$, so that the process can be repeated.

        More precisely, (\ref{|x_j - x_{j - 1}|_p = ... = |f'(x_0)|_p^{-1} 
|f(x_{j - 1})|_p}) and (\ref{|f(x_j)|_p le |x_j - x_{j - 1}|_p^2, 2})
imply that
\begin{equation}
\label{|f(x_j)|_p le |f'(x_0)|_p^{-2} |f(x_{j - 1})|_p^2}
        |f(x_j)|_p \le |f'(x_0)|_p^{-2} \, |f(x_{j - 1})|_p^2
\end{equation}
for each $j \ge 1$.   Thus
\begin{equation}
        |f(x_j)|_p \le (|f'(x_0)|_p^{-2} \, |f(x_0)|_p) \, |f(x_{j - 1})|_p
\end{equation}
for each $j \ge 1$, by (\ref{|f(x_{j - 1})|_p le |f(x_0)|_p}).  This
implies that $|f(x_j)|_p \to 0$ as $j \to \infty$, by (\ref{|f(x_0)|_p
  < |f'(x_0)|_p^2}).

        It follows from this and (\ref{|x_j - x_{j - 1}|_p = ... = 
|f'(x_0)|_p^{-1} |f(x_{j - 1})|_p}) that $|x_j - x_{j - 1}|_p \to 0$ as 
$j \to \infty$.  Thus $\{x_j\}_{j = 1}^\infty$ is a Cauchy sequence
in ${\bf Z}_p$, as in Section \ref{sequences, series}, which converges
to an element $x$ of ${\bf Z}_p$, by completeness.  Note that
\begin{equation}
\label{|x - x_0|_p < |f'(x_0)|_p}
        |x - x_0|_p < |f'(x_0)|_p,
\end{equation}
because of the analogous condition for the $x_j$'s, and the fact that
open balls in ultrametric spaces are closed sets.  We also have that
$f(x) = 0$, as desired, because $f$ is continuous on ${\bf Q}_p$,
and $|f(x_j)|_p \to 0$ as $j \to \infty$.

\section{Contractions}
\label{contractions}

        Let $f(x)$ be a polynomial with coefficients in ${\bf Z}_p$,
and suppose that $x_0 \in {\bf Z}_p$ satisfies (\ref{|f(x_0)|_p <
  |f'(x_0)|_p^2}).  Consider
\begin{equation}
\label{g(x) = x - (f'(x_0))^{-1} f(x + x_0) - x}
        g(x) = x - (f'(x_0))^{-1} \, f(x + x_0),
\end{equation}
which is a polynomial with coefficients in ${\bf Q}_p$ that satisfies
\begin{equation}
\label{g'(0) = 1 - (f'(x_0))^{-1} f'(x_0) = 0}
        g'(0) = 1 - (f'(x_0))^{-1} \, f'(x_0) = 0.
\end{equation}
More precisely, $|f'(x_0)|_p = p^{-k}$ for some nonnegative integer
$k$, which implies that the coefficients of $g$ are in $p^{-k} \, {\bf
  Z}_p$.  Using (\ref{|f(x_0)|_p < |f'(x_0)|_p^2}), we also get that
\begin{equation}
\label{|g(0)| = |f'(x_0)|_p^{-1} |f(x_0)| < |f'(x_0)|_p = p^{-k}}
        |g(0)| = |f'(x_0)|_p^{-1} \, |f(x_0)| < |f'(x_0)|_p = p^{-k}.
\end{equation}

        Let us now start over, and let $k$ be a nonnegative integer
and $g(x)$ be a polynomial with coefficients in $p^{-k} \, {\bf Z}_p$
such that
\begin{equation}
\label{g(0), g'(0) in p^{k + 1} {bf Z}_p}
        g(0), \, g'(0) \in p^{k + 1} \, {\bf Z}_p.
\end{equation}
If $x \in p^{k + 1} \, {\bf Z}_p$, then it is easy to see that
\begin{equation}
\label{g(x) in p^{k + 1} {bf Z}_p}
        g(x) \in p^{k + 1} \, {\bf Z}_p
\end{equation}
and
\begin{equation}
\label{g'(x) in p {bf Z}_p}
        g'(x) \in p \, {\bf Z}_p.
\end{equation}
Observe that
\begin{equation}
\label{|g(y) - g(x) - g'(x) (y - x)|_p le p^k |x - y|_p^2}
        |g(y) - g(x) - g'(x) \, (y - x)|_p \le p^k \, |x - y|_p^2
\end{equation}
for every $x, y \in {\bf Z}_p$, by applying (\ref{|f(x + h) - f(x) -
  f'(x) h|_p le |h|_p^2}) to $p^k \, g$.  This implies that
\begin{equation}
\label{|g(y) - g(x)|_p le max(|g'(x)|_p, p^k |x - y|_p) |x - y|_p}
        |g(y) - g(x)|_p \le \max(|g'(x)|_p, p^k \, |x - y|_p) \, |x - y|_p 
                         \le p^{-1} \, |x - y|_p
\end{equation}
when $x, y  \in p^{k + 1} \, {\bf Z}_p$, by (\ref{g'(x) in p {bf Z}_p}).

        Thus $g$ maps $p^{k + 1} \, {\bf Z}_p$ into itself under these 
conditions, and the restriction of $g$ to $p^{k + 1} \, {\bf Z}_p$ is
a strict contraction, by (\ref{|g(y) - g(x)|_p le max(|g'(x)|_p, p^k
  |x - y|_p) |x - y|_p}).  The contraction mapping principle implies
that $g$ has a unique fixed point in $p^{k + 1} \, {\bf Z}_p$, because
$p^{k + 1} \, {\bf Z}_p$ is complete as a metric space.  If $g$ is as
in (\ref{g(x) = x - (f'(x_0))^{-1} f(x + x_0) - x}), then this is the
same as saying that there is a unique $x \in p^{k + 1} \, {\bf Z}_p$
such that $f(x + x_0) = 0$.

\section{Local geometry}
\label{local geometry}

        Let $f(x)$ be a polynomial with coefficients in ${\bf Z}_p$,
and suppose that $x_0 \in {\bf Z}_p$ satisfies $f'(x_0) \ne 0$.  Let
$k$ be a nonnegative integer such that $|f'(x_0)|_p = p^{-k}$, as
before.  If $x \in x_0 + p^{k + 1} \, {\bf Z}_p$, then
\begin{equation}
\label{|f'(x) - f'(x_0)|_p le |x - x_0|_p le p^{-k - 1}}
        |f'(x) - f'(x_0)|_p \le |x - x_0|_p \le p^{-k - 1},
\end{equation}
by (\ref{|f'(x + h) - f'(x)| le |h|_p}).  This implies that
\begin{equation}
\label{|f'(x)|_p = |f'(x_0)|_p = p^{-k}}
        |f'(x)|_p = |f'(x_0)|_p = p^{-k},
\end{equation}
by the ultrametric version of the triangle inequality.  If $x, y \in
x_0 + p^{k + 1} \, {\bf Z}_p$, then
\begin{equation}
\label{|f(y) - f(x) - f'(x) (y - x)|_p le |x - y|_p^2 le p^{-k - 1} |x - y|_p}
 |f(y) - f(x) - f'(x) \, (y - x)|_p \le |x - y|_p^2 \le p^{-k - 1} \, |x - y|_p,
\end{equation}
as in (\ref{|f(x + h) - f(x) - f'(x) h|_p le |h|_p^2}).  It follows that
\begin{eqnarray}
\label{|f(x) - f(y)|_p le ... = p^{-k} |x - y|_p}
 |f(x) - f(y)|_p & \le & \max(|f'(x)|_p \, |x - y|_p, p^{-k - 1} \, |x - y|_p) \\
                   & = & p^{-k} \, |x - y|_p,                    \nonumber
\end{eqnarray}
by (\ref{|f'(x)|_p = |f'(x_0)|_p = p^{-k}}).  Similarly,
(\ref{|f'(x)|_p = |f'(x_0)|_p = p^{-k}}) and (\ref{|f(y) - f(x) -
  f'(x) (y - x)|_p le |x - y|_p^2 le p^{-k - 1} |x - y|_p}) also imply that
\begin{eqnarray}
\label{p^{-k} |x - y|_p = ... le max(|f(x) - f(y)|_p, p^{k - 1} |x - y|_p)}
        p^{-k} \, |x - y|_p & = & |f'(x)|_p \, |x - y|_p  \\
                & \le & \max(|f(x) - f(y)|_p, p^{k - 1} \, |x - y|_p), \nonumber
\end{eqnarray}
and hence that
\begin{equation}
\label{p^{-k} |x - y|_p le |f(x) - f(y)|_p}
        p^{-k} \, |x - y|_p \le |f(x) - f(y)|_p.
\end{equation}
This shows that
\begin{equation}
\label{|f(x) - f(y)|_p = p^{-k} |x - y|_p}
        |f(x) - f(y)|_p = p^{-k} \, |x - y|_p
\end{equation}
for every $x, y \in x_0 + p^{k + 1} \, {\bf Z}_p$ under these
conditions.

        Now let $f$ be any mapping from $x_0 + p^{k + 1} \, {\bf Z}_p$
into ${\bf Q}_p$ that satisfies (\ref{|f(x) - f(y)|_p = p^{-k} |x - y|_p}) 
for every $x, y \in p^{p + 1} \, {\bf Z}_p$, where $k$ is a nonnegative
integer.  In particular,
\begin{equation}
\label{f(x_0 + p^{k + 1} {bf Z}_p) subseteq f(x_0) + p^{2 k + 1} {bf Z}_p}
 f(x_0 + p^{k + 1} \, {\bf Z}_p) \subseteq f(x_0) + p^{2 k + 1} \, {\bf Z}_p.
\end{equation}
Of course,
\begin{equation}
\label{x mapsto p^{-k} (x - x_0) + f(x_0)}
        x \mapsto p^{-k} \, (x - x_0) + f(x_0)
\end{equation}
is a similarity from $x_0 + p^{k + 1} \, {\bf Z}_p$ onto $f(x_0) +
p^{2 \, k + 1} \, {\bf Z}_p$ with respect to the $p$-adic metric, with
the same similarity constant $p^{-k}$.  Because $x_0 + p^{k + 1} \,
{\bf Z}_p$ is compact, one can use this to show that
\begin{equation}
\label{f(x_0 + p^{k + 1} {bf Z}_p) = f(x_0) + p^{2 k + 1} {bf Z}_p}
        f(x_0 + p^{k + 1} \, {\bf Z}_p) = f(x_0) + p^{2 k + 1} \, {\bf Z}_p,
\end{equation}
as in Section \ref{isometries, similarities}.

        Remember that the one-dimensional Hausdorff measure of 
$x_0 + p^{k + 1} \, {\bf Z}_p$ with respect to the $p$-adic metric is
equal to $p^{-k - 1}$, as in Section \ref{some special cases}.  Using
this and (\ref{|f(x) - f(y)|_p = p^{-k} |x - y|_p}), it is easy to see
that the one-dimensional Hausdorff measure of $f(x_0 + p^{k + 1} \,
{\bf Z}_p)$ is equal to $p^{-2 k - 1}$, as in Section \ref{isometries,
  similarities}.  Of course, this is the same as the one-dimensional
Hausdorff measure of $f(x_0) + p^{2 k + 1} \, {\bf Z}_p$.  It follows
that $f(x_0 + p^k \, {\bf Z}_p)$ is dense in $f(x_0) + p^{2 k + 1} \,
{\bf Z}_p$, since every ball in ${\bf Q}_p$ with positive radius has
positive one-dimensional Hausdorff measure.  Note that $f(x_0 + p^k \,
{\bf Z}_p)$ is a compact set in ${\bf Q}_p$, because $x_0 + p^k \,
{\bf Z}_p$ is compact, and $f$ is continuous.  Thus the density of
$f(x_0 + p^k \, {\bf Z}_p)$ in $f(x_0) + p^{2 k + 1} \, {\bf Z}_p$
implies that (\ref{f(x_0 + p^{k + 1} {bf Z}_p) = f(x_0) + p^{2 k + 1}
  {bf Z}_p}) holds, as before.  This type of argument was also mentioned
in Section \ref{isometries, similarities}, using the properties of
one-dimensional Hausdorff measure in this case.

        Here is an analogous but more elementary approach, which
is a more explicit version of another argument in Section
\ref{isometries, similarities} in this situation.  If $n$ is a
positive integer, then there is a set $A_n \subseteq p^{k + 1} \, {\bf
  Z}_p$ with exactly $p^n$ elements such that
\begin{equation}
\label{|a - b|_p ge p^{-k - n}}
        |a - b|_p \ge p^{-k - n}
\end{equation}
for every $a, b \in A_n$ with $a \ne b$.  Equivalently, this means
that the restriction of the natural quotient mapping from $p^{k + 1}
\, {\bf Z}_p$ onto $p^{k + 1} \, {\bf Z}_p / p^{k + n + 1} \, {\bf
  Z}_p$ to $A_n$ is injective.  If $f$ is as in the previous
two paragraphs, then
\begin{equation}
\label{|f(x_0 + a) - f(x_0 + b)|_p = p^{-k} |a - b|_p ge p^{-2 k - n}}
        |f(x_0 + a) - f(x_0 + b)|_p = p^{-k} \, |a - b|_p \ge p^{-2 k - n}
\end{equation}
for every $a, b \in A_n$ with $a \ne b$.  However, $f(x_0) + p^{2 k +
  1} \, {\bf Z}_p$ can be expressed as the union of $p^n$
pairwise-disjoint closed balls of radius $p^{-2 k - n - 1}$ for each
positive integer $n$.  Each of these balls can contain at most one
element of $f(x_0 + A_n)$, by (\ref{|f(x_0 + a) - f(x_0 + b)|_p =
  p^{-k} |a - b|_p ge p^{-2 k - n}}).  It follows that each of these
balls must contain an element of $f(x_0 + A_n)$, which implies that
$f(x_0 + p^{k + 1} \, {\bf Z}_p)$ is dense in $f(x_0) + p^{2 k + 1} \,
{\bf Z}_p$ under these conditions.

\section{Power series}
\label{power series}

        Let $\sum_{j = 0}^\infty a_j \, x^j$ be a power series with
coefficients in ${\bf Q}_p$, where $x^j$ is interpreted as being equal
to $1$ for every $x \in {\bf Q}_p$ when $j = 0$, as usual.  As in
Section \ref{sequences, series}, $\sum_{j = 0}^\infty a_j \, x^j$
converges for some $x \in {\bf Q}_p$ if and only if $\{a_j \, x^j\}_{j
  = 0}^\infty$ converges to $0$ in ${\bf Q}_p$, which is the same as
saying that
\begin{equation}
\label{|a_j x^j|_p = |a_j|_p |x|_p^j to 0 as j to infty}
        |a_j \, x^j|_p = |a_j|_p \, |x|_p^j \to 0 \hbox{ as } j \to \infty.
\end{equation}
In this case, $\sum_{j = 0}^\infty a_j \, y^j$ also converges in ${\bf
  Q}_p$ when $y \in {\bf Q}_p$ satisfies $|y|_p \le |x|_p$.  More precisely,
\begin{eqnarray}
\label{|sum_{j = 0}^infty a_j y^j - sum_{j = 0}^n a_j y^j| = ...}
 \biggl|\sum_{j = 0}^\infty a_j \, y^j - \sum_{j = 0}^n a_j \, y^j\biggr|_p
        & = & \biggl|\sum_{j = n + 1}^\infty a_j \, y^j\biggr|_p  \\
 & \le & \max_{j \ge n + 1} |a_j \, y^j|_p \le \max_{j \ge n + 1} |a_j|_p |x|_p^j
                                                        \nonumber
\end{eqnarray}
when $|y|_p \le |x|_p$, which implies that the partial sums $\sum_{j =
  0}^\infty a_j \, y^j$ converge to $\sum_{j = 0}^\infty a_j \, y^j$
uniformly as $n \to \infty$ on the set of $y \in {\bf Q}_p$ with
$|y|_p \le |x|_p$.  It follows that $\sum_{j = 0}^\infty a_j \, x^j$
defines a continuous ${\bf Q}_p$-valued function on the set of $x \in
{\bf Q}_p$ for which the series converges, which is either a closed
disk centered at $0$ or all of ${\bf Q}_p$.

          Now let $\sum_{j = 0}^\infty a_j \, x^j$ and
$\sum_{k = 0}^\infty b_k \, x^k$ be infinite series with coefficients
in ${\bf Q}_p$, and let
\begin{equation}
\label{c_l = sum_{j = 0}^l a_j b_{l - j}, 2}
        c_l = \sum_{j = 0}^l a_j \, b_{l - j}
\end{equation}
be the Cauchy product\index{Cauchy products} of their coefficients,
as in Section \ref{sequences, series}.  Observe that
\begin{equation}
\label{c_l x^l = sum_{j = 0}^n (a_j x^j) (b_{l - j} x^{l - j})}
        c_l \, x^l = \sum_{j = 0}^n (a_j \, x^j) \, (b_{l - j} \, x^{l - j})
\end{equation}
for each $x \in {\bf Q}_p$, so that $\sum_{l = 0}^\infty c_l \, x^l$
is the Cauchy product of $\sum_{j = 0}^\infty a_j \, x^j$ and $\sum_{k
  = 0}^\infty b_k \, x^k$.  In particular,
\begin{equation}
\label{sum_{l = 0}^infty c_l x^l = ...}
        \sum_{l = 0}^\infty c_l \, x^l = \Big(\sum_{j = 0}^\infty a_j \, x^j\Big)
                                      \, \Big(\sum_{k = 0}^\infty b_k \, x^k\Big)
\end{equation}
formally, collecting all of the terms which are multiples of $x^l$ for
each $l \ge 0$.  If $\sum_{j = 0}^\infty a_j \, x^j$ and $\sum_{k =
  0}^\infty b_k \, x^k$ both converge for some $x \in {\bf Q}_p$, then
it follows that $\sum_{l = 0}^\infty c_l \, x^l$ converges and
satisfies (\ref{sum_{l = 0}^infty c_l x^l = ...}), as in Section
\ref{sequences, series}.

        Suppose that
\begin{equation}
\label{f(x) = sum_{j = 0}^infty a_j x^j}
        f(x) = \sum_{j = 0}^\infty a_j \, x^j
\end{equation}
is a power series with coefficients $a_j \in {\bf Z}_p$ for each $j
\ge 0$, and that $\{a_j\}_{j = 0}^\infty$ converges to $0$ in ${\bf Q}_p$.
This implies that (\ref{f(x) = sum_{j = 0}^infty a_j x^j}) converges
for every $x \in {\bf Z}_p$, and that $f(x)$ defines a continuous 
${\bf Q}_p$-valued function on ${\bf Z}_p$, as before.  Moreover,
the formal derivative
\begin{equation}
\label{f'(x) = sum_{j = 1}^infty j a_j x^j}
        f'(x) = \sum_{j = 1}^\infty j \, a_j \, x^j
\end{equation}
also has coefficients in ${\bf Z}_p$ that converge to $0$, and hence
defines a continuous function on ${\bf Z}_p$ as well.  It is easy to
see that $f(x)$ and $f'(x)$ satisfy the same estimates (\ref{|f(x + h)
  - f(x)|_p le |h|_p}), (\ref{|f'(x + h) - f'(x)| le |h|_p}), and
(\ref{|f(x + h) - f(x) - f'(x) h|_p le |h|_p^2}) as for polynomials
with coefficients in ${\bf Z}_p$, by approximating the corresponding
infinite series by their partial sums.  It follows that the results
for polynomials with coefficients in ${\bf Z}_p$ discussed in the
previous sections also work for power series of this type.

        In particular, the analogue of (\ref{|f(x + h) - f(x) - f'(x) h|_p 
le |h|_p^2}) in this context implies that the derivative of $f$ at any
point $x \in {\bf Z}_p$ exists and is equal to $f'(x)$.  There are analogous
statements for any convergent power series with coefficients in ${\bf Q}_p$.

\section{Linear mappings on ${\bf Q}_p^n$}
\label{linear mappings on {bf Q}_p^n}

        Let $n$ be a positive integer, and let ${\bf Q}_p^n$ be the set
of $n$-tuples $v = (v_1, \ldots, v_n)$ of elements of ${\bf Q}_p$.
As usual, this is a vector space over ${\bf Q}_p$ with respect to
coordinatewise addition and scalar multiplication.  Put
\begin{equation}
\label{||v||= max (|v_1|_p, ldots, |v_n|_p)}
        \|v\| = \max (|v_1|_p, \ldots, |v_n|_p)
\end{equation}
for each $v \in {\bf Q}_p$, and observe that
\begin{equation}
\label{||t v|| = |t|_p ||v||}
        \|t \, v\| = |t|_p \, \|v\|
\end{equation}
for every $v \in {\bf Q}_p^n$ and $t \in {\bf Q}_p$, and that
\begin{equation}
\label{||v + w|| le max(||v||, ||w||)}
        \|v + w\| \le \max(\|v\|, \|w\|)
\end{equation}
for every $v, w \in {\bf Q}_p^n$.  Thus $\|v\|$ is an
\emph{ultranorm}\index{ultranorms} on ${\bf Q}_p^n$, which is like a
norm on a real or complex vector space, except that it satifies the
ultrametric version of the triangle inequality (\ref{||v + w|| le
  max(||v||, ||w||)}).  It follows that
\begin{equation}
\label{d(v, w) = ||v - w||}
        d(v, w) = \|v - w\|
\end{equation}
defines an ultrametric on ${\bf Q}_p^n$, for which the corresponding
topology on ${\bf Q}_p^n$ is the same as the product topology associated
to the standard topology on ${\bf Q}_p$.

        Let $e_1, \ldots, e_n$ be the standard basis vectors for ${\bf Q}_p^n$,
so that the $j$th coordinate of $e_k$ is equal to $1$ when $j = k$ and to $0$
otherwise.  If $T$ is a linear mapping from ${\bf Q}_p^n$ into itself, then put
\begin{equation}
\label{||T||_{op} = max (||T(e_1)||, ldots, ||T(e_n)||)}
        \|T\|_{op} = \max (\|T(e_1)\|, \ldots, \|T(e_n)\|).
\end{equation}
The space of linear mappings from ${\bf Q}_p^n$ into itself is also a
vector space over ${\bf Q}_p$ with respect to the usual addition and
scalar multiplication of linear mappings, and it is easy to see that
$\|T\|_{op}$ defines an ultranorm on this vector space.  Each $v \in
{\bf Q}_p^n$ can be expressed as $v = \sum_{j = 1}^n v_j \, e_j$, and hence
\begin{equation}
\label{||T(v)|| le max_{1 le j le n} (|v_j|_p ||T(e_j)||) le ||T||_{op} ||v||}
        \|T(v)\| \le \max_{1 \le j \le n} (|v_j|_p \, \|T(e_j)\|)
                  \le \|T\|_{op} \, \|v\|.
\end{equation}
Thus $\|T\|_{op}$ is the same as the operator norm of $T$ associated
to the ultranorm $\|v\|$ on ${\bf Q}_p^n$, and
\begin{equation}
\label{||T_2 circ T_1||_{op} le ||T_1||_{op} ||T_2||_{op}}
        \|T_2 \circ T_1\|_{op} \le \|T_1\|_{op} \, \|T_2\|_{op}
\end{equation}
for any two linear mappings $T_1$, $T_2$ from ${\bf Q}_p^n$ into itself.

        If $\{a_{j, k}\}_{j, k = 1}^n$ is an $n \times n$ matrix with entries
in ${\bf Q}_p$, then
\begin{equation}
\label{(T(v))_j = sum_{k = 1}^n a_{j, k} v_k}
        (T(v))_j = \sum_{k = 1}^n a_{j, k} \, v_k
\end{equation}
defines a linear mapping from ${\bf Q}_p^n$ into itself, where
$(T(v))_j$ is the $j$th coordinate of $T(v)$.  Of course, every linear
mapping from ${\bf Q}_p^n$ into itself can be expressed in this way,
and one can check that
\begin{equation}
\label{||T||_{op} = max_{1 le j, k le n} |a_{j, k}|_p}
        \|T\|_{op} = \max_{1 \le j, k \le n} |a_{j, k}|_p
\end{equation}
when $T$ is as in (\ref{(T(v))_j = sum_{k = 1}^n a_{j, k} v_k}).  Note
that $\|T\|_{op} \le 1$ if and only if $a_{j, k} \in {\bf Z}_p$ for
each $j, k = 1, \ldots n$, which happens if and only if $T$ maps ${\bf
  Z}_p^n$ into itself.  If $T$ is as in (\ref{(T(v))_j = sum_{k = 1}^n
  a_{j, k} v_k}), then the determinant of $T$ as a linear mapping on
${\bf Q}_p^n$ is the same as the determinant of the corresponding
matrix $\{a_{j, k}\}_{j, k = 1}^n$, and hence
\begin{equation}
\label{|det T|_p le ||T||_{op}^n}
        |\det T|_p \le \|T\|_{op}^n.
\end{equation}

        Suppose that $T$ is a linear mapping from ${\bf Q}_p^n$ into
itself that satisfies
\begin{equation}
\label{||T(v)|| = ||v||}
        \|T(v)\| = \|v\|
\end{equation}
for every $v \in {\bf Q}_p^n$.  In particular, this implies that $T$
is one-to-one, and hence that $T$ maps ${\bf Q}_p^n$ onto itself, by
linear algebra.  Thus $T$ is an invertible linear mapping on ${\bf
  Q}_p^n$, and the operator norms of $T$ and $T^{-1}$ are both equal
to $1$.  Conversely, if $T$ is an invertible linear mapping on ${\bf
  Q}_p^n$ such that $\|T\|_{op}, \, \|T^{-1}\|_{op} \le 1$, then $T$
satisfies (\ref{||T(v)|| = ||v||}).  In this case, $T$ and $T^{-1}$
both correspond to matrices with entries in ${\bf Z}_p$, whose
determinants are in ${\bf Z}_p$ as well.  It follows that
\begin{equation}
\label{|det T|_p = 1}
        |\det T|_p = 1,
\end{equation}
because $\det T$ and $(\det T)^{-1} = \det T^{-1}$ are both in ${\bf
  Z}_p$.  Conversely, suppose that $T$ is a linear mapping on ${\bf
  Q}_p^n$ that corresponds to an $n \times n$ matrix with entries in
${\bf Z}_p$, and that $T$ satisfies (\ref{|det T|_p = 1}).  This
implies that $T$ is invertible, where the matrix associated to the
inverse of $T$ can be expressed in terms of determinants in the usual
way.  More precisely, the entries of the matrix associated to $T^{-1}$
are in ${\bf Z}_p$ too, because of (\ref{|det T|_p = 1}).

\chapter{Commutative topological groups}
\label{commutative topological groups}

\section{Haar measure}
\label{haar measure}

        Let $G$ be a group, in which the group operations are expressed
multiplicatively.  If $G$ is also equipped with a topology with respect
to which the group operations are continuous, then $G$ is said to be a
\emph{topological group}.\index{topological groups}  More precisely,
this means that multiplication in the group should be continuous as
a mapping from $G \times G$ into $G$, where $G \times G$ is equipped with
the product topology associated to the given topology on $G$.  Similarly,
$x \mapsto x^{-1}$ should be continuous as a mapping from $G$ onto itself.
It is customary to ask also that the set containing only the identity element
$e$ in $G$ be a closed set in $G$.  This implies that every one-element
subset of $G$ is closed, using the continuity of translations on $G$, which
follows from continuity of multiplication on $G$.  One can show that $G$ 
is Hausdorff under these conditions, and in fact regular as a topological
space.

        Let $G$ be a topological group which is locally compact as a
topological space.  It is well known that there is a nonnegative Borel
measure on $G$ with suitable regularity properties that is invariant
under left translations, known as \emph{Haar measure}.\index{Haar measure} 
In particular, the Haar measure of a nonempty open set in $G$ should
be positive, and the Haar measure of a compact set in $G$ should be
finite.  This measure is unique up to multiplication by a positive
real number.  Similarly, there is a nonnegative Borel measure on $G$
with suitable regularity properties that is invariant under right
translations, with the same type of uniqueness property.  Of course,
one can use the mapping $x \mapsto x^{-1}$ to switch between left and
right-invariant Haar measures on $G$.  If $G$ is compact, then one can
show that left-invariant Haar measure on $G$ is invariant under right
translations too.  This is trivial when $G$ is commutative, and one
can check that Haar measure on $G$ is invariant under the mapping $x
\mapsto x^{-1}$ when $G$ is compact or commutative.

        Using Haar measure on $G$, one gets a nonnegative linear functional
on the space of continuous real or complex-valued functions with compact
support on $G$ which is invariant under left or right translations,
as appropriate.  This type of linear functional is known as a Haar
integral on $G$, and it is strictly positive in the sense that the
integral of a nonnegative real-valued continuous function with compact
support on $G$ is positive when the function is positive somewhere on $G$.
Of course, any nonnegative linear functional on the space of continuous
functions with compact support on $G$ determines a unique nonnegative
Borel measure with certain regularity properties, by the Riesz representation
theorem.  If such a linear functional is invariant under left or right 
translations and strictly positive in the sense mentioned earlier,
then the corresponding measure is a Haar measure on $G$.  One can also
deal with uniqueness directly in terms of these linear functionals.

        The real line is a locally compact commutative topological group
with respect to addition, and Lebesgue measure on ${\bf R}$ satisfies the 
requirements of Haar measure.  If $p$ is any prime number, then the
$p$-adic numbers ${\bf Q}_p$ form a locally compact commutative topological
group with respect to addition as well, and Haar measure on ${\bf Q}_p$
was discussed in Section \ref{haar measure on Q_p}.  Any group $G$ is a
locally compact topological group with respect to the discrete topology,
with counting measure as Haar measure that is invariant under both left
and right translations.  The unit circle ${\bf T}$ in the complex plane
is a compact commutative topological group with respect to multiplication
and the topology induced on ${\bf T}$ by the standard topology on ${\bf C}$,
and the usual arc-length measure on ${\bf T}$ satisfies the requirements
of Haar measure.  Haar measure on a real Lie group can be given in terms
of a smooth volume form that is invariant under left or right translations,
as appropriate.

\section{Dual groups}
\label{dual groups}

        Let $A$ be a commutative topological group, with the group
operations expressed additively.  A continuous homomorphism from $A$
into the multiplicative group ${\bf T}$ of complex numbers with
modulus equal to $1$ is said to be a
\emph{character}\index{characters} on $A$.  The collection of
characters on $A$ forms a commutative group $\widehat{A}$ with respect
to pointwise multiplication, known as the \emph{dual group}\index{dual
  groups} associated to $A$.  In particular, the identity element in
$\widehat{A}$ is the trivial character on $A$, which is the constant
function equal to $1$ at every point in $A$.  Note that the
multiplicative inverse of $\phi \in \widehat{A}$ is the same as the
complex conjugate of $\phi$.

        Of course, the group ${\bf Z}$ of integers is a commutative
topological group with respect to addition the discrete topology.  
If $z \in {\bf T}$, then
\begin{equation}
\label{j mapsto z^j}
        j \mapsto z^j
\end{equation}
defines a homomorphism from ${\bf Z}$ into ${\bf T}$, and every
homomorphism from ${\bf Z}$ into ${\bf T}$ is of this form.
Similarly,
\begin{equation}
\label{z mapsto z^j}
        z \mapsto z^j
\end{equation}
is a continuous homomorphism from ${\bf T}$ into itself for every
integer $j$, and it is well known that every character on ${\bf T}$ is
of this form.  If $y \in {\bf R}$, then
\begin{equation}
\label{x mapsto exp (i x y)}
        x \mapsto \exp (i \, x \, y)
\end{equation}
is a character on ${\bf R}$, where $\exp$ refers to the complex
exponential function.  It is also well known that every character on
${\bf R}$ is of this form.

        If $A$ is any commutative topological group, then one can
consider $\widehat{A}$ equipped with the topology associated to
uniform convergence on nonempty compact subsets of $A$.  In
particular, it is easy to see that $\widehat{A}$ is also a topological
group with respect to this topology.  This is especially nice when $A$
is locally compact, in which case it can be shown that $\widehat{A}$
is locally compact too.  If $A = {\bf R}$ as a commutative topological
group with respect to addition and the standard topology, then $\widehat{A}$
is isomorphic to ${\bf R}$ as a commutative group, as in the previous
paragraph.  One can check that the dual topology on $\widehat{A}$
corresponds exactly to the standard topology on ${\bf R}$ in this case as
well.

        If $z \in {\bf T}$ has nonnegative real part and $z \ne 1$, then
the real part of $z^j$ is negative for some integer $j$.  This implies
that the only subgroup of ${\bf T}$ consisting of $z \in {\bf T}$ with
nonnegative real part is the trivial subgroup $\{1\}$.  If $A$ is a
commutative topological group and $\phi \in \widehat{A}$ has the
property that the real part of $\phi(x)$ is nonnegative for every $x
\in A$, then it follows that $\phi$ is the trivial character on $A$.
In particular, if $\phi \in \widehat{A}$ satisfies
\begin{equation}
\label{|phi(x) - 1| le 1}
        |\phi(x) - 1| \le 1
\end{equation}
for every $x \in A$, then $\phi$ is the trivial character on $A$.
If $\phi, \psi \in \widehat{A}$ satisfy
\begin{equation}
\label{|phi(x) - psi(x)| le 1}
        |\phi(x) - \psi(x)| \le 1
\end{equation}
for every $x \in A$, then $\phi(x) = \psi(x)$ for every $x \in A$,
since we can apply the previous argument to $\phi / \psi$.
Equivalently, the distance between any two distinct elements of
$\widehat{A}$ with respect to the supremum metric is greater than $1$.
If $A$ is compact, then the topology on $\widehat{A}$ mentioned in the
previous paragraph is the same as the topology determined by the
supremum metric on $\widehat{A}$, and hence $\widehat{A}$ is discrete
with respect to this topology.

\section{Compact commutative groups}
\label{compact commutative groups}

        Let $A$ be a compact commutative topological group, with the
group operations expressed additively.  As in Section \ref{haar
  measure}, there is a unique translation-invariant nonnegative
regular Borel measure $H$ on $A$ that satisfies $H(A) = 1$, which is
the normalized Haar measure on $A$.  Let $L^2(A)$ be the usual space
of complex-valued square-integrable functions on $A$ with respect to
$H$, with the inner product
\begin{equation}
\label{langle f, g rangle = int_A f(x) overline{g(x)} dH(x)}
        \langle f, g \rangle = \int_A f(x) \, \overline{g(x)} \, dH(x).
\end{equation}
If $\phi$ is a character on $A$, then
\begin{equation}
\label{int_A phi(x) dH(x) = ... = phi(a) int_A phi(x) dH(x)}
        \int_A \phi(x) \, dH(x) = \int_A \phi(x + a) \, dH(x) 
                             = \phi(a) \, \int_A \phi(x) \, dH(x)
\end{equation}
for every $a \in A$, using the translation-invariance of $H$ in the
first step.  If $\phi(a) \ne 1$ for some $a \in A$, then it follows
that
\begin{equation}
\label{int_A phi(x) dH(x) = 0}
        \int_A \phi(x) \, dH(x) = 0.
\end{equation}
This implies that
\begin{equation}
\label{langle phi, psi rangle = int_A phi(x) overline{psi(x)} dH(x) = 0}
 \langle \phi, \psi \rangle = \int_A \phi(x) \, \overline{\psi(x)} \, dH(x) = 0
\end{equation}
when $\phi$ and $\psi$ are distinct elements of $\widehat{A}$, by
applying the previous argument to $\phi(x) \, \overline{\psi(x)}$,
which is a nontrivial character on $A$.  The normalization $H(A) = 1$
implies that each element of $\widehat{A}$ has $L^2$ norm equal to
$1$, so that the elements of $\widehat{A}$ are orthonormal in $L^2(A)$.

        Let $C(A)$ be the algebra of complex-valued continuous functions
on $A$, equipped with the supremum norm.  If $\mathcal{E}$ is the
linear span of $\widehat{A}$ in $C(A)$, then it is easy to see that
$\mathcal{E}$ is a sub-algebra of $C(A)$ which is invariant under
complex conjugation and contains the constant functions on $A$.  It is
well known that $\widehat{A}$ separates points in $A$, which implies
that $\mathcal{E}$ separates points in $A$.  It follows that
$\mathcal{E}$ is dense in $C(A)$ with respect to the supremum norm, by
the Stone--Weierstrass theorem.  In particular, $\mathcal{E}$ is dense
in $L^2(A)$, so that $\widehat{A}$ is an orthonormal basis for
$L^2(A)$.

        Similarly, if $E_1$ is a subgroup of $\widehat{A}$, then the
linear span $\mathcal{E}_1$ of $E_1$ in $C(A)$ is a sub-algebra of $C(A)$
that is invariant under complex-conjugation and contains the constant
functions.  If $E_1$ separates points in $A$, then $\mathcal{E}_1$
separates points in $A$ too, and hence $\widehat{E}_1$ is dense in
$C(A)$, by the Stone--Weierstrass theorem again.  If $\phi$ is any
character on $A$ not in $E_1$, then $\phi$ is orthogonal to every element
of $E_1$ with respect to the $L^2$ inner product, which implies that
$\phi$ is orthogonal to every element of $\mathcal{E}_1$.  This implies
that $\phi = 0$ when $\mathcal{E}_1$ is dense in $C(A)$, which is a
contradiction.  It follows that $E_1 = \widehat{A}$ when $E_1$ is a
subgroup of $\widehat{A}$ that separates points in $A$.

        If $\phi \in \widehat{A}$ and the real part of $\phi(x)$ is
nonnegative for each $x \in A$, then $\phi$ is the trivial character
on $A$, as in the previous section.  Alternatively, if $\phi$ is a
nontrivial character on $A$, then the integral of $\phi$ with respect
to $H$ is equal to $0$, as in (\ref{int_A phi(x) dH(x) = 0}).  If the
real part of $\phi(x)$ is nonnegative for every $x \in A$, then it
follows that the real part of $\phi(x)$ is equal to $0$ for every $x
\in A$, contradicting the fact that $\phi(0) = 1$.  One can also use
the orthonormality of characters on $A$ to get that
\begin{equation}
\label{int_A |phi(x) - psi(x)|^2 dH(x) = 2}
        \int_A |\phi(x) - \psi(x)|^2 \, dH(x) = 2
\end{equation}
when $\phi$, $\psi$ are distinct elements of $\widehat{A}$.

\section{Cartesian products}
\label{cartesian products}

        Let $A_1, \ldots, A_n$ be finitely many commutative topological
groups, and consider their Cartesian product $A = \prod_{j = 1}^n
A_j$.  It is easy to see that $A$ is a commutative topological group
as well, where the group operations are defined coordinatewise, and
using the corresponding product topology.  If $\phi_1, \ldots, \phi_n$
are characters on $A_1, \ldots, A_n$, respectively, then
\begin{equation}
\label{phi(x) = prod_{j = 1}^n phi_j(x_j)}
        \phi(x) = \prod_{j = 1}^n \phi_j(x_j)
\end{equation}
defines a character on $A$.  Conversely, one can check that every
character on $A$ is of this form.

        Now let $I$ be an infinite set, and suppose that $A_j$ is a
commutative topological group for each $j \in I$.  As before, $A =
\prod_{j \in I} A_j$ is a commutative topological group with respect
to coordinatewise addition and the product topology.  Let $j_1,
\ldots, j_n$ be finitely many distinct elements of $I$, and let
$\phi_{j_l}$ be a character on $A_{j_l}$ for each $l = 1, \ldots, n$.
Clearly
\begin{equation}
\label{phi(x) = prod_{l = 1}^n phi_{j_l}(x_{j_l})}
        \phi(x) = \prod_{l = 1}^n \phi_{j_l}(x_{j_l})
\end{equation}
defines a character on $A$, where $x_j \in A_j$ denotes the $j$th
coordinate of $x \in A$ for each $j \in I$.  Conversely, suppose that
$\phi$ is a character on $A$.  Thus the set $V$ of $x \in A$ such that
the real part of $\phi(x)$ is positive is an open set in $A$ that
contains $0$.  By definition of the product topology on $A$, there are
open sets $U_j \subseteq A_j$ for each $j \in I$ such that $0 \in U_j$
for each $j$, $U_j = A_j$ for all but finitely many $j \in I$, and
$\prod_{j \in I} U_j \subseteq V$.  Put $B_j = A_j$ when $U_j = A_j$,
and $B_j = \{0\}$ otherwise.  If $B = \prod_{j \in I} B_j$, then $B$
is a subgroup of $A$ contained in $V$, so that the real part of
$\phi(x)$ is positive when $x \in B$.  It follows that $\phi(x) = 1$
for every $x \in B$, as in Section \ref{dual groups}.  This implies
that $\phi(x)$ depends only on the finitely many coordinates $x_j$ of
$x$ such that $B_j = \{0\}$ for each $x \in A$, and hence that $\phi$
can be expressed as in (\ref{phi(x) = prod_{l = 1}^n
  phi_{j_l}(x_{j_l})}), as in the case of finite products.

        Let $A = \prod_{j = 1}^n A_j$ be the product of finitely many
commutative topological groups again.  If $A_j$ is locally compact for
each $j = 1, \ldots, n$, then $A$ is locally compact too, and Haar
measure on $A$ basically corresponds to the product of the Haar
measures on $A_1, \ldots, A_n$.  More precisely, if there is a
countable base for the topology of $A_j$ for each $j$, then one can
use the standard construction of product measures.  Otherwise, one
should use a version of product measures for Borel measures with
suitable regularity properties.  Equivalently, one can get a Haar
integral on $A$ using Haar integrals on the $A_j$'s.  

        If $A = \prod_{j \in I} A_j$ is the product of infinitely many
compact commutative topological groups, then $A$ is also a compact
commutative topological group with respect to the product topology, by
Tychonoff's theorem.  Haar measure on $A$ again basically corresponds
to the product of the Haar measures on the $A_j$'s, normalized so that
the measure of $A_j$ is equal to $1$ for each $j \in I$.  As before,
this is simpler when $I$ is countably infinite, and there is a base
for the topology of $A_j$ with only finitely or countably many elements
for each $j \in I$, which implies that there is a base for the topology
of $A$ with only finitely or countably many elements.  At any rate,
one can look at the Haar integral on $A$, in terms of the Haar integrals
on the $A_j$'s.  Using compactness, one can show that continuous functions
on $A$ can be approximated uniformly by continuous functions on $A$ that
depend on only finitely many coordinates, for which the Haar integral
is much easier to define.

\section{Discrete commutative groups}
\label{discrete commutative groups}

        Let $A$ be a commutative group equipped with the discrete topology,
so that every homomorphism from $A$ into ${\bf T}$ is continuous and
hence a character.  Note that the collection ${\bf T}^A$ of all
mappings from $A$ into ${\bf T}$ is a commutative group with respect
to pointwise multiplication, and that $\widehat{A}$ is a subgroup of
${\bf T}^A$.  More precisely, ${\bf T}^A$ can be considered as a
Cartesian product of copies of ${\bf T}$ indexed by $A$, equipped with
the product topology corresponding to the standard topology on ${\bf
  T}$, and ${\bf T}^A$ is a compact topological group with respect to
this topology, as in the previous section.  One can check that
$\widehat{A}$ is a closed subgroup of ${\bf T}^A$ with respect to the
product topology, so that $\widehat{A}$ becomes a compact topological
group with respect to the induced topology.  This topology on
$\widehat{A}$ is the same as the one mentioned in Section \ref{dual
  groups} in this case, because compact subsets of $A$ are finite when
$A$ is equipped with the discrete topology.

        Similarly, if $E$ is any nonempty set, then the collection 
${\bf T}^E$ of mappings from $E$ into ${\bf T}$ is a compact commutative
topological group with respect to pointwise multiplication and the
product topology that corresponds to the standard topology on ${\bf T}$,
as in the previous paragraph.  If $E \subseteq A$, then there is an
obvious homomorphism from ${\bf T}^A$ onto ${\bf T}^E$, that sends
each mapping from $A$ into ${\bf T}$ to its restriction to $E$.  This
homomorphism is continuous with respect to the corresponding product
topologies, and the restriction of this homomorphism to $\widehat{A}$
is a continuous homomorphism from $\widehat{A}$ into ${\bf T}^E$.
Suppose that $E$ is a set of generators of $A$, in the sense that
every element of $A$ can be expressed as a sum of finitely many
elements of $E$ and their inverses, where elements of $E$ may be
repeated.  Under these conditions, the homomorphism from $\widehat{A}$
into ${\bf T}^E$ just mentioned is a homeomorphism of $\widehat{A}$
onto its image in ${\bf T}^E$, with respect to the topology on the
image of $\widehat{A}$ in ${\bf T}^E$ induced by the product topology
on ${\bf T}^E$.

        Let $B$ be a subgroup of $A$, let $x \in A \backslash B$ be given, 
and let $B(x)$ be the subgroup of $A$ generated by $B$ and $x$.  If
$\phi$ is a homomorphism from $B$ into ${\bf T}$, then it is well known
that $\phi$ can be extended to a homomorphism from $B(x)$ into ${\bf T}$.
If $A$ is generated by $B$ and finitely or countably many other elements
of $A$, then one can repeat the process to get an extension of $\phi$
to a homomorphism from $A$ into ${\bf T}$, and otherwise one can use
Zorn's lemma or the Hausdorff maximality principle.  Using this, one can
show that for each $a \in A$ with $a \ne 0$ there is a homomorphism $\phi$
from $A$ into ${\bf T}$ such that $\phi(a) \ne 1$.  This implies that
characters on $A$ separate points in $A$.

        Put
\begin{equation}
\label{Psi_a(phi) = phi(a)}
        \Psi_a(\phi) = \phi(a)
\end{equation}
for each $a \in A$ and $\phi \in \widehat{A}$, so that $\Psi_a$ maps
$\widehat{A}$ into ${\bf T}$.  More precisely, $\Psi_a$ is a
continuous homomorphism from $\widehat{A}$ into ${\bf T}$, with
respect to the topology on $\widehat{A}$ discussed earlier.  Thus
$\Psi_a$ is an element of the dual $\widehat{\widehat{A}}$ of the dual
$\widehat{A}$ of $A$, and it is easy to see that
\begin{equation}
\label{a mapsto Psi_a}
        a \mapsto \Psi_a
\end{equation}
defines a homomorphism from $A$ into $\widehat{\widehat{A}}$.  Note
that $\widehat{\widehat{A}}$ should be equipped with the discrete
topology, because $\widehat{A}$ is compact, so that (\ref{a mapsto
  Psi_a}) is automatically continuous.  This mapping (\ref{a mapsto
  Psi_a}) is also one-to-one, because $\widehat{A}$ separates points
in $A$, as in the previous paragraph.  Of course, the collection of
elements of $\widehat{\widehat{A}}$ of the form $\Psi_a$ for some $a
\in A$ is a subgroup of $\widehat{\widehat{A}}$.  This subgroup of
$\widehat{\widehat{A}}$ automatically separates points in
$\widehat{A}$, because $\phi \in \widehat{A}$ is not the trivial
character exactly when (\ref{Psi_a(phi) = phi(a)}) is not equal to $1$
for some $a \in A$.  It follows that every element of
$\widehat{\widehat{A}}$ is of the form $\Psi_a$ for some $a \in A$, as
in Section \ref{compact commutative groups}, because $\widehat{A}$ is
compact.

\section{Characters on ${\bf Z}_p$}
\label{characters on Z_p}

        Let us begin with some remarks about cyclic groups.  Let $n$
be a positive integer, let $n \, {\bf Z}$ be the subgroup of ${\bf Z}$
consisting of integer multiples of $n$, and let ${\bf Z} / n \, {\bf Z}$
be the corresponding quotient group, which is a cyclic group of order
$n$.  Also let $w \in {\bf C}$ be an $n$th root of unity, so that $w^n = 1$,
which implies that the modulus of $w$ is equal to $1$.  The mapping
from $j \in {\bf Z}$ to $w^j \in {\bf T}$ is equal to $1$ on $n \, {\bf Z}$,
and hence determines a group homomorphism from ${\bf Z} / n \, {\bf Z}$
into ${\bf T}$.  Every homomorphism from ${\bf Z} / n \, {\bf Z}$ into
${\bf T}$ is of this form, which implies that the dual group associated
to ${\bf Z} / n \, {\bf Z}$ is also a cyclic group of order $n$.

        Now let $p$ be a prime number, and let $\phi$ be a continuous 
homomorphism from ${\bf Z}_p$ as a commutative topological group with
respect to addition into ${\bf T}$.  The continuity of $\phi$ implies
that there is a nonnegative integer $k$ such that the real part of
$\phi(x)$ is positive for every $x \in p^k \, {\bf Z}_p$.  This
implies that $\phi(x) = 1$ for every $x \in p^k \, {\bf Z}_p$, as in
Section \ref{dual groups}, because $p^k \, {\bf Z}_p$ is a subgroup of
${\bf Z}_p$.  It follows that $\phi$ determines a homomorphism from
${\bf Z}_p / p^k \, {\bf Z}_p$ into ${\bf T}$.  We have also seen in
Section \ref{p-adic numbers} that ${\bf Z} / p^k \, {\bf Z}_p$ is
isomorphic as a group to ${\bf Z} / p^k \, {\bf Z}$, so that the
induced homomorphism from ${\bf Z}_p / p^k \, {\bf Z}_p$ into ${\bf
  T}$ is of the form described in the preceding paragraph.
Conversely, every homomorphism from ${\bf Z}_p / p^k \, {\bf Z}_p$
into ${\bf T}$ leads to a homomorphism from ${\bf Z}_p$ into ${\bf
  T}$, by composition with the canonical quotient mapping from ${\bf
  Z}_p$ onto ${\bf Z}_p / p^k \, {\bf Z}_p$.  Any homomorphism from
${\bf Z}_p$ into ${\bf T}$ of this type is automatically continuous,
because $p^k \, {\bf Z}_p$ is an open subgroup of ${\bf Z}_p$ for each
$k \ge 0$.

        Let $n$ be a positive integer again, and consider the space
of complex-valued functions on ${\bf Z} / n \, {\bf Z}$.  This is an
$n$-dimensional vector space, which may be equipped with a
translation-invariant inner product as in Section \ref{compact
  commutative groups}.  Of course, normalized Haar measure on ${\bf Z}
/ n \, {\bf Z}$ is simply the measure that assigns the value $1/n$ to
each element of ${\bf Z} / n \, {\bf Z}$.  As before, characters on
${\bf Z} / n \, {\bf Z}$ are orthonormal with respect to this inner
product.  It follows that the characters on ${\bf Z} / n \, {\bf Z}$
form an orthonormal basis for the space of complex-valued functions on
${\bf Z} / n \, {\bf Z}$, since there are exactly $n$ characters on
${\bf Z} / n \, {\bf Z}$.

        Similarly, characters on ${\bf Z}_p$ are orthonormal with respect
to the $L^2$ inner product associated to normalized Haar measure on
${\bf Z}_p$.  Note that there are $p^k$ characters on ${\bf Z}_p$
obtained from characters on ${\bf Z}_p / p^k \, {\bf Z}_p$ for each
nonnegative integer $k$.  The linear span of these $p^k$ characters
consists of the complex-valued functions on ${\bf Z}_p$ that are
constant on the cosets of $p^k \, {\bf Z}_p$ in ${\bf Z}_p$, which is
a vector space of dimension $p^k$.  The linear span of the set of all
characters on ${\bf Z}_p$ is the space of complex-valued functions on
${\bf Z}_p$ that are constant on the cosets of $p^k \, {\bf Z}_p$ in
${\bf Z}_p$ for some nonnegative integer $k$.  In particular, this implies
that the linear span of the characters on ${\bf Z}_p$ is dense in the
space of all continuous complex-valued functions on ${\bf Z}_p$ with
respect to the supremum norm, and hence in $L^2({\bf Z}_p)$ as well.

\section{The quotient group ${\bf Q}_p / {\bf Z}_p$}
\label{the quotient group Q_p / Z_p}

        Let $p$ be a prime number again, and consider the quotient 
${\bf Q}_p / {\bf Z}_p$ of ${\bf Q}_p$ as a commutative group with
respect to addition by its subgroup ${\bf Z}_p$.  Also let ${\bf
  Z}[1/p]$ be the collection of rational numbers of the form $p^{-j}
\, x$, where $x \in {\bf Z}$, and $j$ is a nonnegative integer.  This
is a dense subgroup of ${\bf Q}_p$ with respect to addition and the
$p$-adic metric, because ${\bf Z}$ is dense in ${\bf Z}_p$, as in
Section \ref{p-adic numbers}.  It follows that the image of ${\bf
  Z}[1/p]$ under the canonical quotient mapping from ${\bf Q}_p$ onto
${\bf Q}_p / {\bf Z}_p$ is all of ${\bf Q}_p / {\bf Z}_p$, since ${\bf
  Z}_p$ is an open subgroup of ${\bf Q}_p$.  Thus we get a
homomorphism from ${\bf Z}[1/p]$ onto ${\bf Q}_p / {\bf Z}_p$ whose
kernel is equal to
\begin{equation}
\label{{bf Z}[1/p] cap {bf Z}_p = {bf Z}}
        {\bf Z}[1/p] \cap {\bf Z}_p = {\bf Z}.
\end{equation}
This leads to a group isomorphism from ${\bf Z}[1/p] / {\bf Z}$ onto
${\bf Q}_p / {\bf Z}_p$.  Because ${\bf Z}_p$ is an open subgroup of
${\bf Q}_p$, we take ${\bf Q}_p / {\bf Z}_p$ to be equipped with the
discrete topology.

        Alternatively, observe that
\begin{equation}
\label{{bf Z}[1/p] = bigcup_{j = 0}^infty p^{-j} {bf Z}}
        {\bf Z}[1/p] = \bigcup_{j = 0}^\infty p^{-j} \, {\bf Z}
\end{equation}
and
\begin{equation}
\label{{bf Q}_p = bigcup_{j = 0}^infty p^{-j} {bf Z}_p}
        {\bf Q}_p = \bigcup_{j = 0}^\infty p^{-j} \, {\bf Z}_p,
\end{equation}
which imply that
\begin{equation}
\label{{bf Z}[1/p] / {bf Z} = bigcup_{j = 0}^infty ((p^{-j} {bf Z}) / {bf Z})}
 {\bf Z}[1/p] / {\bf Z} = \bigcup_{j = 0}^\infty ((p^{-j} \, {\bf Z}) / {\bf Z})
\end{equation}
and
\begin{equation}
\label{{bf Q}_p / {bf Z}_p = bigcup_{j = 0}^infty ((p^{-j} {bf Z}_p)/{bf Z}_p)}
{\bf Q}_p / {\bf Z}_p = \bigcup_{j = 0}^\infty ((p^{-j} \, {\bf Z}_p) / {\bf Z}_p).
\end{equation}
Of course, $(p^{-j} \, {\bf Z}) / {\bf Z}$ is isomorphic as a group to
${\bf Z} / p^j \, {\bf Z}$ for each nonnegative integer $j$.  We also
have that $p^{-j} \, {\bf Z}_p / {\bf Z}_p$ is isomorphic as a group
to ${\bf Z}_p / p^j \, {\bf Z}_p$, which is isomorphic to ${\bf Z} /
p^j \, {\bf Z}$ when $j \ge 0$, as in Section \ref{p-adic numbers}.
The obvious inclusion of $p^{-j} \, {\bf Z}$ in $p^{-j} \, {\bf Z}_p$
leads more directly to a group homomorphism from $(p^{-j} \, {\bf Z})
/ {\bf Z}$ into $(p^{-j} \, {\bf Z}_p) / {\bf Z}_p$, which is actually
an isomorphism, for the usual reasons.  It is easy to see that the
isomorphism from ${\bf Z}[1/p] / {\bf Z}$ onto ${\bf Q}_p / {\bf Z}_p$
described in the previous paragraph sends $(p^{-j} \, {\bf Z}) / {\bf
  Z}$ onto $(p^{-j} \, {\bf Z}_p) / {\bf Z}_p$ in this way for each
nonnegative integer $j$.
  
        We can also consider ${\bf Z}[1/p] \subseteq {\bf Q}$ as a 
subgroup of ${\bf R}$ with respect to addition, so that ${\bf Z}[1/p]
/ {\bf Z}$ can be identified with a subgroup of ${\bf R} / {\bf Z}$.
If $\exp z$ is the complex exponential function, then
\begin{equation}
\label{r mapsto exp (2 pi i r)}
        r \mapsto \exp (2 \, \pi \, i \, r)
\end{equation}
defines a continuous homomorphism from ${\bf R}$ as a commutative
topological group with respect to addition onto ${\bf T}$ as a compact
commutative group with respect to multiplication.  The kernel of this
homomorphism is equal to ${\bf Z}$, which leads to an isomorphism from
${\bf R} / {\bf Z}$ onto ${\bf T}$.  This isomorphism sends ${\bf
  Z}[1/p] / {\bf Z}$ onto the subgroup of ${\bf T}$ consisting of all
roots of unity with order equal to $p^l$ for some nonnegative integer
$l$.  Using the isomorphism between ${\bf Z}[1/p] / {\bf Z}$ and ${\bf
  Q}_p / {\bf Z}_p$ described earlier, we get a homomorphism from
${\bf Q}_p$ into ${\bf T}$ with kernel ${\bf Z}_p$.

        Equivalently,
\begin{equation}
\label{E_p(x') = exp (2 pi i x')}
        E_p(x') = \exp (2 \, \pi \, i \, x')
\end{equation}
defines a homomorphism from ${\bf Z}[1/p]$ as a commutative group with
respect to addition into ${\bf T}$ as a commutative group with respect
to multiplication, with kernel equal to ${\bf Z}$.  If $x \in {\bf
  Q}_p$, then there is an $x' \in {\bf Z}[1/p]$ such that $x - x' \in
{\bf Z}_p$, because ${\bf Z}[1/p]$ is dense in ${\bf Q}_p$, as before.
If $x'' \in {\bf Z}[1/p]$ also satisfies $x - x'' \in {\bf Z}_p$, then
$x' - x'' \in {\bf Z}_p$, and hence $x' - x'' \in {\bf Z}$, as in
(\ref{{bf Z}[1/p] cap {bf Z}_p = {bf Z}}).  This implies that $E_p(x')
= E_p(x'')$, so that we can extend $E_p$ to a mapping from ${\bf Q}_p$
into ${\bf T}$ by putting
\begin{equation}
\label{E_p(x) = E_p(x')}
        E_p(x) = E_p(x')
\end{equation}
when $x \in {\bf Q}_p$, $x' \in {\bf Z}[1/p]$, and $x - x' \in {\bf
  Z}_p$.  If $x, y \in {\bf Q}_p$ and $x', y' \in {\bf Z}[1/p]$
satisfy $x - x', y - y' \in {\bf Z}_p$, then $x' + y' \in {\bf
  Z}[1/p]$ and
\begin{equation}
\label{(x + y) - (x' + y') = (x - x') + (y - y') in {bf Z}_p}
        (x + y) - (x' + y') = (x - x') + (y - y') \in {\bf Z}_p,
\end{equation}
so that
\begin{equation}
\label{E_p(x + y) = E_p(x' + y') = E_p(x') E_p(y') = E_p(x) E_p(y)}
        E_p(x + y) = E_p(x' + y') = E_p(x') \, E_p(y') = E_p(x) \, E_p(y).
\end{equation}
Thus the extension of $E_p$ to ${\bf Q}_p$ is a homomorphism from
${\bf Q}_p$ as a commutative group with respect to addition into ${\bf
  T}$ as a commutative group with respect to multiplication.  It is
easy to see that the kernel of this homomorphism is equal to ${\bf
  Z}_p$, since the kernel of $E_p$ on ${\bf Z}[1/p]$ is equal to ${\bf
  Z}$.

\section{Characters on ${\bf Q}_p$}
\label{characters on Q_p}

        Let $p$ be a prime number, and let $\phi$ be a continuous 
homomorphism from ${\bf Q}_p$ into ${\bf T}$.  Thus the restriction of
$\phi$ to ${\bf Z}_p$ is a continuous homomorphism from ${\bf Z}_p$
into ${\bf T}$, and hence there is a nonnegative integer $k$ such that
$\phi(x) = 1$ for every $x \in p^k \, {\bf Z}_p$, as in Section
\ref{characters on Z_p}.  This implies that $\phi$ determines a
homomorphism from ${\bf Q}_p / p^k \, {\bf Z}_p$ into ${\bf T}$.
Conversely, every homomorphism from ${\bf Q}_p / p^k \, {\bf Z}_p$
into ${\bf T}$ leads to a homomorphism from ${\bf Q}_p$ into ${\bf
  T}$, by composition with the canonical quotient mapping from ${\bf
  Q}_p$ onto ${\bf Q}_p / p^k \, {\bf Z}_p$.  Any homomorphism from
${\bf Q}_p$ into ${\bf T}$ obtained in this way is continuous, because
$p^k \, {\bf Z}_p$ is an open subgroup of ${\bf Q}_p$.

        Let $E_p$ be the group homomorphism from ${\bf Q}_p$ into ${\bf T}$ 
discussed in the previous section, and put
\begin{equation}
\label{phi_y(x) = E_p(x y)}
        \phi_y(x) = E_p(x \, y)
\end{equation}
for each $x, y \in {\bf Q}_p$.  Thus $\phi_y$ is a group homomorphism
from ${\bf Q}_p$ into ${\bf T}$ for each $y \in {\bf Q}_p$, which is
trivial when $y = 0$.  Otherwise, if $y \ne 0$, so that $|y|_p = p^k$
for some integer $k$, then the kernel of $\phi_y$ is equal to $p^k \,
{\bf Z}_p$.  In particular, $\phi_y$ is continuous as a mapping from
${\bf Q}_p$ into ${\bf T}$ for every $y \in {\bf Q}_p$.  This implies
that $\phi_y$ is an element of the dual $\widehat{{\bf Q}_p}$ of ${\bf
  Q}_p$ as a commutative topological group with respect to addition,
and it is easy to see that the mapping from $y \in {\bf Q}_p$ to
$\phi_y \in \widehat{{\bf Q}_p}$ is a group homomorphism.

        Similarly, the restriction of $\phi_y$ to ${\bf Z}_p$ is a 
continuous homomorphism from ${\bf Z}_p$ as a commutative topological
group with respect to addition into ${\bf T}$, and the mapping from $y
\in {\bf Q}_p$ to the restriction of $\phi_y$ to ${\bf Z}_p$ is a
group homomorphism from ${\bf Q}_p$ into the dual $\widehat{{\bf
    Z}_p}$ of ${\bf Z}_p$.  As before, the restriction of $\phi_y$ to
${\bf Z}_p$ is the trivial character on ${\bf Z}_p$ if and only if $y
\in {\bf Z}_p$.  This implies that the mapping from $y \in {\bf Q}_p$
to the restriction of $\phi_y$ to ${\bf Z}_p$ leads to an injective
homomorphism from ${\bf Q}_p / {\bf Z}_p$ into $\widehat{{\bf Z}_p}$.
One can check that every continuous group homomorphism from ${\bf
  Z}_p$ into ${\bf T}$ is equal to the restriction of $\phi_y$ to
${\bf Z}_p$ for some $y \in {\bf Q}_p$, so that $\widehat{{\bf Z}_p}$
is isomorphic to ${\bf Q}_p / {\bf Z}_p$ as a group.  More precisely,
if $\phi$ is a homomorphism from ${\bf Z}_p$ into ${\bf T}$ whose
kernel contains $p^k \, {\bf Z}_p$ for some nonnegative integer $k$,
then $\phi$ is equal to the restriction of $\phi_y$ to ${\bf Z}_p$ for
some $y \in p^{-k} \, {\bf Z}_p$.

        If $y \in {\bf Z}_p$, then the kernel of $\phi_y : {\bf Q}_p \to 
{\bf T}$ contains ${\bf Z}_p$, and hence $\phi_y$ determines a homomorphism
$\psi_y$ from ${\bf Q}_p / {\bf Z}_p$ into ${\bf T}$.  As in the previous
section, we consider ${\bf Q}_p / {\bf Z}_p$ to be equipped with the
discrete topology, so that every homomorphism from ${\bf Q}_p / {\bf Z}_p$
into ${\bf T}$ is automatically continuous.  Thus $\psi_y$ is an element
of the dual $\widehat{({\bf Q}_p / {\bf Z}_p)}$ of ${\bf Q}_p / {\bf Z}_p$
as a commutative topological group with respect to the discrete topology
for each $y \in {\bf Z}_p$.  Note that $\psi_y$ is the trivial character
on ${\bf Q}_p / {\bf Z}_p$ if and only if $\phi_y$ is the trivial character
on ${\bf Q}_p$, which happens if and only if $y = 0$.  It is easy to see
that $y \mapsto \psi_y$ defines a group homomorphism from ${\bf Z}_p$
into $\widehat{({\bf Q}_p / {\bf Z}_p)}$, as usual.

        If $\psi$ is any homomorphism from ${\bf Q}_p / {\bf Z}_p$ into
${\bf T}$, then one can check that there is a $y \in {\bf Z}_p$ such that
$\psi = \psi_y$.  This is the same as saying that if $\phi$ is a
homomorphism from ${\bf Q}_p$ into ${\bf T}$ whose kernel contains ${\bf Z}_p$,
then there is a $y \in {\bf Z}_p$ such that $\phi = \phi_y$.  To see this,
one can start by showing that for each nonnegative integer $j$, there is a
$y_j \in {\bf Z}_p$ such that $\phi = \phi_{y_j}$ on $p^{-j} \, {\bf Z}_p$.
The image of $y_j$ in ${\bf Z}_p / p^j \, {\bf Z}_p$ is uniquely determined
by this property, which implies that $\{y_j\}_{j = 1}^\infty$ is a Cauchy
sequence in ${\bf Z}_p$.  Thus $\{y_j\}_{j = 1}^\infty$ converges to an
element $y$ of ${\bf Z}_p$, by completeness of the $p$-adic metric,
and one can verify that $\phi = \phi_y$.

        It follows that $y \mapsto \psi_y$ is a group isomorphism of
${\bf Z}_p$ onto $\widehat{({\bf Q}_p / {\bf Z}_p)}$.  Because ${\bf Q}_p / 
{\bf Z}_p$ is equipped with discrete topology, $\widehat{({\bf Q}_p / 
{\bf Z}_p)}$ is compact with respect to the usual topology on the dual
group.  One can also check that $y \mapsto \psi_y$ is a homeomorphism
with respect to the topology on ${\bf Z}_p$ determined by the $p$-adic
metric and the usual topology on the dual group $\widehat{({\bf Q}_p / 
{\bf Z}_p)}$.

        If $\phi$ is any continuous homomorphism from ${\bf Q}_p$
into ${\bf T}$, then the kernel of $\phi$ contains $p^k \, {\bf Z}_p$
for some integer $k$.  Under these conditions, there is a $y \in
p^{-k} \, {\bf Z}_p$ such that $\phi = \phi_y$ on ${\bf Q}_p$.  This
follows from the previous discussion when $k = 0$, and otherwise it is
easy to reduce to that case.  This implies that $y \mapsto \phi_y$
defines a group isomorphism from ${\bf Q}_p$ onto its dual.  It is not
too difficult to verify that this mapping is also a homeomorphism
with respect to the topology on ${\bf Q}_p$ defined by the $p$-adic metric
and the corresponding topology on $\widehat{{\bf Q}_p}$.

\chapter{$r$-Adic integers}
\label{r-adic integers}

\section{$r$-Adic absolute values}
\label{r-adic absolute values}

        Let $r = \{r_j\}_{j = 0}^\infty$ be a sequence of positive integers,
with $r_j \ge 2$ for each $j$.  Put
\begin{equation}
\label{R_l = prod_{j = 1}^l r_j}
        R_l = \prod_{j = 1}^l r_j
\end{equation}
for each positive integer $l$, and $R_0 = 1$, so that $\{R_l\}_{l =
  0}^\infty$ is a strictly increasing sequence of positive integers.
Note that $R_{l + 1} \, {\bf Z}$ is a proper subset of $R_l \, {\bf
  Z}$ for each $l \ge 0$, and that $\bigcap_{l = 0}^\infty R_l \, {\bf
  Z} = \{0\}$.  If $a$ is a nonzero integer, then let $l_r(a)$ be the
largest nonnegative integer such that $a \in R_{l_r(a)} \, {\bf Z}$,
and put $l_r(0) = +\infty$.  Equivalently, $l_r(a) + 1$ is the
smallest positive integer such that $a \not\in R_{l_r(a) + 1} \, {\bf
  Z}$ when $a \ne 0$.  It is easy to see that
\begin{equation}
\label{l_r(a + b) ge min(l_r(a), l_r(b))}
        l_r(a + b) \ge \min(l_r(a), l_r(b))
\end{equation}
and
\begin{equation}
\label{l_r(a b) ge max(l_r(a), l_r(b))}
        l_r(a \, b) \ge \max(l_r(a), l_r(b))
\end{equation}
for every $a, b \in {\bf Z}$.  In particular, $l_r(-a) = l_r(a)$ for
each $a \in {\bf Z}$.  If $r$ is a constant sequence, then
\begin{equation}
\label{l_r(a b) ge l_r(a) + l_r(b)}
        l_r(a \, b) \ge l_r(a) + l_r(b)
\end{equation}
for every $a, b \in {\bf Z}$, and equality holds when $r_1$ is a prime
number.

        Let $t = \{t_l\}_{l = 0}^\infty$ be a strictly decreasing sequence
of positive real numbers that converges to $0$, with $t_0 = 1$.  Put
\begin{equation}
\label{|a|_r = t_{l_r(a)}}
        |a|_r = t_{l_r(a)}
\end{equation}
for each nonzero integer $a$, and $|0|_r = 0$, which corresponds to
(\ref{|a|_r = t_{l_r(a)}}) with $t_\infty = 0$.  Let us call $|a|_r$
the \emph{$r$-adic absolute value}\index{r-adic absolute
  value@$r$-adic absolute value} of $a \in {\bf Z}$, although it also
depends on $t$.  If $p$ is a prime number, $r_j = p$ for each $j \ge
1$, and $t_l = p^{-l}$ for each $l \ge 0$, then this reduces to the
usual $p$-adic absolute value on ${\bf Z}$, as in Section \ref{p-adic
  absolute value}.

        Observe that
\begin{equation}
\label{|a + b|_r le max(|a|_r, |b|_r)}
        |a + b|_r \le \max(|a|_r, |b|_r)
\end{equation}
and
\begin{equation}
\label{|a b|_r le min(|a|_r, |b|_r)}
        |a \, b|_r \le \min(|a|_r, |b|_r)
\end{equation}
for every $a, b \in {\bf Z}$, by (\ref{l_r(a + b) ge min(l_r(a),
  l_r(b))}) and (\ref{l_r(a b) ge max(l_r(a), l_r(b))}).  If $r$ is a
constant sequence, and if $t$ is submultiplicative in the sense that
\begin{equation}
\label{t_{k + l} le t_k t_l}
        t_{k + l} \le t_k \, t_l
\end{equation}
for every $k, l \ge 0$, then
\begin{equation}
\label{|a b|_r le |a|_r |b|_r}
        |a \, b|_r \le |a|_r \, |b|_r
\end{equation}
for every $a, b \in {\bf Z}$, by (\ref{l_r(a b) ge l_r(a) + l_r(b)}).
If $p$ is a prime number, $r_j = p$ for each $j \ge 1$, and $t_l =
(t_1)^l$ for each $l \ge 0$, then equality holds in (\ref{|a b|_r le
  |a|_r |b|_r}) for each $a, b \in {\bf Z}$.  Of course, this reduces
to the case of the $p$-adic absolute value when $t_1 = 1/p$, and
otherwise $|a|_r$ would be the same as a positive power of the
$p$-adic absolute value of $a$.

        Put
\begin{equation}
\label{d_r(a, b) = |a - b|_r}
        d_r(a, b) = |a - b|_r
\end{equation}
for every $a, b \in {\bf Z}$, which we shall call the \emph{$r$-adic
  metric}\index{r-adic metric@$r$-adic metric} on ${\bf Z}$, although
it also depends on $t$, as before.  This is symmetric in $a$ and $b$,
because $l_r(-c) = l_r(c)$ for every $c \in {\bf Z}$, and hence
$|-c|_r = |c|_r$.  Using (\ref{|a + b|_r le max(|a|_r, |b|_r)}), we
get that
\begin{equation}
\label{d_r(a, c) le max(d_r(a, b), d_r(b, c))}
        d_r(a, c) \le \max(d_r(a, b), d_r(b, c))
\end{equation}
for every $a, b, c \in {\bf Z}$, so that $d_r(\cdot, \cdot)$ is an
ultrametric on ${\bf Z}$.  As usual, this reduces to the $p$-adic metric
on ${\bf Z}$ when $r_j = p$ for some prime number $p$ and every $j \ge 1$,
and $t_l = p^{-l}$ for each $l \ge 0$.

\section{An embedding}
\label{an embedding}

        Let $r = \{r_j\}_{j = 1}^\infty$ and $R_l$ be as in the previous section,
and consider the Cartesian product
\begin{equation}
\label{X = prod_{l = 1}^infty ({bf Z} / R_l {bf Z})}
        X = \prod_{l = 1}^\infty ({\bf Z} / R_l \, {\bf Z}),
\end{equation}
consisting of the sequences $x = \{x_l\}_{l = 1}^\infty$ with $x_l \in
{\bf Z} / R_l \, {\bf Z}$ for each $l$.  As in Section \ref{abstract
  cantor sets}, this is a compact Hausdorff topological space with
respect to the product topology that corresponds to the discrete
topology on $X_l = {\bf Z} / R_l \, {\bf Z}$ for each $l$.  This is
also a commutative ring with respect to coordinatewise addition and
multiplication, using the standard ring structure on ${\bf Z} / R_l \,
{\bf Z}$ for each $l$.  It is easy to see that the ring operations are
continuous on $X$, so that $X$ is a topological ring.

        Let $q_l$ be the canonical quotient mapping from ${\bf Z}$
onto ${\bf Z} / R_l \, {\bf Z}$ for each $l \ge 1$, which is a ring
homomorphism with kernel $R_l \, {\bf Z}$.  Thus
\begin{equation}
\label{q(a) = {q_l(a)}_{l = 1}^infty}
        q(a) = \{q_l(a)\}_{l = 1}^\infty
\end{equation}
is an element of $X$ for each $a \in {\bf Z}$, so that $q$ defines a
mapping from ${\bf Z}$ into $X$.  This mapping is an injective ring
homomorphism, because $\bigcap_{l = 1}^\infty R_l \, {\bf Z} = \{0\}$.
If $l(x, y)$ is defined for $x, y \in X$ as in Section \ref{abstract
cantor sets}, and $l_r(a)$ is as in the previous section, then
\begin{equation}
\label{l(q(a), q(b)) = l_r(a - b)}
        l(q(a), q(b)) = l_r(a - b)
\end{equation}
for every $a, b \in {\bf Z}$.

        Let $t = \{t_l\}_{l = 0}^\infty$ be a strictly decreasing sequence
of positive real numbers that converges to $0$ and with $t_0 = 1$, as
before.  This leads to an ultrametric $d(x, y)$ on $X$ as in
(\ref{d(x, y) = t_{l(x, y)}}), and to an $r$-adic metric $d_r(a, b)$
on ${\bf Z}$, as in (\ref{d_r(a, b) = |a - b|_r}).  Under these conditions,
\begin{equation}
        d(q(a), q(b)) = d_r(a, b)
\end{equation}
for every $a, b \in {\bf Z}$, using also (\ref{|a|_r = t_{l_r(a)}}).
By construction, $d_r(a, b)$ is invariant under translations on ${\bf
  Z}$, and in fact $d(x, y)$ is invariant under translations on $X$ as
a commutative group with respect to addition as well.

        Note that $X$ is complete with respect to $d(x, y)$, because
$X$ is compact.  One can also check this directly from the definitions.
It follows that the completion of ${\bf Z}$ with respect to $d_r(a, b)$
can be identified with the closure of $q({\bf Z})$ in $X$.  We shall
discuss this further in the next section.

\section{Coherent sequences}
\label{coherent sequences}

        Let us continue with the notation and hypotheses in the 
previous sections.  Observe that there is a natural ring homomorphism
from ${\bf Z} / R_{l + 1} \, {\bf Z}$ onto ${\bf Z} / R_l \, {\bf Z}$
for each $l \ge 1$, because $R_{l + 1} \, {\bf Z} \subseteq R_l \,
{\bf Z}$.  An element $x = \{x_l\}_{l = 1}^\infty$ of $X$ is said to
be a \emph{coherent sequence}\index{coherent sequences} if the image
of $x_{l + 1}$ under the natural homomorphism from ${\bf Z} / R_{l +
  1} \, {\bf Z}$ onto ${\bf Z} / R_l \, {\bf Z}$ is equal to $x_l$
for each $l$.  Let $Y$ be the subset of $X$ consisting of all coherent
sequences, which is a sub-ring of $X$ with respect to termwise
addition and multiplication.  It is easy to see that $Y$ is also
a closed set in $X$ with respect to the product topology, which implies
that $Y$ is compact, since $X$ is compact.

        If $a \in {\bf Z}$, then the image of $q_{l + 1}(a)$ under the
natural mapping from ${\bf Z} / R_{l + 1} \, {\bf Z}$ onto ${\bf Z} /
R_l \, {\bf Z}$ is equal to $q_l(a)$ for each $l$, so that $q(a) =
\{q_l(a)\}_{l = 1}^\infty$ is a coherent sequence.  Thus $q({\bf Z})
\subseteq Y$, and one can check that
\begin{equation}
\label{overline{q({bf Z})} = Y}
        \overline{q({\bf Z})} = Y,
\end{equation}
where $\overline{q({\bf Z})}$ is the closure of $q({\bf Z})$ with
respect to the product topology on $X$.  More precisely, let $x \in Y$
and $n \in {\bf Z}_+$ be given, and let $a$ be an integer such that
$q_n(a) = x_n$.  This implies that $q_l(a) = x_l$ when $l \le n$,
because $x = \{x_l\}_{l = 1}^\infty$ is a coherent sequence.  It
follows that $x \in \overline{q({\bf Z})}$, as desired, since $n$ is
arbitrary.

        The space ${\bf Z}_r$ of \emph{$r$-adic integers}\index{r-adic 
integers@$r$-adic integers} can be obtained initially as a metric space 
by completing ${\bf Z}$ with respect to the $r$-adic metric.  Using
the isometric embedding $q$ of ${\bf Z}$ in $X$, ${\bf Z}_r$ can be
identified with the set $Y$ of coherent sequences in $X$, equipped with
the restriction of the metric $d(x, y)$ on $X$ to $Y$.  This identification
is very convenient for showing that addition and multiplication on ${\bf Z}$
can be extended continuously to ${\bf Z}_r$, so that ${\bf Z}_r$ is a
compact commutative topological group.  Similarly, the $r$-adic absolute
value can be extended to ${\bf Z}_r$, by taking the distance to $0$ in
${\bf Z}_r$, and it satisfies properties like those on ${\bf Z}$.

        If $t' = \{t_l'\}_{l = 0}^\infty$ is another sequence of positive
real numbers that converges to $0$, then we get another $r$-adic
absolute value function $|a|_r'$ on ${\bf Z}$ as in (\ref{|a|_r =
  t_{l_r(a)}}), a corresponding $r$-adic metric $d_r'(a, b)$ on ${\bf
  Z}$ as in (\ref{d_r(a, b) = |a - b|_r}), and another metric $d'(x,
y)$ on $X$ as in (\ref{d(x, y) = t_{l(x, y)}}).  However, the
embedding $q$ of ${\bf Z}$ into $X$ and the set $Y$ of coherent
sequences in $X$ do not depend on $t$, and the metric $d'(x, y)$ also
determines the product topology on $X$.  The identity mapping on ${\bf
  Z}$ is uniformly continuous as a mapping from ${\bf Z}$ equipped
with $d_r(a, b)$ onto ${\bf Z}$ equipped with $d_r'(a, b)$, as well as
in the other direction, and there are analogous statements for the
identity mapping on $X$ and the metrics $d(x, y)$ and $d'(x, y)$.  In
particular, the completion ${\bf Z}_r$ of ${\bf Z}$ does not depend on
the choice of $t$ as a topological ring.

\section{Haar measure on ${\bf Z}_r$}
\label{haar measure on Z_r}

        Let us continue with the same notation and hypotheses as before.
In particular, let us identify the ring ${\bf Z}_r$ of $r$-adic integers
with the set $Y$ of coherent sequences in $X$.  Let $n$ be a positive integer,
and put
\begin{eqnarray}
\label{Y_n = {x = {x_l}_{l = 1}^infty in Y : x_n = 0} = ...}
        Y_n & = & \{x = \{x_l\}_{l = 1}^\infty \in Y : x_n = 0\} \\
 & = & \{x = \{x_l\}_{l = 1}^\infty \in Y : x_l = 0 \hbox{ for each } l \le n\},
                                                               \nonumber
\end{eqnarray}
where the second step uses the fact that $x \in Y$ is a coherent
sequence.  This is a closed set in $X$ with respect to the product
topology, and an ideal in $Y$ as a commutative ring.  This is also a
relatively open set in $Y$, because ${\bf Z} / R_n \, {\bf Z}$ is
equipped with the discrete topology.  It is easy to see that
\begin{equation}
\label{overline{q(R_n {bf Z})} = Y_n}
        \overline{q(R_n \, {\bf Z})} = Y_n,
\end{equation}
for essentially the same reasons as in (\ref{overline{q({bf Z})} =
  Y}).  Let $\pi_n$ be the mapping from $Y$ into ${\bf Z} / R_n \,
{\bf Z}$ defined by
\begin{equation}
\label{pi_n(x) = x_n}
        \pi_n(x) = x_n,
\end{equation}
which is a ring homomorphism from $Y$ into ${\bf Z} / R_n \, {\bf Z}$
whose kernel is equal to $Y_n$.  Of course,
\begin{equation}
\label{pi_n(q(a)) = q_n(a)}
        \pi_n(q(a)) = q_n(a)
\end{equation}
for every $a \in {\bf Z}$, so that $\pi_n$ maps $Y$ onto ${\bf Z} /
R_n \, {\bf Z}$.

        Let $H$ be Haar measure on $Y$, normalized so that $H(Y) = 1$.
Observe that
\begin{equation}
\label{H(Y_n) = 1/R_n}
        H(Y_n) = 1/R_n
\end{equation}
for each positive integer $n$, because $Y$ can be expressed as the
disjoint union of $R_n$ translates of $Y_n$, by the discussion in the
preceding paragraph.  In this situation, it is easy to define the Haar
integral of a continuous real or complex-valued function on $Y$
directly as a limit of Riemann sums, with the measure of $Y_n$ and its
translates equal to $1/R_n$.  This leads to a translation-invariant
regular Borel measure on $Y$, by the Riesz representation theorem,
which is Haar measure on $Y$.

        Let $t = \{t_l\}_{l = 0}^\infty$ be a strictly decreasing sequence
of positive real numbers that converges to $0$ and with $t_0 = 1$, as
in Section \ref{r-adic absolute values}.  This leads to an $r$-adic
absolute value function $|a|_r$ on ${\bf Z}$ as in (\ref{|a|_r =
  t_{l_r(a)}}), an $r$-adic metric $d_r(a, b)$ on ${\bf Z}$ as in
(\ref{d_r(a, b) = |a - b|_r}), and a metric $d(x, y)$ on $X$ as in
(\ref{d(x, y) = t_{l(x, y)}}).  The $r$-adic absolute value function
and metric can be extended to ${\bf Z}_r$ in a natural way, as in the
previous section, and the extension of the $r$-adic metric on ${\bf
  Z}_r$ corresponds exactly to the restriction of $d(x, y)$ to $Y$.
By construction, these metrics are invariant under translations, and
the diameter of $Y$ is equal to $t_0 = 1$.  Similarly, the
diameter of $Y_n$ is equal to $t_n$ for each positive integer $n$.

        Let $H^\alpha_{con}(E)$, $H^\alpha_\delta(E)$, and $H^\alpha(E)$
be defined for $\alpha \ge 0$, $0 < \delta \le \infty$, and $E \subseteq
Y$ as in Chapter \ref{hausdorff measures}, using the restriction of
$d(x, y)$ to $Y$.  Because this is an ultrametric on $Y$, one may as
well use coverings of $E \subseteq Y$ by closed balls in $Y$ in the
definitions of $H^\alpha_{con}(E)$ and $H^\alpha_\delta(E)$, as in Section
\ref{some special cases}.  More precisely, one should consider the empty
set as a closed ball in $Y$ when $\alpha = 0$, but we are mostly interested
in $\alpha > 0$ here.  Otherwise, the closed balls in $Y$ are $Y$ itself
and the translates of $Y_n$ for each positive integer $n$.  

        Let us now restrict our attention for the rest of this section
to the case where $\alpha = 1$ and
\begin{equation}
\label{t_l = 1/R_l}
        t_l = 1/R_l
\end{equation}
for each $l \ge 0$, which satisfies the usual conditions on $t$.
Remember that $H^1_{con}(E) \le H^1_\delta(E)$ for every $E \subseteq
Y$ and $\delta > 0$, by construction.  In the present situation, one
can check that
\begin{equation}
\label{H^1_delta(E) = H^1_{con}(E), 2}
        H^1_\delta(E) = H^1_{con}(E)
\end{equation}
for every $E \subseteq Y$ and $\delta > 0$, as in Section \ref{some
  special cases}.  This uses the fact that $Y_n$ can be expressed as
the union of $R_k / R_n$ translates of $Y_k$ when $k \ge n$.  It
follows that
\begin{equation}
\label{H^1(E) = H^1_{con}(E), 2}
        H^1(E) = H^1_{con}(E)
\end{equation}
for every $E \subseteq Y$ under these conditions, as before.

        In particular,
\begin{equation}
\label{H^1(Y) = H^1_{con}(Y) le diam Y = 1}
        H^1(Y) = H^1_{con}(Y) \le \diam Y = 1.
\end{equation}
In order to show that
\begin{equation}
\label{H^1(Y) = 1}
        H^1(Y) = 1,
\end{equation}
it suffices to verify that $H^1_\delta(Y) \ge 1$ for every $\delta >
0$.  Because $Y$ is compact, one might as well consider only coverings
of $Y$ by finitely many sets in the definition of $H^1_\delta(Y)$, as
in Section \ref{restricting the diameters}.  In fact, it is enough to
consider only coverings of $Y$ by closed balls, as mentioned earlier.
One can then reduce to coverings of $Y$ by finitely many balls of the
same diameter, by subdividing the balls in a covering of $Y$ as necessary.
Thus one gets either a covering of $Y$ by itself, or by finitely many
translates of $Y_k$ for some $k \ge 1$.  The first case is trivial,
and in the second case, we have that $Y$ cannot be covered by fewer than
$R_k$ translates of $Y_k$.  This implies that $H^1_\delta(Y) \ge 1$
for every $\delta > 0$, and hence that (\ref{H^1(Y) = 1}) holds.

        Similarly,
\begin{equation}
\label{H^1(Y_n) = 1/R_n}
        H^1(Y_n) = 1/R_n
\end{equation}
for each $n \ge 1$, which implies that $H^1(U) > 0$ when $U$ is a
nonempty open subset of $Y$.  Of course, any Hausdorff measure on $Y$
with respect to a translation-invariant metric on $Y$ is invariant
under translations as well.

\section{Some related groups}
\label{some related groups}

        If $k$ is a positive integer, then let $k^{-1} \, {\bf Z}$ be the
set of integer multiples of $1/k$, which is a subgroup of the group
${\bf Q}$ of rational numbers with respect to addition.  Note that
${\bf Z} \subseteq k^{-1} \, {\bf Z}$, so that the quotient group
\begin{equation}
\label{(k^{-1} {bf Z}) / {bf Z}}
        (k^{-1} \, {\bf Z}) / {\bf Z}
\end{equation}
can be defined in the usual way.  Of course, (\ref{(k^{-1} {bf Z}) /
  {bf Z}}) is isomorphic to ${\bf Z} / k \, {\bf Z}$.

        Let $r = \{r_j\}_{j = 1}^\infty$ and $R_l$ be as in Section 
\ref{r-adic absolute values}, and observe that
\begin{equation}
\label{R_l^{-1} {bf Z} subseteq R_{l + 1}^{-1} {bf Z}}
        R_l^{-1} \, {\bf Z} \subseteq R_{l + 1}^{-1} \, {\bf Z}
\end{equation}
for each $l \ge 0$.  Thus
\begin{equation}
\label{bigcup_{l = 0}^infty R_l^{-1} {bf Z}}
        \bigcup_{l = 0}^\infty R_l^{-1} \, {\bf Z}
\end{equation}
is also a subgroup of ${\bf Q}$ that contains ${\bf Z}$, so that the
quotient group
\begin{equation}
\label{(bigcup_{l = 0}^infty R_l^{-1} {bf Z}) / {bf Z}}
        \Big(\bigcup_{l = 0}^\infty R_l^{-1} \, {\bf Z}\Big) / {\bf Z}
\end{equation}
is defined.  If we consider $(R_l^{-1} \, {\bf Z} / {\bf Z})$ as a
subgroup of (\ref{(bigcup_{l = 0}^infty R_l^{-1} {bf Z}) / {bf Z}}),
then
\begin{equation}
\label{(R_l^{-1} {bf Z}) / {bf Z} subseteq (R_{l + 1}^{-1} {bf Z}) / {bf Z}}
(R_l^{-1} \, {\bf Z}) / {\bf Z} \subseteq (R_{l + 1}^{-1} \, {\bf Z}) / {\bf Z}
\end{equation}
for every $l$, because of (\ref{R_l^{-1} {bf Z} subseteq R_{l +
    1}^{-1} {bf Z}}), and (\ref{(bigcup_{l = 0}^infty R_l^{-1} {bf Z})
  / {bf Z}}) is the same as
\begin{equation}
\label{bigcup_{l = 0}^infty (R_l^{-1} {bf Z}) / {bf Z}}
        \bigcup_{l = 0}^\infty (R_l^{-1} \, {\bf Z}) / {\bf Z}.
\end{equation}
If $p$ is a prime number, and $r_j = p$ for each $j$, then
(\ref{bigcup_{l = 0}^infty R_l^{-1} {bf Z}}) is the same as ${\bf
  Z}[1/p]$, as in Section \ref{the quotient group Q_p / Z_p}.

        Remember that $\exp (2 \pi i w)$ defines a continuous homomorphism
from ${\bf R}$ as a commutative topological group with respect to
addition onto ${\bf T}$ with kernel ${\bf Z}$, which leads to an
isomorphism from ${\bf R} / {\bf Z}$ onto ${\bf T}$.  The image of
(\ref{(bigcup_{l = 0}^infty R_l^{-1} {bf Z}) / {bf Z}}) under this
isomorphism consists of the $z \in {\bf T}$ such that
\begin{equation}
\label{z^{R_l} = 1}
        z^{R_l} = 1
\end{equation}
for some $l \ge 0$.  In fact, every homomorphism from (\ref{(bigcup_{l
    = 0}^infty R_l^{-1} {bf Z}) / {bf Z}}) into ${\bf T}$ takes values
in this subgroup of ${\bf T}$.  It follows that homomorphisms from
(\ref{(bigcup_{l = 0}^infty R_l^{-1} {bf Z}) / {bf Z}}) into ${\bf T}$
correspond exactly to homomorphisms from (\ref{(bigcup_{l = 0}^infty
  R_l^{-1} {bf Z}) / {bf Z}}) into itself composed with the embedding
of (\ref{(bigcup_{l = 0}^infty R_l^{-1} {bf Z}) / {bf Z}}) into ${\bf
  T}$ obtained from the complex exponential function.

        Let $\theta$ be a homomorphism from (\ref{(bigcup_{l = 0}^infty 
R_l^{-1} {bf Z}) / {bf Z}}) into itself.  Note that $(R_l^{-1} \, {\bf Z}) / 
{\bf Z}$ consists of exactly the elements $a$ of (\ref{(bigcup_{l = 0}^infty 
R_l^{-1} {bf Z}) / {bf Z}}) such that $R_l \cdot a$ is equal to $0$ in 
(\ref{(bigcup_{l = 0}^infty R_l^{-1} {bf Z}) / {bf Z}}).  This implies that
\begin{equation}
\label{theta((R_l^{-1} {bf Z}) / {bf Z}) subseteq (R_l^{-1} {bf Z}) / {bf Z}}
        \theta((R_l^{-1} \, {\bf Z}) / {\bf Z})
                            \subseteq (R_l^{-1} \, {\bf Z}) / {\bf Z}
\end{equation}
for each $l$.  Let $\theta_l$ be the restriction of $\theta$ to
$(R_l^{-1} \, {\bf Z}) / {\bf Z}$ for each $l \ge 1$.

        Because $(R_l^{-1} \, {\bf Z}) / {\bf Z}$ is a cyclic group,
$\theta_l$ is determined by its value at the generator of 
$(R_l^{-1} \, {\bf Z}) / {\bf Z}$ for each $l$.  This permits $\theta_l$
to be expressed in terms of multiplication by an integer.  This integer
is determined by $\theta_l$ modulo $R_l$, so that $\theta_l$ corresponds
to an element $x_l$ of ${\bf Z} / R_l \, {\bf Z}$.  Conversely, every 
element of ${\bf Z} / R_l \, {\bf Z}$ determines a homomorphism from 
$(R_l^{-1} \, {\bf Z}) / {\bf Z}$ into itself in this way.

        By construction, $\theta_l$ is equal to the restriction of
$\theta_{l + 1}$ to $(R_l^{-1} \, {\bf Z}) / {\bf Z}$ for each $l$.
This means exactly that $x_l$ is the image of $x_{l + 1}$ under the
natural homomorphism from ${\bf Z} / R_{l + 1} \, {\bf Z}$ onto
${\bf Z} / R_l \, {\bf Z}$.  Thus $x = \{x_l\}_{l = 1}^\infty$
is a coherent sequence, which is to say that $x$ is an element of the
group $Y$ discussed in Section \ref{coherent sequences}.  Conversely,
every element of $Y$ leads to a sequence of homomorphisms $\theta_l$
from $(R_l^{-1} \, {\bf Z}) / {\bf Z}$ into itself such that $\theta_l$
is the restriction of $\theta_{l + 1}$ to $(R_l^{-1} \, {\bf Z}) / {\bf Z}$
for each $l$.  This leads in turn to a homomorphism $\theta$ from
(\ref{(bigcup_{l = 0}^infty R_l^{-1} {bf Z}) / {bf Z}}) into itself,
since (\ref{(bigcup_{l = 0}^infty R_l^{-1} {bf Z}) / {bf Z}}) is the same 
as (\ref{bigcup_{l = 0}^infty (R_l^{-1} {bf Z}) / {bf Z}}).

        The collection of homomorphisms from (\ref{(bigcup_{l = 0}^infty 
R_l^{-1} {bf Z}) / {bf Z}}) into itself is a commutative group with respect
to addition.  The discussion in the previous paragraphs determines a
one-to-one correspondence between this group and $Y$, which is a group
isomorphism.  Here we consider (\ref{(bigcup_{l = 0}^infty R_l^{-1}
  {bf Z}) / {bf Z}}) to be equipped with the discrete topology, so that
the corresponding dual group consists of all homomorphisms from 
(\ref{(bigcup_{l = 0}^infty R_l^{-1} {bf Z}) / {bf Z}}) into ${\bf T}$.
Because of the correspondence between homomorphisms from
(\ref{(bigcup_{l = 0}^infty R_l^{-1} {bf Z}) / {bf Z}}) into ${\bf T}$
and homomorphisms from (\ref{(bigcup_{l = 0}^infty R_l^{-1} {bf Z}) /
  {bf Z}}) into itself mentioned earlier, we get an isomorphism
between $Y$ and the dual group associated to (\ref{(bigcup_{l =
    0}^infty R_l^{-1} {bf Z}) / {bf Z}}).  One can check that this
isomorphism is also a homeomorphism with respect to the usual topology
on the dual of (\ref{(bigcup_{l = 0}^infty R_l^{-1} {bf Z}) / {bf Z}})
as a discrete commutative group.

\section{Characters on ${\bf Z}_r$}
\label{characters on Z_r}

        Let us continue with the same notation and hypotheses as before,
and let $\phi$ be a continuous homomorphism from $Y$ as a commutative
topological group with respect to addition into ${\bf T}$.  Thus the
set of $x \in Y$ such that the real part of $\phi(x)$ is positive is
an open set in $Y$ that contains $0$, and hence contains $Y_n$ for
some positive integer $n$.  This implies that $\phi(x) = 1$ for every
$x \in Y_n$, as in Section \ref{dual groups}, because $Y_n$ is a
subgroup of $Y$.  If $\pi_n$ is the homomorphism from $Y$ onto
${\bf Z} / R_n \, {\bf Z}$ in (\ref{pi_n(x) = x_n}), then there is a
homomorphism $\psi$ from ${\bf Z} / R_n \, {\bf Z}$ as a commutative
group with respect to addition into ${\bf T}$ such that
\begin{equation}
\label{phi = psi circ pi_n}
        \phi = \psi \circ \pi_n,
\end{equation}
because the kernel of $\pi_n$ is equal to $Y_n$.  Conversely, if
$\psi$ is a homomorphism from ${\bf Z} / R_n \, {\bf Z}$ as a
commutative group with respect to addition into ${\bf T}$, then
(\ref{phi = psi circ pi_n}) defines a continuous group homomorphism
from $Y$ into ${\bf T}$.

        Alternatively, we have seen in the previous section that $Y$
is isomorphic as a commutative topological group to the dual group
associated to (\ref{(bigcup_{l = 0}^infty R_l^{-1} {bf Z}) / {bf Z}}),
where (\ref{(bigcup_{l = 0}^infty R_l^{-1} {bf Z}) / {bf Z}}) is
equipped with the discrete topology.  As in Section \ref{discrete
  commutative groups}, it follows that each element of
(\ref{(bigcup_{l = 0}^infty R_l^{-1} {bf Z}) / {bf Z}}) determines a
character on $Y$, and in fact that this defines an isomorphism between
(\ref{(bigcup_{l = 0}^infty R_l^{-1} {bf Z}) / {bf Z}}) and the dual of $Y$.
The natural topology on the dual of $Y$ is discrete, because $Y$ is compact,
so that this isomorphism is automatically a homeomorphism.  In this case,
one can also check that the dual of $Y$ is isomorphic to
(\ref{(bigcup_{l = 0}^infty R_l^{-1} {bf Z}) / {bf Z}}) using the remarks
in the previous paragraph, and the descriptions of the group homomorphisms 
from ${\bf Z} / R_n \, {\bf Z}$ into ${\bf T}$ at the beginning of Section
\ref{characters on Z_p}.

        As in Section \ref{compact commutative groups}, characters on
$Y$ are orthonormal with respect to the usual $L^2$ inner product
associated to Haar measure on $Y$.  There are $R_n$ characters on $Y$
of the form (\ref{phi = psi circ pi_n}) for each positive integer $n$,
which are constant on the cosets of $Y_n$ in $Y$.  The linear span
of these characters consists of all functions on $Y$ that are constant
on the cosets of $Y_n$ in $Y$, since every function on ${\bf Z} / R_n \, 
{\bf Z}$ can be expressed as a linear combination of characters on
${\bf Z} / R_n , {\bf Z}$.  The linear span of all characters on $Y$
consists of functions on $Y$ that are constant on the cosets of $Y_n$ in $Y$
for some $n$.

\section{Topological equivalence}
\label{topological equivalence}

        Let $r = \{r_j\}_{j = 1}^\infty$ be as in Section \ref{r-adic 
absolute values}, and let $r' = \{r'_j\}_{j = 1}^\infty$ be another
sequence of integers with $r'_j \ge 2$ for each $j$.  Also let $R_l$
be associated to $r$ as before, and put
\begin{equation}
\label{R'_l = prod_{j = 1}^l r'_j}
        R'_l = \prod_{j = 1}^l r'_j
\end{equation}
when $l \ge 1$, and $R'_0 = 1$.  If for each $l \in {\bf Z}_+$
there is an $n \in {\bf Z}_+$ such that $R'_n$ is an integer multiple
of $R_l$, then we put
\begin{equation}
\label{r prec r'}
        r \prec r'.
\end{equation}
It is easy to see that this relation is reflexive and transitive, and
that (\ref{r prec r'}) holds if and only if every open subset of ${\bf
  Z}$ with respect to the $r$-adic topology is an open set with
respect to the $r'$-adic topology as well.  Similarly, if we put
\begin{equation}
\label{r sim r'}
        r \sim r'
\end{equation}
when $r \prec r'$ and $r'\prec r$, then we get an equivalence relation
on the set of these sequences, which holds exactly when the $r$-adic
and $r'$-adic topologies on ${\bf Z}$ are the same.

        As in Section \ref{an embedding}, consider the Cartesian product
\begin{equation}
\label{X' = prod_{l = 1}^infty ({bf Z} / R'_l {bf Z})}
        X' = \prod_{l = 1}^\infty ({\bf Z} / R'_l \, {\bf Z})
\end{equation}
associated to $r'$, which is a compact commutative topological ring
with respect to coordinatewise addition and multiplication, and using
the product topology corresponding to the discrete topology on ${\bf
  Z} / R'_l \, {\bf Z}$ for each $l$.  Let $q'_l$ be the canonical
quotient mapping from ${\bf Z}$ onto ${\bf Z} / R'_l \, {\bf Z}$ for
each $l$, and put
\begin{equation}
\label{q'(a) = {q'_l(a)}_{l = 1}^infty}
        q'(a) = \{q'_l(a)\}_{l = 1}^\infty
\end{equation}
for each $a \in {\bf Z}$, which defines an injective ring homomorphism
from ${\bf Z}$ into $X'$.  As in Section \ref{coherent sequences},
there is a natural ring homomorphism from ${\bf Z} / R'_{l + 1} \,
{\bf Z}$ onto ${\bf Z} / R_l \, {\bf Z}$ for each $l$, and $x' =
\{x'_l\}_{l = 1}^\infty \in X'$ is said to be a coherent sequence if
$x'_l$ is the image under this homomorphism of $x'_{l + 1}$ for each $l$.
Let $Y'$ be the set of coherent sequences in $X'$, which is a closed
sub-ring of $X'$.  This is the same as the closure of $q'({\bf Z})$
in $X'$, and the topological ring ${\bf Z}_{r'}$ of $r'$-adic integers
can be identified with $Y'$.

        If $r \prec r'$, then there is a natural continuous ring
homomorphism from $Y'$ onto $Y$, defined as follows.  Let $x' \in Y'$
and $l \in {\bf Z}_+$ be given, and remember that there is an $n = n(l)
\in {\bf Z}_+$ such that $R'_n$ is an integer multiple of $R_l$.  Of
course, this implies that $R'_k$ is an integer multiple of $R_l$ for
every $k \in {\bf Z}_+$ with $k \ge n$, and hence that there is a
natural ring homomorphism from ${\bf Z} / R'_k \, {\bf Z}$ onto
${\bf Z} / R_l \, {\bf Z}$ when $k \ge n$.  Let $x_l$ be the image of
$x'_k$ under this homomorphism, which one can check is the same for
all $k \ge n$, because $x'$ is a coherent sequence.  One can also check
that $x = \{x_l\}_{l = 1}^\infty$ is a coherent sequence in $X$,
so that
\begin{equation}
\label{x' mapsto x}
        x' \mapsto x
\end{equation}
leads to a natural mapping from $Y'$ into $Y$.  This mapping is a
continuous ring homomorphism, with respect to the topologies induced
on $Y$ and $Y'$ by the product topologies on $X$ and $X'$,
respectively.  If $a \in {\bf Z}$, then (\ref{x' mapsto x}) sends
$q'(a)$ to $q(a)$, so that (\ref{x' mapsto x}) may be considered as an
extension of the identity mapping on ${\bf Z}$ to a continuous ring
homomorphism from ${\bf Z}_{r'}$ into ${\bf Z}_r$.  Because $Y'$ is
compact, (\ref{x' mapsto x}) maps $Y'$ onto a compact set in $Y$, and
onto a closed set in $Y$ in particular.  This implies that (\ref{x'
  mapsto x}) maps $Y'$ onto $Y$, since $q'({\bf Z})$ is mapped onto
$q({\bf Z})$, which is dense in $Y$.  If $r \sim r'$, then (\ref{x'
  mapsto x}) is an isomorphism from $Y'$ onto $Y$ as topological
rings.

        As in Section \ref{some related groups},
\begin{equation}
\label{bigcup_{l = 0}^infty (R'_l)^{-1} {bf Z}}
        \bigcup_{l = 0}^\infty (R'_l)^{-1} \, {\bf Z}
\end{equation}
is a subgroup of ${\bf Q}$ with respect to addition that contains
${\bf Z}$.  Observe that $r \prec r'$ if and only if the analogous
subgroup (\ref{bigcup_{l = 0}^infty R_l^{-1} {bf Z}}) of ${\bf Q}$
associated to $r$ is contained in (\ref{bigcup_{l = 0}^infty
  (R'_l)^{-1} {bf Z}}), and that $r \sim r'$ if and only if
(\ref{bigcup_{l = 0}^infty R_l^{-1} {bf Z}}) is the same as
(\ref{bigcup_{l = 0}^infty (R'_l)^{-1} {bf Z}}).  Similarly, the
quotient group
\begin{equation}
\label{(bigcup_{l = 0}^infty (R'_l)^{-1} {bf Z}) / {bf Z}}
        \Big(\bigcup_{l = 0}^\infty (R'_l)^{-1} \, {\bf Z}\Big) / {\bf Z}
\end{equation}
and its analogue (\ref{(bigcup_{l = 0}^infty R_l^{-1} {bf Z}) / {bf
    Z}}) for $r$ may be considered as subgroups of ${\bf Q} / {\bf
  Z}$.  Clearly (\ref{bigcup_{l = 0}^infty R_l^{-1} {bf Z}}) is
contained in (\ref{bigcup_{l = 0}^infty (R'_l)^{-1} {bf Z}}) if and
only if (\ref{(bigcup_{l = 0}^infty R_l^{-1} {bf Z}) / {bf Z}}) is
contained in (\ref{(bigcup_{l = 0}^infty (R'_l)^{-1} {bf Z}) / {bf
    Z}}), and (\ref{bigcup_{l = 0}^infty R_l^{-1} {bf Z}}) is equal to
(\ref{bigcup_{l = 0}^infty (R'_l)^{-1} {bf Z}}) if and only if
(\ref{(bigcup_{l = 0}^infty R_l^{-1} {bf Z}) / {bf Z}}) is equal to
(\ref{(bigcup_{l = 0}^infty (R'_l)^{-1} {bf Z}) / {bf Z}}).
Thus $r \prec r'$ if and only if (\ref{(bigcup_{l = 0}^infty R_l^{-1}
  {bf Z}) / {bf Z}}) is contained in (\ref{(bigcup_{l = 0}^infty
  (R'_l)^{-1} {bf Z}) / {bf Z}}), and $r \sim r'$ if and only if
(\ref{(bigcup_{l = 0}^infty R_l^{-1} {bf Z}) / {bf Z}}) is the same as
(\ref{(bigcup_{l = 0}^infty (R'_l)^{-1} {bf Z}) / {bf Z}}).

\chapter{Some geometric conditions}
\label{some geometric conditions}

\section{A class of isometries}
\label{a class of isometries}

        Let $X_1, X_2, X_3, \ldots$ be a sequence of sets, each of which
has at least two elements, and let $X = \prod_{j = 1}^\infty X_j$ be
their Cartesian product, as in Section \ref{abstract cantor sets}.
Also let $\{t_l\}_{l = 0}^\infty$ be a strictly decreasing sequence
of positive real numbers, and let $d(x, y)$ be the ultrametric on $X$
defined as in (\ref{d(x, y) = t_{l(x, y)}}).  Thus $x, y \in X$ satisfy
\begin{equation}
\label{d(x, y) le t_k}
        d(x, y) \le t_k
\end{equation}
for some nonnegative integer $k \ge 0$ if and only if $x_j = y_j$ when
$j \le k$.  Suppose that $\phi : X \to X$ is a Lipschitz mapping of
order $1$ with constant $C = 1$ with respect to $d(\cdot, \cdot)$, so
that
\begin{equation}
\label{d(phi(x), phi(y)) le d(x, y)}
        d(\phi(x), \phi(y)) \le d(x, y)
\end{equation}
for every $x, y \in X$.  Let us express $\phi(x)$ as
\begin{equation}
\label{phi(x) = {phi_j(x)}_{j = 1}^infty}
        \phi(x) = \{\phi_j(x)\}_{j = 1}^\infty,
\end{equation}
where $\phi_j : X \to X_j$ for each $j$.  If $x, y \in X$ satisfy $x_j
= y_j$ for $j \le k$ and some $k$, then it follows that
\begin{equation}
\label{d(phi(x), phi(y)) le t_k}
        d(\phi(x), \phi(y)) \le t_k,
\end{equation}
and hence that $\phi_j(x) = \phi_j(y)$ for $j \le k$.  Equivalently,
this means that for each $k \ge 1$, $\phi_k(x)$ only depends on $x_j$
with $j \le k$.  Conversely, if $\phi_k : X \to X_k$ has this property
for each $k \ge 1$, then $\phi : X \to X$ defined as in (\ref{phi(x) =
  {phi_j(x)}_{j = 1}^infty}) satisfies (\ref{d(phi(x), phi(y)) le d(x,
  y)}) for every $x, y \in X$.

        If $x, y \in X$ satisfy
\begin{equation}
\label{x_j = y_j for j < k and x_k ne y_k}
        x_j = y_j \hbox{ for } j < k \hbox{ and } x_k \ne y_k
\end{equation}
for some $k \in {\bf Z}_+$, then $d(x, y) = t_{k - 1}$, by
construction.  Suppose that $\phi_k$ has the property mentioned in the
preceding paragraph for each $k \ge 1$, and that
\begin{equation}
\label{phi_k(x) ne phi_k(y)}
        \phi_k(x) \ne \phi_k(y)
\end{equation}
when $x, y \in X$ satisfy (\ref{x_j = y_j for j < k and x_k ne y_k}).
This implies that
\begin{equation}
\label{d(phi(x), phi(y)) = t_{k - 1}}
        d(\phi(x), \phi(y)) = t_{k - 1}
\end{equation}
when $x, y \in X$ satisfy (\ref{x_j = y_j for j < k and x_k ne y_k}),
because $\phi_j(x) = \phi_j(y)$ when $j < k$.  It follows that $\phi :
X \to X$ is an isometry with respect to $d(\cdot, \cdot)$ under these
conditions.  Conversely, if $\phi : X \to X$ is an isometry with
respect to $d(\cdot, \cdot)$, then one can check that $\phi$ has these
properties.

        In particular, if $\phi_k(x)$ depends only on $x_k$ for each $k$,
then $\phi$ satisfies (\ref{d(phi(x), phi(y)) le d(x, y)}).  In this
case, $\phi : X \to X$ is an isometry with respect to $d(\cdot,
\cdot)$ if and only if $\phi_k(x)$ corresponds to a one-to-one mapping
from $X_k$ into itself for each $k$.  Similarly, if $\phi_k(x)$
corresponds to a mapping from $X_k$ onto itself for each $k$, then
$\phi$ maps $X$ onto itself.

\section{Some isometric equivalences}
\label{some isometric equivalences}

        Let $r = \{r_j\}_{j = 1}^\infty$ be a sequence of positive
integers with $r_j \ge 2$ for each $j$, and let $R_l$ be as in Section
\ref{r-adic absolute values}.  Also let $X$ and $Y$ be as in Sections
\ref{an embedding} and \ref{coherent sequences}, respectively.  Put
\begin{equation}
\label{widetilde{X}_j = {0, 1, ldots, r_j - 1}}
        \widetilde{X}_j = \{0, 1, \ldots, r_j - 1\}
\end{equation}
for each $j \ge 1$, and
\begin{equation}
\label{widetilde{X} = prod_{j = 1}^infty widetilde{X}_j}
        \widetilde{X} = \prod_{j = 1}^\infty \widetilde{X}_j.
\end{equation}
If $t = \{t_l\}_{l = 0}^\infty$ is a strictly decreasing sequence of
positive real numbers, then we get corresponding ultrametrics $d(x,
y)$ on $X$ and $d'(x', y')$ on $\widetilde{X}$, as in Section
\ref{abstract cantor sets}.  As before, a mapping $\psi : Y \to
\widetilde{X}$ corresponds exactly to a sequence of mappings $\psi_j :
Y \to \widetilde{X}_j$, $j \in {\bf Z}_+$, with
\begin{equation}
\label{psi(x) = {psi_j(x)}_{j = 1}^infty}
        \psi(x) = \{\psi_j(x)\}_{j = 1}^\infty
\end{equation}
for each $x \in Y$.

        Suppose that $\psi : Y \to \widetilde{X}$ is Lipschitz of order $1$
with constant $C = 1$ with respect to the restriction of $d(x, y)$ to
$x, y \in Y$ and $d'(x', y')$ on $\widetilde{X}$, so that
\begin{equation}
\label{d'(psi(x), psi(y)) le d(x, y)}
        d'(\psi(x), \psi(y)) \le d(x, y)
\end{equation}
for every $x, y \in Y$.  If $x, y \in Y$ satisfy $x_j = y_j$ for $j
\le k$ and some $k$, then we get that $\psi_j(x) = \psi_j(y)$ when $j
\le k$, as in the previous section.  This is the same as saying that
$\psi_k(x)$ depends only on $x_k$ for each $k \ge 1$, because the
elements of $Y$ are coherent sequences.  Conversely, if $\psi_k(x)$
depends only on $x_k$ for each $k \ge 1$, then $\psi : Y \to
\widetilde{X}$ satisfies (\ref{d'(psi(x), psi(y)) le d(x, y)}).

        Suppose now that $\psi_k(x)$ depends only on $x_k$ for each $k \ge 1$,
and that
\begin{equation}
\label{psi_k(x) ne psi_k(y)}
        \psi_k(x) \ne \psi_k(y)
\end{equation}
when $x, y \in Y$ satisfy
\begin{equation}
\label{x_{k - 1} = y_{k - 1} and x_k ne y_k}
        x_{k - 1} = y_{k - 1} \hbox{ and } x_k \ne y_k.
\end{equation}
If $k = 1$, then we interpret (\ref{x_{k - 1} = y_{k - 1} and x_k ne
  y_k}) as meaning simply that $x_1 \ne y_1$.  Under these conditions,
$\psi : Y \to \widetilde{X}$ is an isometry, for the same reasons as
before.  Conversely, any isometry from $Y$ into $\widetilde{X}$ has
these properties.

        If $\theta_k$ is a mapping from ${\bf Z} / R_k \, {\bf Z}$
into $\widetilde{X}_k$, then
\begin{equation}
\label{psi_k(x) = theta_k(x_k)}
        \psi_k(x) = \theta_k(x_k)
\end{equation}
defines a mapping from $Y$ into $\widetilde{X}_k$.  We would like to
choose such a mapping $\theta_k$ for each $k \in {\bf Z}_+$ so that
the corresponding mapping $\psi_k$ satisfies (\ref{psi_k(x) ne
  psi_k(y)}) for every $x, y \in Y$ for which (\ref{x_{k - 1} = y_{k -
    1} and x_k ne y_k}) holds.  If $k = 1$, then we can use any
one-to-one mapping from ${\bf Z} / R_1 \, {\bf Z}$ onto
$\widetilde{X}_1$, because $R_1 = r_1$ and $\widetilde{X}_1$ has
exactly $r_1$ elements.  Suppose now that $k \ge 2$, and remember that
there is a natural ring homomorphism from ${\bf Z} / R_k \, {\bf Z}$
onto ${\bf Z} / R_{k - 1} , {\bf Z}$, because $R_k \, {\bf Z}
\subseteq R_{k - 1} \, {\bf Z}$.  The kernel of this homomorphism
is equal to
\begin{equation}
\label{R_{k - 1} {bf Z} / R_k {bf Z}}
        R_{k - 1} \, {\bf Z} / R_k \, {\bf Z},
\end{equation}
which has exactly $r_k$ elements.  Of course, ${\bf Z} / R_k \, {\bf
  Z}$ can be partitioned into translates of (\ref{R_{k - 1} {bf Z} /
  R_k {bf Z}}).  The property of $\psi_k$ that we want is equivalent
to saying that the restriction of $\theta_k$ to any translate of
(\ref{R_{k - 1} {bf Z} / R_k {bf Z}}) in ${\bf Z} / R_k \, {\bf Z}$ is
injective.  It is easy to choose $\theta_k$ in this way, because
(\ref{R_{k - 1} {bf Z} / R_k {bf Z}}) has exactly $r_k$ elements,
which is the same as the number of elements of $\widetilde{X}_k$.
This leads to a sequence of mappings $\psi_k : Y \to \widetilde{X}_k$
as in (\ref{psi_k(x) = theta_k(x_k)}), and hence a mapping $\psi : Y
\to \widetilde{X}$ as in (\ref{psi(x) = {psi_j(x)}_{j = 1}^infty}),
which is an isometry.  Note that $\theta_k$ also maps every translate
of (\ref{R_{k - 1} {bf Z} / R_k {bf Z}}) in ${\bf Z} / R_k \, {\bf Z}$
onto $\widetilde{X}_k$, by construction.  Using this, one can check
that $\psi$ maps $Y$ onto $\widetilde{X}$ as well.

\section{Doubling metrics}
\label{doubling metrics}

        A metric $d(x, y)$ on a set $M$ is said to be 
\emph{doubling}\index{doubling metrics} if there is a positive real 
number $C$ such that every open ball in $M$ with radius $r > 0$ can be
covered by $\le C$ open balls of radius $r/2$.  In this case, one
might also say that the metric space $(M, d(x, y))$ is doubling, or
simply that $M$ is doubling, if the choice of the metric is clear.  If
$M$ is doubling, then we can apply the condition repeatedly to get
that every open ball in $M$ with radius $r$ can be covered by $\le
C^k$ open balls of radius $2^{-k} \, r$ for every $k \in {\bf Z}_+$.
In particular, this implies that bounded subsets of $M$ are totally
bounded.  If $M$ is doubling and complete, then it follows that
subsets of $M$ that are both closed and bounded are compact as well.

        If $M$ is doubling, then the iterated condition mentioned in the
previous paragraph implies that every closed ball in $M$ with radius $r > 0$
can be covered by a bounded number of closed balls of radius $r/2$.
Similarly, every subset of $M$ with diameter $\le r$ can be covered
by a bounded number of sets with diameter $\le r/2$.  This implies that
the restriction of $d(x, y)$ to any subset of $M$ is a doubling metric
too.  One can also use the iterated version of the doubling condition
to show that if $M$ is is bilipschitz equivalent to a metric space that
is doubling, then $M$ is doubling.

        The doubling condition can be defined in the same way for 
quasi-metrics, with the same type of properties as before.  If $d(x, y)$
is a quasi-metric on a set $M$ and $a$ is a positive real number, then
we have seen that $d(x, y)^a$ is quasi-metric on $M$ too, as in
Section \ref{snowflake metrics, quasi-metrics}.  It is easy to see that
$d(x, y)$ is doubling if and only if $d(x, y)^a$ is doubling.

        The real line ${\bf R}$ is doubling with respect to the standard
metric, and similarly ${\bf R}^n$ is doubling with respect to the standard
metric for every positive integer $n$.  More precisely, one can use the
invariance of the standard metric under translations and dilations
to reduce the doubling condition to the case of the unit ball, which
then follows from the fact that the unit ball is totally bounded.
Similarly, the set ${\bf Q}_p$ of $p$-adic numbers is doubling with
respect to the $p$-adic metric, for every prime number $p$.

        Let $X_1, X_2, X_3, \ldots$ be a sequence of sets, each of
which has at least two elements, and let $X = \prod_{j = 1}^\infty X_j$
be their Cartesian product.  Also let $\{t_l\}_{l =0}^\infty$ be a strictly
decreasing sequence of positive real numbers, and let $d(x, y)$ be the 
corresponding ultrametric on $X$, as in (\ref{d(x, y) = t_{l(x, y)}}).
If $X_j$ has only finitely many elements for each $j$, and if 
$\{t_l\}_{l = 0}^\infty$ converges to $0$, then $X$ is totally bounded
with respect to $d(x, y)$.  Conversely, one can check that these
conditions are necessary for $X$ to be doubling with respect to $d(x, y)$.

        Of course, $X$ is bounded with respect to $d(x, y)$ by construction.
If $X$ is doubling with respect to $d(x, y)$, then $X$ is totally
bounded in particular, and hence  the number of elements of $X_j$ has to be
finite for each $j \ge 1$, as in the previous paragraph.  In fact, the
number of elements of $X_j$ has to be uniformly bounded in this case.

        Similarly, if $X$ is totally bounded with respect to $d(x, y)$, 
then we have seen that $\{t_l\}_{l = 0}^\infty$ converges to $0$, which 
implies that for each $l \ge 0$, the number of $j \ge l$ such that
$t_j \ge t_l/2$ is finite.  If $X$ is doubling with respect to $d(x, y)$,
then the number of $j \ge l$ such that $t_j \ge t_l/2$ is uniformly bounded
over $l$.  Conversely, if the number of elements of $X_j$ is uniformly
bounded in $j$, and if the number of $j \ge l$ such that $t_j \ge t_l/2$
is uniformly bounded in $l$, then $X$ is doubling with respect to $d(x, y)$.

        Now let $r = \{r_j\}_{j = 1}$ be a sequence of positive integers
with $r_j \ge 2$ for each $j$, and let $t = \{t_l\}_{l = 0}^\infty$
be a strictly decreasing sequence of positive real numbers.  This leads
to an $r$-adic ultrametric on ${\bf Z}$, as in Section \ref{r-adic
absolute values}.  As before, ${\bf Z}$ is totally bounded with respect
to this $r$-adic metric if and only if $\{t_l\}_{l = 0}^\infty$ converges
to $0$.  One can check that ${\bf Z}$ is doubling with respect to this
$r$-adic metric if and only if the $r_j$'s are bounded and the number
of $j \ge l$ such that $t_j \ge t_l/2$ is uniformly bounded in $l$.
Remember that the set ${\bf Z}_r$ of $r$-adic integers is obtained
by completing ${\bf Z}$ as a metric space with respect to the $r$-adic
metric.  Under these same conditions on $r$ and $t$, ${\bf Z}_r$ is
doubling with respect to the corresponding extension of the $r$-adic metric.
This also follows from the earlier discussion of Cartesian products,
using the isometric equivalence described at the end of the preceding
section.

\section{Doubling measures}
\label{doubling measures}

        A nonnegative Borel measure $\mu$ on a metric space $(M, d(x, y))$
is said to be \emph{doubling}\index{doubling measures} if the measure
of every open ball in $M$ is positive and finite, and if there is a
positive real number $C$ such that
\begin{equation}
\label{mu(B(x, 2 r)) le C mu(B(x, r))}
        \mu(B(x, 2 \, r)) \le C \, \mu(B(x, r))
\end{equation}
for every $x \in M$ and $r > 0$.  It is easy to see that Lebesgue
measure on ${\bf R}^n$ is doubling with respect to the standard metric
on ${\bf R}^n$ for each positive integer $n$, and that Haar measure on
${\bf Q}_p$ is doubling with respect to the $p$-adic metric for every
prime number $p$.  Some other examples will be discussed later in the
section.  If $\mu$ satisfies (\ref{mu(B(x, 2 r)) le C mu(B(x, r))}) on
$M$, then
\begin{equation}
\label{mu(B(x, 2^k r)) le C^k mu(B(x, r))}
        \mu(B(x, 2^k \, r)) \le C^k \, \mu(B(x, r))
\end{equation}
for every $x \in M$, $r > 0$, and $k \in {\bf Z}_+$.  Using this, one
can check that if the measure of some open ball in $M$ is positive or
finite with respect to $\mu$, then every open ball in $M$ has the same
property, because of (\ref{mu(B(x, 2 r)) le C mu(B(x, r))}).

        Suppose that $\mu$ is a doubling measure on a metric space
$(M, d(x, y))$, and let $x \in M$ and $r > 0$ be given.  Also let
$y_1, \ldots, y_n$ be finitely many elements of $B(x, r)$ such that
\begin{equation}
\label{d(y_j, y_l) ge r/2}
        d(y_j, y_l) \ge r/2
\end{equation}
when $j \ne l$.  We would like to show that
\begin{equation}
\label{n le C_1}
        n \le C_1
\end{equation}
for some positive real number $C_1$ that depends only on the doubling
constant for $\mu$.  It follows from (\ref{d(y_j, y_l) ge r/2}) and
the triangle inequality that
\begin{equation}
\label{B(y_j, r/4) cap B(y_l, r/4) = emptyset}
        B(y_j, r/4) \cap B(y_l, r/4) = \emptyset
\end{equation}
when $j \ne l$, and hence
\begin{equation}
\label{sum_{j = 1}^n mu(B(y_j, r/4)) = mu(bigcup_{j = 1}^n B(y_j, r/4))}
 \sum_{j = 1}^n \mu(B(y_j, r/4)) = \mu\Big(\bigcup_{j = 1}^n B(y_j, r/4)\Big).
\end{equation}
Using the triangle inequality again, we have that
\begin{equation}
\label{B(y_j, r/4) subseteq B(x, 5r/4)}
        B(y_j, r/4) \subseteq B(x, 5r/4)
\end{equation}
for each $j$, so that
\begin{equation}
\label{mu(bigcup_{j = 1}^n B(y_j, r/4)) le mu(B(x, 5r/4))}
        \mu\Big(\bigcup_{j = 1}^n B(y_j, r/4)\Big) \le \mu(B(x, 5r/4)).
\end{equation}
In the other direction,
\begin{equation}
\label{B(x, 5r/4) subseteq B(y_j, 9r/4)}
        B(x, 5r/4) \subseteq B(y_j, 9r/4)
\end{equation}
for each $j$, because $d(x, y_j) \le r$ by hypothesis.  Thus $\mu(B(x,
5r/4))$ is bounded by a constant times $\mu(B(y_j, r/4)$ for each $j$,
by the doubling condition.  Combining this with (\ref{sum_{j = 1}^n
  mu(B(y_j, r/4)) = mu(bigcup_{j = 1}^n B(y_j, r/4))}) and
(\ref{mu(bigcup_{j = 1}^n B(y_j, r/4)) le mu(B(x, 5r/4))}), we get
(\ref{n le C_1}), as desired.

        Suppose now that $n$ is the largest positive integer for which
there are $n$ elements $y_1, \ldots, y_l$ of $B(x, r)$ satisfying 
(\ref{d(y_j, y_l) ge r/2}).  If $y$ is any element of $B(x, r)$,
then it follows that
\begin{equation}
\label{d(y, y_j) < r/2}
        d(y, y_j) < r/2
\end{equation}
for some $j = 1, \ldots, n$, since otherwise $y_1, \ldots, y_n$
together with $y$ would be $n + 1$ elements of $B(x, r)$ with the same
property.  This shows that
\begin{equation}
\label{B(x, r) subseteq bigcup_{j = 1}^n B(y_j, r/2),}
        B(x, r) \subseteq \bigcup_{j = 1}^n B(y_j, r/2),
\end{equation}
and hence that $M$ is doubling as a metric space, because $x \in M$
and $r > 0$ are arbitrary, and $n$ is uniformly bounded.

        If $d(x, y)$ is an ultrametric on $M$, then the proof of 
(\ref{n le C_1}) can be improved somewhat.  In this case, we can replace
(\ref{B(y_j, r/4) cap B(y_l, r/4) = emptyset}) with
\begin{equation}
\label{B(y_j, r/2) cap B(y_l, r/2) = emptyset}
        B(y_j, r/2) \cap B(y_l, r/2) = \emptyset
\end{equation}
when $j \ne l$.  Of course, we should then use the analogue of
(\ref{sum_{j = 1}^n mu(B(y_j, r/4)) = mu(bigcup_{j = 1}^n B(y_j,
  r/4))}) with $r/4$ replaced by $r/2$.  We also have that
\begin{equation}
\label{B(y_j, r/2) subseteq B(x, r)}
        B(y_j, r/2) \subseteq B(x, r)
\end{equation}
for each $j$, by the ultrametric version of the triangle inequality, so that
\begin{equation}
\label{mu(bigcup_{j = 1}^n B(y_j, r/2)) le mu(B(x, r))}
        \mu\Big(\bigcup_{j = 1}^n B(y_j, r/2)\Big) \le \mu(B(x, r)),
\end{equation}
which is the substitute for (\ref{mu(bigcup_{j = 1}^n B(y_j, r/4)) le
  mu(B(x, 5r/4))}).  Instead of (\ref{B(x, 5r/4) subseteq B(y_j, 9r/4)}),
we can use the fact that
\begin{equation}
\label{B(x, r) subseteq B(y_j, r)}
        B(x, r) \subseteq B(y_j, r)
\end{equation}
for each $j$, by the ultrametric version of the triangle inequality,
and then continue as before.

        Let $X_1, X_2, X_3, \ldots$ be a sequence of finite sets, each of
which has at least two elements, and let $X = \prod_{j = 1}^\infty
X_j$ be their Cartesian product.  Also let $\{t_l\}_{l = 0}^\infty$ be
a strictly decreasing sequence of positive real numbers that converges
to $0$, and let $d(x, y)$ be the corresponding ultrametric on $X$, as
in (\ref{d(x, y) = t_{l(x, y)}}).  As in the previous section, $X$ is
doubling with respect to $d(x, y)$ if and only if the number of
elements of $X_j$ is uniformly bounded in $j$, and the number of $j
\ge l$ such that $t_j \ge t_l/2$ is uniformly bounded in $l$.  Let
$\mu_j$ be a probability measure on $X_j$ for each $j$, where all
subsets of $X_j$ are measurable, and let $\mu$ be the corresponding
product measure on $X$, as in Section \ref{abstract cantor sets}.  If
$\mu$ is a doubling measure on $X$ with respect to $d(x, y)$, then
there is a $c > 0$ such that
\begin{equation}
\label{mu_j({x_j}) ge c}
        \mu_j(\{x_j\}) \ge c
\end{equation}
for every $j \ge 1$ and $x_j \in X_j$.  This implies that $X_j$ has
$\le 1/c$ elements for each $j$, because $\mu_j(X_j) = 1$.
Conversely, if there is a $c > 0$ such that (\ref{mu_j({x_j}) ge c})
holds for every $j \ge 1$ and $x_j \in X_j$, and if the number of $j
\ge l$ such that $t_j \ge t_l/2$ is uniformly bounded in $l$, then one
can check that $\mu$ is a doubling measure on $X$.

        Similarly, let $r = \{r_j\}_{j = 1}^\infty$ be a sequence of 
positive integers with $r_j \ge 2$ for each $j$, and let $t =
\{t_l\}_{l = 0}^\infty$ be a strictly decreasing sequence of positive
real numbers that converges to $0$.  If Haar measure on ${\bf Z}_r$ is
doubling with respect to the $r$-adic metric on ${\bf Z}_r$ associated
to $r$ and $t$, then it is easy to see that the $r_j$'s have to be
uniformly bounded in $j$.  Conversely, if the $r_j$'s are uniformly
bounded in $j$, and if the number of $j \ge l$ such that $t_j \ge
t_l/2$ is uniformly bounded in $l$, then one can check that Haar
measure on ${\bf Z}_r$ is doubling with respect to the $r$-adic metric
associated to $r$ and $t$.  One can also look at this in terms of an
isometric equivalence of ${\bf Z}_r$ with a Cartesian product
$\widetilde{X}$, as in Section \ref{some isometric equivalences}.
More precisely, Haar measure on ${\bf Z}_r$ corresponds to a product
measure $\widetilde{\mu}$ on $\widetilde{X}$ with respect to this
isometric equivalence, using the probability measures
$\widetilde{\mu}_j$ that are uniformly distributed on each factor
$\widetilde{X}_j$ in (\ref{widetilde{X} = prod_{j = 1}^infty
  widetilde{X}_j}), in the sense that $\mu_j(\{x_j\})$ is the same for
each $x_j \in X_j$.

        If $d(\cdot, \cdot)$ is a quasi-metric on $M$, then one can define
the notion of a doubling measure on $M$ in the same way as before, at
least if open balls in $M$ are Borel sets.  In particular, open balls
in $M$ with respect to $d(\cdot, \cdot)$ are open sets when $d(\cdot,
\cdot)$ is continuous with respect to the topology on $M$ that it
determines.  At any rate, this is normally not a problem, and there
are various ways to deal with it.  One can check that the arguments in
this and the next sections have suitable versions for quasi-metrics,
with different constants, as appropriate.

        If $d(\cdot, \cdot)$ and $d'(\cdot, \cdot)$ are quasi-metrics on $M$
such that each is bounded by a constant multiple of the other, then it is
easy to see that doubling measures on $M$ with respect to $d(\cdot, \cdot)$
are the same as doubling measures on $M$ with respect to $d'(\cdot, \cdot)$,
aside from measurability issues as in the previous paragraph.  Similarly,
if $d(\cdot, \cdot)$ is a quasi-metric on $M$ and $a$ is a positive real
number, then $d(\cdot, \cdot)^a$ is a quasi-metric on $M$, as in Section
\ref{snowflake metrics, quasi-metrics}, and doubling measures on $M$
with respect to $d(\cdot, \cdot)$ are the same as doubling measures on $M$
with respect to $d(\cdot, \cdot)^a$, aside from measurability issues
again.  As in Section \ref{snowflake metrics, quasi-metrics}, if 
$d(\cdot, \cdot)$ is a quasi-metric on $M$,it is shown in \cite{m-s-1}
that there is a metric $\widetilde{d}(\cdot, \cdot)$ on $M$ and a positive
real number $a$ such that $d(\cdot, \cdot)$ and $\widetilde{d}(\cdot, \cdot)^a$
are each bounded by constant multiples of the other.  It follows that
doubling measures on $M$ with respect to $d(\cdot, \cdot)$ are the same
as doubling measures with respect to $\widetilde{d}(\cdot, \cdot)$,
aside from the usual measurability issues.

\section{Another doubling condition}
\label{another doubling condition}

        Let $h(r)$ be a monotone increasing nonnegative real-valued
function on the set $[0, +\infty)$ of nonnegative real numbers.  If
there is a nonnegative real number $C$ such that
\begin{equation}
\label{h(2 r) le C h(r)}
        h(2 \, r) \le C \, h(r)
\end{equation}
for every $r \ge 0$, then we say that $h$ satisfies a doubling
condition.  Using the monotonicity of $h$, we can reformulate
(\ref{h(2 r) le C h(r)}) as saying that
\begin{equation}
\label{h(r + t) le h(2 max(r, t)) le C h(max(r, t)) = C max(h(r), h(t))}
 h(r + t) \le h(2 \, \max(r, t)) \le C \, h(\max(r, t)) = C \, \max(h(r), h(t))
\end{equation}
for every $r, t \ge 0$.  As usual, we can also iterate (\ref{h(2 r) le
  C h(r)}), to get that
\begin{equation}
\label{h(2^k r) le C^k h(r)}
        h(2^k \, r) \le C^k \, h(r)
\end{equation}
for every $r \ge 0$ and positive integer $k$.  Note that $h(r) = r^a$
satisfies these conditions with $C = 2^a$ for each $a \ge 0$.

        Let $(M, d(x, y))$ be a metric space, and let $\mu$ be a nonnegative
Borel measure on $M$.  Suppose that the measure of every open ball in $M$ 
with respect to $\mu$ is finite, and put $h_x(0) = 0$ and
\begin{equation}
\label{h_x(r) = mu(B(x, r))}
        h_x(r) = \mu(B(x, r))
\end{equation}
for every $x \in M$ and $r > 0$.  Thus $h_x(r)$ is a monotone
increasing nonnegative real-valued function on $[0, +\infty)$ for each
  $x \in M$, and in fact $h_x(r)$ is also left-continuous at each $r >
  0$ from the left for every $x \in M$, because of the countable
  additivity of $\mu$.  Clearly $\mu$ satisfies the doubling condition
(\ref{mu(B(x, 2 r)) le C mu(B(x, r))}) for some $C \ge 0$ and every
$x \in M$ and $r > 0$ if and only if $h_x(r)$ satisfies the doubling
condition (\ref{h(2 r) le C h(r)}) with the same constant $C$ for
every $x \in M$ and $r \ge 0$.

        As in \cite{mat}, $\mu$ is said to be uniformly distributed on $M$
if (\ref{h_x(r) = mu(B(x, r))}) does not depend on $x$, so that there
is a function $h(r)$ on $[0, +\infty)$ such that $h(0) = 0$ and
\begin{equation}
\label{mu(B(x, r)) = h(r)}
        \mu(B(x, r)) = h(r)
\end{equation}
for every $x \in M$ and $r > 0$.  If $X = \prod_{j = 1}^\infty X_j$
and $\mu$ are as in the previous section, for instance, then $\mu$ has
this property when $\mu_j$ is uniformly distributed on $X_j$ for each
$j$, in the sense that $\mu_j(\{x_j\})$ is the same for each $x_j \in
X_j$.  If there is a transitive group of isometries on $M$ that
preserve $\mu$, then $\mu$ is uniformly distributed on $M$ in the
sense of (\ref{mu(B(x, r)) = h(r)}).  In particular, this includes the
case of Haar measure on a topological group equipped with a
translation-invariant metric.  If $\mu$ is uniformly distributed on
$M$, and if $M$ is doubling as a metric space, then one can check that
$\mu$ is a doubling measure on $M$.

        Now let $d(x, y)$ be a quasi-metric on a set $M$, and let $h(r)$
be a monotone increasing nonnegative real-valued function on $[0, +\infty)$
such that $h(0) = 0$ and $h(r) > 0$ when $r > 0$.  If $h(r)$ also satisfies
a doubling condition as in (\ref{h(2 r) le C h(r)}), then it is easy to see
that $h(d(x, y))$ is a quasi-metric on $M$ as well.  If, in addition,
\begin{equation}
\label{lim_{r to 0+} h(r) = 0}
        \lim_{r \to 0+} h(r) = 0,
\end{equation}
then $h(d(x, y))$ determines the same topology on $M$ as $d(x, y)$,
and indeed they determine the same uniform structure on $M$.  This is
a variant of the situation in Section \ref{subadditive functions}.

\section{Some variants}
\label{some variants}

        Let $\mu$ be a doubling measure on a metric space $(M, d(x, y))$,
and let $x, y \in M$ be given, with $x \ne y$.  Put $t = d(x, y) > 0$,
and observe that
\begin{equation}
\label{B(x, t/2) cap B(y, t/2) = emptyset}
        B(x, t/2) \cap B(y, t/2) = \emptyset
\end{equation}
and
\begin{equation}
\label{B(x, t/2) cup B(y, t/2) subseteq B(x, 3t/2)}
        B(x, t/2) \cup B(y, t/2) \subseteq B(x, 3t/2),
\end{equation}
by the triangle inequality.  Thus
\begin{equation}
\label{mu(B(x, t/2)) + mu(B(y, t/2)) le mu(B(x, 3t/2))}
        \mu(B(x, t/2)) + \mu(B(y, t/2)) \le \mu(B(x, 3t/2)).
\end{equation}
We also have that $B(x, t/2) \subseteq B(y, 3t/2)$, which implies that
\begin{equation}
\label{mu(B(x, t/2)) le mu(B(y, 3t/2))}
        \mu(B(x, t/2)) \le \mu(B(y, 3t/2)).
\end{equation}
Because $\mu$ is a doubling measure on $M$, $\mu(B(y, 3t/2))$ is
bounded by a constant multiple of $\mu(B(y, t/2))$, and hence
$\mu(B(x, t))$ is bounded by a constant multiple of $\mu(B(y, t/2))$.
It follows that there is a positive real number $c_1 < 1$, depending only
on the doubling constant associated to $\mu$, such that
\begin{equation}
\label{mu(B(x, t/2)) le c_1 mu(B(x, 3t/2))}
        \mu(B(x, t/2)) \le c_1 \, \mu(B(x, 3t/2))
\end{equation}
under these conditions.  In particular, if $x$ is a limit point of $M$,
then one can use this to show that $\mu(\{x\}) = 0$.  Similarly, if $M$
is unbounded, then one can check that $\mu(M) = +\infty$.

        Suppose now that $d(\cdot, \cdot)$ is an ultrametric on $M$,
and let $\mu$ be a nonnegative Borel measure on $M$ such that the
measure of every open ball is open and finite.  Instead of the
doubling condition (\ref{mu(B(x, 2 r)) le C mu(B(x, r))}), let us ask that
\begin{equation}
\label{mu(overline{B}(w, r)) le C_2 mu(B(w, r))}
        \mu(\overline{B}(w, r)) \le C_2 \, \mu(B(w, r))
\end{equation}
for some $C_2 \ge 1$ and every $w \in M$ and $r > 0$.  Let $w \in M$
and $r > 0$ be given, and let $z_1, \ldots, z_n$ be finitely many
elements of $\overline{B}(w, r)$ such that
\begin{equation}
\label{d(z_j, z_l) = r}
        d(z_j, z_l) = r
\end{equation}
when $j \ne l$.  Thus the open balls $B(z_j, r)$ are pairwise-disjoint
subsets of $\overline{B}(w, r)$, so that
\begin{equation}
\label{sum_{j = 1}^n mu(B(z_j, r)) le mu(overline{B}(w, r))}
        \sum_{j = 1}^n \mu(B(z_j, r)) \le \mu(\overline{B}(w, r)).
\end{equation}
We also have that
\begin{equation}
\label{overline{B}(z_j, r) = overline{B}(w, r)}
        \overline{B}(z_j, r) = \overline{B}(w, r)
\end{equation}
for each $j$, because $d(w, z_j) \le r$ for each $j$, and using the
ultrametric version of the triangle inequality.  This implies that
\begin{equation}
\label{mu(overline{B}(w, r)) = mu(overline{B}(z_j, r)) le C_2 mu(B(z_j, r))}
 \mu(\overline{B}(w, r)) = \mu(\overline{B}(z_j, r)) \le C_2 \, \mu(B(z_j, r))
\end{equation}
for each $j$, by hypothesis.  Averaging over $j$, we get that
\begin{equation}
\label{mu(overline{B}(w, r)) le ... le frac{C_2}{n} mu(overline{B}(w, r))}
  \mu(\overline{B}(w, r)) \le \frac{C_2}{n} \, \sum_{j = 1}^n \mu(B(z_j, r))
                           \le \frac{C_2}{n} \, \mu(\overline{B}(w, r)),
\end{equation}
using (\ref{sum_{j = 1}^n mu(B(z_j, r)) le mu(overline{B}(w, r))}) in
the second step.  It follows that
\begin{equation}
\label{n le C_2}
        n \le C_2.
\end{equation}
If we take $n$ to be the largest positive integer for which there are
$n$ elements $z_1, \ldots, z_n$ of $\overline{B}(w, r)$ satisfying
(\ref{d(z_j, z_l) = r}) when $j \ne l$, and if $z$ is any element
of $\overline{B}(w, r)$, then $d(z_j, z) < r$ for some $j$, since otherwise
there would be $n + 1$ elements of $\overline{B}(w, r)$ with this property.
This shows that
\begin{equation}
\label{overline{B}(w, r) subseteq bigcup_{j = 1}^n B(z_j, r)}
        \overline{B}(w, r) \subseteq \bigcup_{j = 1}^n B(z_j, r),
\end{equation}
and hence
\begin{equation}
\label{overline{B}(w, r) = bigcup_{j = 1}^n B(z_j, r)}
        \overline{B}(w, r) = \bigcup_{j = 1}^n B(z_j, r),
\end{equation}
because $B(z_j, r) \subseteq \overline{B}(w, r)$ for each $j$, as in
(\ref{overline{B}(z_j, r) = overline{B}(w, r)}).

        Let $x, y \in M$ be given again, with $x \ne y$, and put 
$t = d(x, y) > 0$.  Because $d(\cdot, \cdot)$ is an ulrametric,
we have that
\begin{equation}
\label{B(x, t) cap B(y, t) = emptyset}
        B(x, t) \cap B(y, t) = \emptyset
\end{equation}
and
\begin{equation}
\label{B(x, t) cup B(y, t) subseteq overline{B}(x, t)}
        B(x, t) \cup B(y, t) \subseteq \overline{B}(x, t).
\end{equation}
instead of (\ref{B(x, t/2) cap B(y, t/2) = emptyset}) and (\ref{B(x,
  t/2) cup B(y, t/2) subseteq B(x, 3t/2)}).  This implies that
\begin{equation}
\label{mu(B(x, t)) + mu(B(y, t)) le mu(overline{B}(x, t))}
        \mu(B(x, t)) + \mu(B(y, t)) \le \mu(\overline{B}(x, t)),
\end{equation}
which replaces (\ref{mu(B(x, t/2)) + mu(B(y, t/2)) le mu(B(x,
  3t/2))}).  Under these conditions, $\overline{B}(x, t)$ is the same
as $\overline{B}(y, t)$, so that
\begin{equation}
\label{mu(overline{B}(x, t)) = mu(overline{B}(y, t)) le C_2 mu(B(y, t))}
 \mu(\overline{B}(x, t)) = \mu(\overline{B}(y, t)) \le C_2 \, \mu(B(y, t)),
\end{equation}
by (\ref{mu(overline{B}(w, r)) le C_2 mu(B(w, r))}).  Combining this with
(\ref{mu(B(x, t)) + mu(B(y, t)) le mu(overline{B}(x, t))}), we get that
\begin{eqnarray}
        \mu(B(x, t)) & \le & \mu(\overline{B}(x, t)) - \mu(B(y, t)) \\
       & \le & \mu(\overline{B}(x, t)) - (1/C_2) \, \mu(\overline{B}(x, t))
                                                               \nonumber \\
        & = & (1 - (1/C_2)) \, \mu(\overline{B}(x, t)), \nonumber
\end{eqnarray}
in place of (\ref{mu(B(x, t/2)) le c_1 mu(B(x, 3t/2))}).

        Of course, the condition (\ref{overline{B}(w, r) subseteq 
bigcup_{j = 1}^n B(z_j, r)}) that a closed ball in $M$ with radius $r$
can be covered by a bounded number of open balls of radius $r$ is
weaker than the usual doubling condition for metrics, as in Section
\ref{doubling metrics}.  If $X$ is a Cartesian product as in Section
\ref{abstract cantor sets}, then this condition corresponds to asking
that the number of elements of the $X_j$'s be bounded, without any
additional condition on the sequence $\{t_l\}_{l = 0}^\infty$ used to
define the ultrametric as in (\ref{d(x, y) = t_{l(x, y)}}).
Similarly, (\ref{mu(overline{B}(w, r)) le C_2 mu(B(w, r))}) is weaker
than the doubling condition (\ref{mu(B(x, 2 r)) le C mu(B(x, r))}) in
Section \ref{doubling measures}.  Let $X$ be as in Section \ref{abstract
cantor sets} again, and let $\mu$ be the probability measure on $X$
corresponding to the product of probability measures $\mu_j$ on $X_j$
for each positive integer $j$, as before.  In this case, it is easy to 
see that $\mu$ satisfies (\ref{mu(overline{B}(w, r)) le C_2 mu(B(w, r))})
if and only if
\begin{equation}
\label{mu_j({x_j}) ge 1/C_2}
        \mu_j(\{x_j\}) \ge 1/C_2
\end{equation}
for every $j \ge 1$ and $x_j \in X_j$, without additional conditions
on the $t_l$'s.  If $r = \{r_j\}_{j = 1}^\infty$ is a sequence of
positive integers, with $r_j \ge 2$ for each $j$, then Haar measure on
the group ${\bf Z}_r$ of $r$-adic integers satisfies
(\ref{mu(overline{B}(w, r)) le C_2 mu(B(w, r))}) with respect to an
$r$ if and only if the $r_j$'s are uniformly bounded.  As usual, this
can also be seen in terms of a suitable isometric equivalence with a
Cartesian product, as in Section \ref{some isometric equivalences}.
If $\mu$ is a uniformly distributed Borel measure on an ultrametric
space $M$, and if $M$ satisfies the covering condition
(\ref{overline{B}(w, r) subseteq bigcup_{j = 1}^n B(z_j, r)}) with
(\ref{n le C_2}), then $\mu$ also satisfies (\ref{mu(overline{B}(w,
  r)) le C_2 mu(B(w, r))}), as in the previous section.

\section{Separability}
\label{separability}

        Let $(M, d(x, y))$ be a metric space.  If the metric $d(x, y)$
is doubling, then bounded subsets of $M$ are totally bounded, as in
Section \ref{doubling metrics}.  This implies that $M$ is separable,
by expressing $M$ as a countable union of balls, each of which
is totally bounded and thus has a countable dense subset.  In particular,
if there is a doubling measure $\mu$ on $M$, then $d(x, y)$ is a doubling
metric on $M$, as in Section \ref{doubling measures}, and hence $M$
is separable.  Suppose now that $\mu$ is a nonnegative Borel measure
on $M$ such that every open ball in $M$ has positive finite measure
with respect to $\mu$, and let us check that $M$ is separable.

        Let $x \in M$ and $r, t > 0$ be given, with $t \le r$, and let
$A$ be a subset of $B(x, r)$ such that
\begin{equation}
\label{d(y, z) ge t}
        d(y, z) \ge t
\end{equation}
for every $y, z \in A$ with $y \ne z$.  Thus the balls $B(y, t/2)$ with
$y \in A$ are pairwise disjoint, and
\begin{equation}
\label{B(y, t/2) subseteq B(x, 3 r/ 2)}
        B(y, t/2) \subseteq B(x, 3 r/ 2)
\end{equation}
for each $y \in A$.  If $y_1, \ldots, y_n$ are finitely many elements
of $A$ such that
\begin{equation}
\label{mu(B(y_j, t/2)) ge a}
        \mu(B(y_j, t/2)) \ge a
\end{equation}
for some $a > 0$ and $j = 1, \ldots, n$, then
\begin{equation}
\label{n a le sum_{j = 1}^n mu(B(y_j, t/2)) = ... le mu(B(x, 3 r / 2))}
        n \, a \le \sum_{j = 1}^n \mu(B(y_j, t/2))
   = \mu\Big(\bigcup_{j = 1}^n B(y_j, t/2)\Big) \le \mu(B(x, 3 r / 2)),
\end{equation}
since $\bigcup_{j = 1}^n B(y_j, t/2) \subseteq B(x, 3 r / 2)$, by
(\ref{B(y, t/2) subseteq B(x, 3 r/ 2)}).  This shows that $n$ is
uniformly bounded under these conditions, and hence that there are
only finitely many $y \in A$ such that
\begin{equation}
\label{mu(B(y, t/2)) ge a}
        \mu(B(y, t/2)) \ge a.
\end{equation}
Applying this to a sequence of $a$'s converging to $0$, we get that
$A$ has only finitely or countably many elements.

        Suppose now that $A$ is a maximal subset of $B(x, r)$ such that
(\ref{d(y, z) ge t}) holds for every $y, z \in A$ with $y \ne z$,
which exists by Zorn's lemma or the Hausdorff maximality principle.
If $w$ is any element of $B(x, r)$, then
\begin{equation}
\label{d(w, y) < t}
        d(w, y) < t
\end{equation}
for some $y \in A$, since otherwise $A \cup \{w\}$ would be a larger
set with the same property.  This implies that
\begin{equation}
\label{B(x, r) subseteq bigcup_{y in A} B(y, t)}
        B(x, r) \subseteq \bigcup_{y \in A} B(y, t).
\end{equation}
where $A$ has finitely or countably many elements, as before.  It
follows that $B(x, r)$ has a dense subset with only finitely or
countably many elements, by considering a sequence of $t$'s converging
to $0$.  Thus $M$ is separable under these conditions, since it can be
expressed as the union of a sequence of open balls.

        Alternatively, let $k$ be a positive integer, and let $A_k$
be a subset of $B(x, r)$ that satisfies (\ref{d(y, z) ge t}) for every 
$y, z \in A_k$ with $y \ne z$, and
\begin{equation}
\label{mu(B(y, t/2)) ge 1/k}
        \mu(B(y, t/2)) \ge 1/k
\end{equation}
for every $y \in A_k$.  The earlier argument shows that $A_k$ is a finite
set with a bounded number of elements, and so we suppose now that $A_k$ is
a maximal set with these properties, for each $k \in {\bf Z}_+$.  Thus
\begin{equation}
\label{A = bigcup_{k = 1}^infty A_k}
        A = \bigcup_{k = 1}^\infty A_k
\end{equation}
has only finitely or countably many elements, although this set $A$ does not
normally satisfy (\ref{d(y, z) ge t}) for every $y, z \in A$ with $y \ne z$.

        If $w$ is any element of $B(x, r)$, then
\begin{equation}
\label{mu(B(w, t/2)) ge 1/k}
        \mu(B(w, t/2)) \ge 1/k
\end{equation}
for some $k \in {\bf Z}_+$, because $\mu(B(w, t/2)) > 0$ by
hypothesis.  It follows that (\ref{d(w, y) < t}) holds for some $y \in
A_k$, since otherwise $A_k \cup \{w\}$ would be a larger set with the
same properties as $A_k$.  In particular, (\ref{d(w, y) < t}) holds
for some $y \in A$, so that (\ref{B(x, r) subseteq bigcup_{y in A}
  B(y, t)}) holds again in this situation.  This implies that $B(x,
r)$ has a dense subset with only finitely or countably many elements,
and hence that $M$ is separable, for the same reasons as before.

        Similarly, if there is an $a > 0$ such that (\ref{mu(B(y, t/2)) ge a})
holds for every $y \in B(x, r)$, then the previous argument shows that
the number of elements of a set $A$ as before is bounded.  This
implies that $B(x, r)$ can be covered by finitely many balls of radius
$t$, as in (\ref{B(x, r) subseteq bigcup_{y in A} B(y, t)}).  If for
each $t \in (0, r]$ there is an $a > 0$ with this property, then it
follows that $B(x, r)$ is totally bounded.

        As usual, these arguments can be simplified when $d(\cdot, \cdot)$
is an ultrametric on $M$.  In this case, if $y, z \in B(x, r)$ and
$0 < t \le r$, then either $d(y, z) < t$, and hence $B(y, t) = B(z, t)$,
or $d(y, z) \ge t$, which implies that
\begin{equation}
\label{B(y, t) cap B(z, t) = emptyset}
        B(y, t) \cap B(z, t) = \emptyset.
\end{equation}
Of course, $B(y, t) \subseteq B(y, r) = B(x, r)$ for every $y \in B(x,
r)$ when $t \le r$.  It follows that for each $a> 0$ and $t \in (0,
r]$, there cannot be more than
\begin{equation}
\label{mu(B(x, r)) / a}
        \mu(B(x, r)) / a
\end{equation}
distinct open balls $B(y, t)$ contained in $B(x, r)$ such that
\begin{equation}
\label{mu(B(y, t)) ge a}
        \mu(B(y, t)) \ge a.
\end{equation}
In particular, for each $t \in (0, r]$, there are only finitely or
countably many distinct open balls $B(y, t)$ contained in $B(x, r)$.

\chapter{Maximal functions}
\label{maximal functions}

\section{Definitions}
\label{definitions}

        Let $(X, d(x, y))$ be a metric space, and let $\mu$ be a nonnegative
Borel measure on $X$ such that the measure of any open ball in $X$ is
positive and finite.  If $f$ is a locally integrable function on $X$
with respect to $\mu$, then put
\begin{equation}
\label{M(f)(x) = sup_{B ni x} frac{1}{mu(B)} int_B |f| d mu}
        M(f)(x) = \sup_{B \ni x} \frac{1}{\mu(B)} \int_B |f| \, d\mu
\end{equation}
for each $x \in X$, which may be infinite.  More precisely, the
supremum is taken over all open balls $B = B(y, r)$ in $X$ that
contain $x$ as an element.  This is the uncentered version of the
Hardy--Littlewood maximal function\index{maximal functions} associated
to $f$ with respect to $\mu$ on $X$.  Similarly, if $\nu$ is a
nonnegative Borel measure on $X$, then the corresponding maximal
function is defined by
\begin{equation}
\label{M(nu)(x) = sup_{B ni x} frac{nu(B)}{mu(B)}}
        M(\nu)(x) = \sup_{B \ni x} \frac{\nu(B)}{\mu(B)}
\end{equation}
for each $x \in X$.  Of course, this reduces to (\ref{M(f)(x) = sup_{B
    ni x} frac{1}{mu(B)} int_B |f| d mu}) when $\nu$ is given by
\begin{equation}
\label{nu(A) = int_A |f| d mu}
        \nu(A) = \int_A |f| \, d\mu
\end{equation}
for every Borel set $A \subseteq X$.  If $\nu$ is a real or complex
Borel measure on $X$, then $M(\nu)$ is defined to be the same as
$M(|\nu|)$, where $|\nu|$ is the total variation measure associated to
$\nu$.

        If $f$ and $g$ are locally integrable functions on $X$ with respect
to $\mu$, then
\begin{equation}
\label{M(f + g)(x) le M(f)(x) + M(g)(x)}
        M(f + g)(x) \le M(f)(x) + M(g)(x)
\end{equation}
for every $x \in X$.  Similarly,
\begin{equation}
\label{M(t f)(x) = |t| M(f)(x)}
        M(t \, f)(x) = |t| \, M(f)(x)
\end{equation}
for every $x \in X$ and real or complex number $t$, as appropriate, so
that the mapping from $f$ to $M(f)$ is sublinear.  There are analogous
statements for maximal functions of Borel measures, as in the previous
paragraph.

        Let $\nu$ be a nonnegative Borel measure on $X$, and let $t$
be a nonnegative real number.  If
\begin{equation}
\label{M(nu)(x) > t}
        M(\nu)(x) > t
\end{equation}
for some $x \in X$, then there is an open ball $B$ in $X$ such that $x
\in B$ and
\begin{equation}
\label{frac{nu(B)}{mu(B)} > t}
        \frac{\nu(B)}{\mu(B)} > t,
\end{equation}
by the definition of $M(\nu)$.  Conversely, if $B$ is an open ball in $X$
that satisfies (\ref{frac{nu(B)}{mu(B)} > t}), then $M(\nu)(y) > t$
for every $y \in B$.  Thus
\begin{equation}
\label{V_t = {x in X : M(nu)(x) > t}}
        V_t = \{x \in X : M(\nu)(x) > t\}
\end{equation}
is the same as the union of the open balls $B$ in $X$ that satisfy
(\ref{frac{nu(B)}{mu(B)} > t}).  In particular, (\ref{V_t = {x in X :
    M(nu)(x) > t}}) is an open set in $X$ for each $t \ge 0$.

        If $f$ is a bounded Borel measurable function on $X$, then
\begin{equation}
\label{sup_{x in X} M(f)(x) le ||f||_infty}
        \sup_{x \in X} M(f)(x) \le \|f\|_\infty,
\end{equation}
where $\|f\|_\infty$ denotes the $L^\infty$ norm of $f$ wih respect to
$\mu$.  We shall consider other estimates for maximal functions in the
next sections.

\section{Three covering arguments}
\label{three covering arguments}

        Let $I$, $I'$, and $I''$ be three intervals in the real line,
which may be open, closed, or half-open and half-closed.  If
\begin{equation}
\label{I cap I' cap I'' ne emptyset}
        I \cap I' \cap I'' \ne \emptyset,
\end{equation}
then it is easy to see that one of these interval is contained in the
union of the other two.  Now let $I_1, I_2, \ldots, I_n$ be finitely
many intervals in ${\bf R}$, which may again be open, closed, or
half-open and half-closed.  Using the previous argument repeatedly,
one can find indices $1 \le j_1 < j_2 < \cdots < j_r \le n$ such that
\begin{equation}
\label{bigcup_{l = 1}^r I_{j_l} = bigcup_{k = 1}^n I_k}
        \bigcup_{l = 1}^r I_{j_l} = \bigcup_{k = 1}^n I_k
\end{equation}
and no element of ${\bf R}$ is contained in more than two of the
$I_{j_l}$'s.

        Suppose instead that $d(x, y)$ is an ultrametric on a set $X$, and
let $B_1, \ldots, B_n$ be finitely many distinct balls in $X$ with
respect to $d(x, y)$, which may be open or closed.  In this case,
there are indices $1 \le j_1 < j_2 < \cdots < j_r \le n$ such that
\begin{equation}
\label{bigcup_{l = 1}^r B_{j_l} = bigcup_{k = 1}^n B_k}
        \bigcup_{l = 1}^r B_{j_l} = \bigcup_{k = 1}^n B_k,
\end{equation}
and the balls $B_{j_l}$ are pairwise disjoint.  To see this, one can
take the $B_{j_l}$'s to be maximal among $B_1, \ldots, B_n$ with
respect to inclusion.  This uses the fact that if $B$ and $B'$ are two
open or closed balls in $X$, then either $B \subseteq B'$, $B'
\subseteq B$, or $B \cap B' = \emptyset$.

        Suppose now that $d(x, y)$ is any metric on a set $X$, and let
$B_j = B(x_j, r_j)$ be the open ball in $X$ centered at a point $x_j \in X$
with radius $r_j > 0$ for $j = 1, \ldots, n$.  By rearranging the
indices if necessary, we may also ask that $r_j$ be monotone
decreasing in $j$.  Put $j_1 = 1$, and let $j_2$ be the smallest
integer such that $2 \le j_2 \le n$ and
\begin{equation}
\label{B_{j_1} cap B_{j_2} = emptyset}
        B_{j_1} \cap B_{j_2} = \emptyset,
\end{equation}
if there is one.  Similarly, if $1 = j_1 < j_2 < \cdots < j_l < n$
have been chosen, then let $j_{l + 1}$ be th smallest integer such
that $j_l < j_{l + 1} \le n$ and
\begin{equation}
\label{B_{j_k} cap B_{j_{l + 1}} = emptyset}
        B_{j_k} \cap B_{j_{l + 1}} = \emptyset
\end{equation}
for each $k = 1, \ldots, l$, if there is one.  This process has to
stop in a finite number $r$ of steps, and the corresponding balls
$B_{j_l}$ are pairwise disjoint, by construction.  If an integer $i$,
$1 \le i \le n$, is not equal to $j_l$ for some $l$, then there is an
$l$ such that $j_l < i$ and $B_i \cap B_{j_l} \ne \emptyset$.  This
implies that
\begin{equation}
\label{B_i subseteq B(x_{j_l}, 3 r_{j_l})}
        B_i \subseteq B(x_{j_l}, 3 \, r_{j_l}),
\end{equation}
since the radius $r_i$ of $B_i$ is less than or equal to $r_{j_l}$.
It follows that
\begin{equation}
\label{bigcup_{i = 1}^n B_i subseteq bigcup_{l = 1}^r B(x_{j_l}, 3 r_{j_l})}
 \bigcup_{i = 1}^n B_i \subseteq \bigcup_{l = 1}^r B(x_{j_l}, 3 \, r_{j_l}).
\end{equation}
Essentially the same argument works when the $B_j$'s are closed balls,
or a mixture of open and closed balls.  If $d(x, y)$ is a quasi-metric
on $X$, then the radius $3 \, r_{j_l}$ in (\ref{B_i subseteq
  B(x_{j_l}, 3 r_{j_l})}) and (\ref{bigcup_{i = 1}^n B_i subseteq
  bigcup_{l = 1}^r B(x_{j_l}, 3 r_{j_l})}) should be replaced with
another constant multiple of $r_{j_l}$, depending on the constant in
the quasi-metric condition for $d(x, y)$.

\section{Weak-type estimates}
\label{weak-type estimates}

        Let $(X, d(x, y))$ be a metric space, and let $\mu$ be a nonnegative
Borel measure on $X$ such that the measure of any open ball in $X$ is
positive and finite.  Also let $\nu$ be a nonnegative Borel measure on
$X$ such that $\nu(X) < +\infty$, and let $V_t$ be as in (\ref{V_t =
  {x in X : M(nu)(x) > t}}) for each $t \ge 0$.  Under suitable conditions,
we would like to show that
\begin{equation}
\label{mu(V_t) le C_1 t^{-1} nu(X)}
        \mu(V_t) \le C_1 \, t^{-1} \, \nu(X)
\end{equation}
for some positive real number $C_1$ and every $t > 0$, where $C_1$
does not depend on $\nu$ or $t$.  As in Section \ref{definitions},
$V_t$ is the same as the union of the open balls $B$ in $X$ that
satisfy (\ref{frac{nu(B)}{mu(B)} > t}), for each $t > 0$.  Note that
$X$ is separable, as in Section \ref{separability}, which implies that
there is a base for the topology of $X$ with only finitely or
countably many elements.  It follows that $V_t$ can be expressed as
the union of finitely or countably many open balls $B$ in $X$ that
satisfy (\ref{frac{nu(B)}{mu(B)} > t}) for each $t > 0$, by
Lindel\"of's theorem in topology.  Let $B_1, \ldots, B_n$ be finitely
many distinct open balls in $X$ that satisfy (\ref{frac{nu(B)}{mu(B)}
  > t}) for some $t > 0$, so that
\begin{equation}
\label{mu(B_j) < t^{-1} nu(B_j)}
        \mu(B_j) < t^{-1} \, \nu(B_j)
\end{equation}
for $j = 1, \ldots, n$.  In order to obtain an estimate of the form
(\ref{mu(V_t) le C_1 t^{-1} nu(X)}), it suffices to show that
\begin{equation}
\label{mu(bigcup_{j = 1}^n B_j) le C_1 t^{-1} nu(X)}
        \mu\Big(\bigcup_{j = 1}^n B_j\Big) \le C_1 \, t^{-1} \, \nu(X),
\end{equation}
where $C_1 > 0$ does not depend on $t$, $\nu$, or $B_1, \ldots, B_n$,
and in particular where $C_1$ does not depend on $n$.

        Suppose first that $d(x, y)$ is an ultrametric on $X$.
In this case, there are indices $1 \le j_1 < j_2 < \cdots < j_r \le n$
such that (\ref{bigcup_{l = 1}^r B_{j_l} = bigcup_{k = 1}^n B_k}) holds,
and the balls $B_{j_l}$ are pairwise disjoint, as in the preceding section.
This implies that
\begin{eqnarray}
\label{mu(bigcup_{k = 1}^n B_k) = ... le t^{-1} nu(X)}
 \mu\Big(\bigcup_{k = 1}^n B_k\Big) = \mu\Big(\bigcup_{l = 1}^r B_{j_l}\Big)
                                  & = & \sum_{l = 1}^r \mu(B_{j_l}) \\
 & < & t^{-1} \, \sum_{l = 1}^r \nu(B_{j_l})  \le t^{-1} \, \nu(X), \nonumber
\end{eqnarray}
using (\ref{mu(B_j) < t^{-1} nu(B_j)}) in the first inequality, and
pairwise-disjointness of the $B_{j_l}$'s in the second inequality.
Thus (\ref{mu(bigcup_{j = 1}^n B_j) le C_1 t^{-1} nu(X)}) holds with
$C_1 = 1$, as desired.

        Now let $X$ be the real line with the standard metric, so that
the open balls $B_j$ are open intervals.  As before, there are indices
$1 \le j_1 < j_2 < \cdots < j_r \le n$ such that (\ref{bigcup_{l =
    1}^r B_{j_l} = bigcup_{k = 1}^n B_k}) holds, and no element of $X
= {\bf R}$ is contained in more than two of the $B_{j_l}$'s.  Let
${\bf 1}_A(x)$ be the characteristic or indicator function associated
to a set $A \subseteq X$, which is equal to $1$ when $x \in A$ and to
$0$ when $x \in X \backslash A$.  The condition that no point belong to
more than two of the $B_{j_l}$'s implies that
\begin{equation}
\label{sum_{l = 1}^r {bf 1}_{B_{j_l}}(x) le 2}
        \sum_{l = 1}^r {\bf 1}_{B_{j_l}}(x) \le 2
\end{equation}
for every $x \in X = {\bf R}$.  It follows that
\begin{equation}
\label{sum_{l = 1}^r nu(B_{j_l}) = ... le 2 nu({bf R})}
        \sum_{l = 1}^r \nu(B_{j_l})
 = \int_{\bf R} \Big(\sum_{l = 1}^r {\bf 1}_{B_{j_l}}(x)\Big) \, d\nu(x)
                \le 2 \, \nu({\bf R}).
\end{equation}
Using this, we get that
\begin{eqnarray}
\label{mu(bigcup_{k = 1}^n B_k) = ... le 2 t^{-1} nu({bf R})}
 \mu\Big(\bigcup_{k = 1}^n B_k\Big) = \mu\Big(\bigcup_{l = 1}^r B_{j_l}\Big)
                                 & \le & \sum_{l = 1}^r \mu(B_{j_l}) \\
 & < & t^{-1} \, \sum_{l = 1}^r \nu(B_{j_l}) \le 2 \, t^{-1} \, \nu({\bf R}),
                                                              \nonumber
\end{eqnarray}
as in the previous situation.  This gives (\ref{mu(bigcup_{j = 1}^n
  B_j) le C_1 t^{-1} nu(X)}), with $C_1 = 2$.

        Suppose that $d(x, y)$ is any metric on a set $X$, and that
$B_j = B(x_j, r_j)$ for some $x_j \in X$ and $r_j > 0$, $j = 1, \ldots, n$.
The third argument in the preceding section implies that there are
indices $1 \le j_1 < j_2 < \cdots < j_r \le n$ such that
(\ref{bigcup_{i = 1}^n B_i subseteq bigcup_{l = 1}^r B(x_{j_l}, 3
  r_{j_l})}) holds and the balls $B_{j_l}$ are pairwise disjoint.  If
$\mu$ is a doubling measure on $X$, then it follows that
\begin{eqnarray}
\label{mu(bigcup_{k = 1}^n B_k) le ... le C_1 sum_{l = 1}^r mu(B_{j_l})}
        \mu\Big(\bigcup_{k = 1}^n B_k\Big) 
            \le \mu\Big(\bigcup_{l = 1}^r B(x_{j_l}, 3 \, r_{j_l})\Big) 
            & \le & \sum_{l = 1}^r \mu(B(x_{j_l}, 3 \, r_{j_l})) \\
            & \le & C_1 \, \sum_{l = 1}^r \mu(B_{j_l}), \nonumber
\end{eqnarray}
for a suitable constant $C_1$.  Combining this with (\ref{mu(B_j) <
  t^{-1} nu(B_j)}), we get that
\begin{equation}
\label{mu(bigcup_{k = 1}^n B_k) le ... le C_1 t^{-1} nu(X)}
 \mu\Big(\bigcup_{k = 1}^n B_k\Big) \le C_1 \, \sum_{l = 1}^r \mu(B_{j_l})
 < C_1 \, t^{-1} \, \sum_{l = 1}^r \nu(B_{j_l}) \le C_1 \, t^{-1} \, \nu(X),
\end{equation}
using also the fact that the $B_{j_l}$'s are pairwise disjoint in the
last step.  Thus the same type of estimate holds when $\mu$ is a
doubling measure on any metric space.

\section{Distribution functions}
\label{distribution functions}

        Let $(X, \mathcal{A}, \mu)$ be a measure space, and let $g$ be 
a measurable function on $X$ with values in the set $[0, +\infty]$ of
nonnegative extended real numbers.  The corresponding distribution
function\index{distribution functions} is defined on $[0, +\infty)$ by
\begin{equation}
\label{lambda(t) = mu({x in X : g(x) > t})}
        \lambda(t) = \mu(\{x \in X : g(x) > t\}),
\end{equation}
which is a monotone decreasing function on $[0, +\infty)$ with values
  in $[0, +\infty]$.  Of course, $\lambda(t) \le \mu(X)$ for every
$t \ge 0$, and
\begin{equation}
\label{t^p lambda(t) le ... le int_X g(x)^p d mu(x)}
        t^p \, \lambda(t) \le \int_{\{x \in X : g(x) > t\}} g(x)^p \, d\mu(x)
                           \le \int_X g(x)^p \, d\mu(x)
\end{equation}
for every $p, t > 0$.  In particular, $\lambda(t) < +\infty$ for every
$t > 0$ when $g \in L^p(X)$ for some $p \in (0, +\infty)$.

        Let us suppose from now on in this section that $X$ is at least
$\sigma$-finite with respect to $\mu$.  Note that the set
\begin{equation}
\label{{x in X : g(x) > 0}}
        \{x \in X : g(x) > 0\}
\end{equation}
is measurable and $\sigma$-finite when $\lambda(t) < +\infty$ for each
$t > 0$, so that we could simply replace $X$ with (\ref{{x in X : g(x)
    > 0}}) in this situation, if necessary.  Let us also consider $[0,
  +\infty)$ as a $\sigma$-finite measure space with respect to
  Lebesgue measure, so that $X \times [0, +\infty)$ is a
    $\sigma$-finite measure space as well, with respect to the usual
    product measure construction.  If $X$ is a topological space too,
    then we may consider $X \times [0, +\infty)$ as a topological
      space with respect to the product topology, using the topology
      induced on $[0, +\infty)$ by the standard topology on ${\bf R}$.

        Put
\begin{eqnarray}
\label{U_r = {(x, t) in X times [0, +infty) : t < r < g(x)} = ...}
        U_r & = & \{(x, t) \in X \times [0, +\infty) : t < r < g(x)\} \\
            & = & [0, r) \times \{x \in X : g(x) > r\} \nonumber
\end{eqnarray}
for each $r \in (0, +\infty)$, and
\begin{equation}
\label{U = {(x, t) in X times [0, +infty) : t < g(x)}}
        U = \{(x, t) \in X \times [0, +\infty) : t < g(x)\}.
\end{equation}
Observe that
\begin{equation}
\label{U = bigcup_{r in {bf Q}_+} U_r}
        U = \bigcup_{r \in {\bf Q}_+} U_r,
\end{equation}
where ${\bf Q}_+ = {\bf Q} \cap (0, +\infty)$ is the set of all
positive rational numbers.  This implies that $U$ is a measurable set
in $X \times [0, +\infty)$, since it can be expressed as a countable union
of measurable rectangles.  If $X$ is a topological space, and if
\begin{equation}
\label{{x in X : g(x) > t}}
        \{x \in X : g(x) > t\}
\end{equation}
is an open set in $X$ for each $t \ge 0$, then it is easy to see that
$U$ is an open set in $X \times [0, +\infty)$.  In fact, (\ref{U =
    bigcup_{r in {bf Q}_+} U_r}) shows that $U$ can be expressed as
a countable union of products of open subsets of $X$ and $[0, +\infty)$
in this case.

        Let $p > 0$ be given, and let ${\bf 1}_U(x, t)$ be the indicator
function associated to $U$ on $X \times [0, +\infty)$.  Observe that
\begin{equation}
\label{p t^{p - 1} {bf 1}_U(x, t)}
        p \, t^{p - 1} \, {\bf 1}_U(x, t)
\end{equation}
is a measurable function on $X \times [0, +\infty)$, because $U$ is a
measurable set.  Clearly
\begin{eqnarray}
\label{... = int_X g(x)^p d mu(x)}
\lefteqn{\int_X\Big(\int_{[0, +\infty)} p \, t^{p - 1} \, {\bf 1}_U(x, t)
                                   \, dt\Big) \, d\mu(x)} \\
 & = & \int_X \Big(\int_0^{g(x)} p \, t^{p - 1} \, dt\Big) \, d\mu(x) 
                                                            \nonumber \\
 & = & \int_X g(x)^p \, d\mu(x), \nonumber
\end{eqnarray}
by elementary calculus, and
\begin{eqnarray}
\label{int_0^infty p t^{p - 1} lambda(t) dt}
\lefteqn{\int_{[0, +\infty)} \Big(\int_X p \, t^{p - 1} \, {\bf 1}_U(x, t) 
                                                 \, d\mu(x)\Big) \, dt} \\
 & = & \int_0^\infty p \, t^{p - 1} \, \Big(\int_{\{x \in X : g(x) > t\}}
                                     \, d\mu(x)\Big) \, dt \nonumber \\
 & = & \int_0^\infty p \, t^{p - 1} \, \lambda(t) \, dt. \nonumber
\end{eqnarray}
Remember that a monotone function on an interval in the real line
continuous at all but at most finitely or countably many points, which
simplifies questions of measurability and integrability.  At any rate,
it follows from (\ref{... = int_X g(x)^p d mu(x)}) and
(\ref{int_0^infty p t^{p - 1} lambda(t) dt}) that
\begin{equation}
\label{int_X g(x)^p d mu(x) = int_0^infty p t^{p - 1} lambda(t) dt}
 \int_X g(x)^p \, d\mu(x) = \int_0^\infty p \, t^{p - 1} \, \lambda(t) \, dt
\end{equation}
for every $p > 0$, by Fubini's theorem.

\section{$L^p$ Estimates}
\label{L^p estimates}

        Let $(X, d(x, y))$ be a metric space again, and let $\mu$ be a 
nonnegative Borel measure on $X$ for which the measure of every open ball
in $X$ is positive and finite.  Suppose that there is a positive real
number $C_1$ such that
\begin{equation}
\label{mu({x in X : M(f)(x) > t}) le C_1 t^{-1} int_X |f(x)| d mu(x)}
 \mu(\{x \in X : M(f)(x) > t\}) \le C_1 \, t^{-1} \, \int_X |f(x)| \, d\mu(x)
\end{equation}
for every integrable function $f$ on $X$ with respect to $\mu$.  This
is the same as (\ref{mu(V_t) le C_1 t^{-1} nu(X)}) in Section
\ref{weak-type estimates}, when $\nu$ corresponds to $f$ as in
(\ref{nu(A) = int_A |f| d mu}) in Section \ref{definitions}.  Let a
real number $p \ge 1$ be given, and suppose now that $f \in L^p(X)$
with respect to $\mu$.  We would like to show that $M(f) \in L^p(X)$
when $p > 1$.

        Put
\begin{eqnarray}
\label{f_t(x) = f(x) when |f(x)| le t, = 0 when |f(x)| > t}
        f_t(x) & = & f(x) \, \hbox{ when } |f(x)| \le t \\
               & = & 0 \qquad\hbox{when } |f(x)| > t \nonumber
\end{eqnarray}
for each $t > 0$.  Thus $f_t$ is a bounded measurable function on $X$,
with $\|f\|_\infty \le t$, so that
\begin{equation}
\label{M(f)(x) le t}
        M(f)(x) \le t
\end{equation}
for every $x \in X$, as in (\ref{sup_{x in X} M(f)(x) le ||f||_infty})
in Section \ref{definitions}.  Observe that
\begin{equation}
\label{M(f)(x) le M(f_{a t})(x) + M(f - f_{a t})(x) le a t + M(f - f_{a t})(x)}
        M(f)(x) \le M(f_{a \, t})(x) + M(f - f_{a \, t})(x)
                 \le a \, t + M(f - f_{a \, t})(x)
\end{equation}
for every $x \in X$ and $a, t > 0$.  If $0 < a < 1$ and $M(f)(x) > t$,
then it follows that $M(f - f_{a \, t})(x) > (1 - a) \, t$, which is
to say that
\begin{equation}
\label{{x in X : M(f)(x) > t} subseteq {x in X : M(f - f_{a t})(x) > (1 - a)t}}
         \{x \in X : M(f)(x) > t\} 
                   \subseteq \{x \in X : M(f - f_{a \, t})(x) > (1- a) \, t\}.
\end{equation}
This implies that
\begin{eqnarray}
\label{mu({x in X : M(f)(x) > t}) le ...}
\lefteqn{\mu(\{x \in X : M(f)(x) > t\})} \\
 & \le & \mu(\{x \in X : M(f - f_{a \, t})(x) > (1 - a) \, t\}) \nonumber \\
 & \le & C_1 \, (1 - a)^{-1} \, t^{-1} \, \int_X |f(x) - f_{a \, t}(x)|
                                             \, d\mu(x) \nonumber \\
 & \le & C_1 \, (1 - a)^{-1} \, t^{-1} \, \int_{\{x \in X : |f(x)| > a \, t\}}
                                         |f(x)| \, d\mu(x). \nonumber
\end{eqnarray}
More precisely, this uses (\ref{mu({x in X : M(f)(x) > t}) le C_1
  t^{-1} int_X |f(x)| d mu(x)}) in the second step, applied to $f -
f_{a \, t}$ instead of $f$, and $(1 - a) \, t$ instead of $t$.  In the
third step, we have used the fact that $f(x) - f_{a \, t}(x)$ is equal
to $f(x)$ when $|f(x)| > a \, t$, and is $0$ otherwise.  Note that $f
- f_{a \, t}$ is integrable on $X$ for every $a, t > 0$, even when $p
> 1$.

        Let us restrict our attention now to the case where $p > 1$.
The integral of $M(f)(x)^p$ with respect to $\mu$ can be expressed as
in (\ref{int_X g(x)^p d mu(x) = int_0^infty p t^{p - 1} lambda(t) dt})
with $g = M(f)$, and we can use (\ref{mu({x in X : M(f)(x) > t}) le
  ...}) to estimate $\lambda(t)$ as in (\ref{lambda(t) = mu({x in X :
    g(x) > t})}).  This implies that
\begin{eqnarray}
\label{int_X M(f)(x)^p d mu(x) le ...}
\lefteqn{\int_X M(f)(x)^p \, d\mu(x)} \\
 & \le & \int_0^\infty p \, t^{p - 1} \, \Big(C_1 \, (1 - a)^{-1} \, t^{-1} \,
 \int_{\{x \in X : |f(x)| > a \, t\}} |f(x)| \, d\mu(x)\Big) \, dt \nonumber \\
 & = & p \, C_1 \, (1 - a)^{-1} \, \int_0^\infty t^{p - 2} \, 
 \Big(\int_{\{x \in X : |f(x)| > a \, t\}} |f(x)| \, d\mu(x)\Big) \, dt.
                                                                    \nonumber
\end{eqnarray}
Interchanging the order of integration, we get that
\begin{eqnarray}
\label{int_X M(f)(x)^p d mu(x) le ..., 2}
\lefteqn{\int_X M(f)(x)^p \, d\mu(x)} \\
 & \le & p \, C_1 \, (1 - a)^{-1} \, \int_X \Big(\int_0^{|f(x)|/a} t^{p - 2}
                                 \, dt\Big) \, |f(x)| \, d\mu(x) \nonumber \\
 & = & p \, C_1 \, (1 - a)^{-1} \, \int_X (p - 1)^{-1} \, (|f(x)|/a)^{p - 1}
                                  \, |f(x)| \, d\mu(x) \nonumber \\
 & = & p \, C_1 \, (1 - a)^{-1} \, (p - 1)^{-1} \, a^{1 - p} \, 
                                      \int_X |f(x)|^p \, d\mu(x). \nonumber
\end{eqnarray}
It follows that $M(f) \in L^p(X)$, with $L^p$ norm bounded by the
$L^p$ norm of $f$ times a constant that depends on $p$ when $p > 1$,
by taking the $p$th root of both sides of (\ref{int_X M(f)(x)^p d
  mu(x) le ..., 2}).  This works using any $a \in (0, 1)$, so that one
can choose an optimal $a$ for each $p$.  In particular, it is better
to take $a$ close to $1$ as $p$ increases.

\section{Conditional expectation}
\label{conditional expectation}

        Let $(X, \mathcal{A}, \mu)$ be a probability space, so that
$X$ is a set, $\mathcal{A}$ is a $\sigma$-algebra of subsets of $X$,
and $\mu$ is a probability measure on $X$, which is to say a
nonnegative countably-additive measure on $\mathcal{A}$ such that
$\mu(X) = 1$.  Suppose that $f$ is a real or complex-valued function
on $X$ that is measurable with respect to $\mathcal{A}$, and
integrable with respect to $\mu$.  Thus
\begin{equation}
\label{nu_f(A) = int_A f d mu}
        \nu_f(A) = \int_A f \, d\mu
\end{equation}
is defined for each $A \in \mathcal{A}$, and determines a
countably-additive real or complex-valued measure on $\mathcal{A}$, as
appropriate.  

        Let $\mathcal{B}$ be another $\sigma$-algebra of subsets of $X$
contained in $\mathcal{A}$, so that $\mathcal{B}$ is a $\sigma$-subalgebra
of $\mathcal{A}$.  The restriction of $\nu_f$ to $\mathcal{B}$ is a
countably-additive real or complex-valued measure on $(X, \mathcal{B})$,
which is absolutely continuous with respect to the restriction of $\mu$
to $\mathcal{B}$.  The Radon--Nikodym theorem implies that there is a
measurable function $f_\mathcal{B}$ on $X$ with respect to $\mathcal{B}$
which is also integrable with respect to $\mu$ such that
\begin{equation}
\label{int_A f_mathcal{B} d mu = nu_f(A) = int_A f d mu}
        \int_A f_\mathcal{B} \, d\mu = \nu_f(A) = \int_A f \, d\mu
\end{equation}
for every $A \in \mathcal{B}$.  This function $f_\mathcal{B}$ is known
as the \emph{conditional expectation}\index{conditional expectation}
of $f$ with respect to $\mathcal{B}$, which may be denoted $E(f \mid
\mathcal{B})$ as well.  If $f'_\mathcal{B}$ is any other function on
$X$ that satisfies the same properties as $f_\mathcal{B}$, then it is
easy to see that $f_\mathcal{B} = f'_\mathcal{B}$ almost everywhere on
$X$ with respect to $\mu$.  Of course, if $\mathcal{B} = \mathcal{A}$,
then we can simply take $f_\mathcal{B} = f$.  The conditional
expectation of $f$ with respect to $\mathcal{B}$ may also be denoted
$E_\mathcal{A}(f \mid \mathcal{B})$, to indicate the initial
$\sigma$-algebra $\mathcal{A}$ explicitly.

        As a basic class of examples, let $\mathcal{P}$ be a partition of
$X$ into finitely or countably many measurable sets with positive measure
with respect to $\mu$.  Thus $\mathcal{P}$ is a collection of finitely
or countably many pairwise-disjoint elements of $\mathcal{A}$ such
that $\mu(A) > 0$ for each $A \in \mathcal{P}$, and the union of the
elements of $\mathcal{P}$ is equal to $X$.  Also let
$\mathcal{B}(\mathcal{P})$ be the collection of subsets of $X$ that
can be expressed as a union of elements of $\mathcal{P}$, which is
interpreted as including the empty set.  It is easy to see that
$\mathcal{B}(\mathcal{P})$ is a $\sigma$-subalgebra of $\mathcal{A}$,
and that a function $f$ on $X$ is measurable with respect to
$\mathcal{B}(\mathcal{P})$ if and only if $f$ is constant on each $A
\in \mathcal{A}$.  If $f$ is a measurable function on $X$ with respect
to $\mathcal{A}$ which is integrable with respect to $\mu$, then the
conditional expectation $f_\mathcal{B}$ of $f$ with respect to
$\mathcal{B}$ is given by
\begin{equation}
\label{f_mathcal{B}(x) = frac{1}{mu(A)} int_A f d mu}
        f_\mathcal{B}(x) = \frac{1}{\mu(A)} \, \int_A f \, d\mu
\end{equation}
for every $A \in \mathcal{B}$ and $x \in A$.

        As another class of examples, let $(X_1, \mathcal{A}_1, \mu_1)$
and $(X_2, \mathcal{A}_2, \mu_2)$ be probability spaces, and consider
their Cartesian product $X = X_1 \times X_2$.  The standard product
measure construction leads to a $\sigma$-algebra $\mathcal{A}$ on $X$,
and a probability measure $\mu$ defined on $\mathcal{A}$.  Let
$\mathcal{B}_1$ be the collection of subsets of $X$ of the form $A
\times X_2$, where $A \in \mathcal{A}_1$ This is a $\sigma$-subalgebra
of $\mathcal{A}$, and a function $f(x) = f(x_1, x_2)$ on $X$ is
measurable with respect to $\mathcal{B}_1$ if and only if $f(x_1,
x_2)$ only depends on $x_1$, and this function of $x_1$ is measurable
with respect to $\mathcal{A}_1$ as a function on $X_1$.  If $f$ is a
function on $X$ which is measurable with respect to $\mathcal{A}$ and
integrable with respect to $\mu$, then the conditional expectation
$f_{\mathcal{B}_1}$ of $f$ with respect to $\mathcal{B}_1$ is given by
\begin{equation}
\label{f_{mathcal{B}_1}(x_1, x_2) = int_{X_2} f(x_1, y_2) d mu_2(y_2)}
        f_{\mathcal{B}_1}(x_1, x_2) = \int_{X_2} f(x_1, y_2) \, d\mu_2(y_2),
\end{equation}
essentially by Fubini's theorem.

        Let $\mathcal{A}$ be any $\sigma$-algebra of subsets of a set $X$
again, and let $\nu$ be a real or complex measure defined on $\mathcal{A}$.
Remember that the corresponding total variation measure $|\nu|$ is defined
on $\mathcal{A}$ by
\begin{equation}
\label{|nu|(A) = sup sum_{j = 1}^infty |nu(A_j)|}
        |\nu|(A) = \sup \sum_{j = 1}^\infty |\nu(A_j)|,
\end{equation}
where the supremum is taken over all sequences $A_1, A_2, A_3, \ldots$
of pairwise-disjoint measurable subsets of $X$ whose union is equal to
$A$.  It is well known that $|\nu|$ is a countably-additive
nonnegative measure defined on $\mathcal{A}$, and that $|\nu|(X) <
+\infty$.

        Suppose that $\mathcal{B}$ is a $\sigma$-subalgebra of $\mathcal{A}$,
and let $\nu_\mathcal{B}$ be the restriction of $\nu$ to $\mathcal{B}$, which 
is a countably-additive real or complex measure defined on $\mathcal{B}$.
Observe that
\begin{equation}
\label{|nu_mathcal{B}|(A) le |nu|(A)}
        |\nu_\mathcal{B}|(A) \le |\nu|(A)
\end{equation}
for every $A \in \mathcal{B}$, where $|\nu|$ is the total variation of
$\nu$ as a measure on $\mathcal{A}$, and $|\nu_\mathcal{B}|$ is the
total variation of $\nu_\mathcal{B}$ as a measure on $\mathcal{B}$.
More precisely, if $A \in \mathcal{B}$, then $|\nu_\mathcal{B}|(A)$ is
the supremum of the same type of sums as in (\ref{|nu|(A) = sup sum_{j
    = 1}^infty |nu(A_j)|}), but where the $A_j$'s are required to be
in $\mathcal{B}$.  Thus $|\nu|(A)$ is given by a supremum of sums that
includes the sums whose supremum is equal to $|\nu_\mathcal{B}|(A)$
when $A \in \mathcal{B}$, which implies (\ref{|nu_mathcal{B}|(A) le
  |nu|(A)}).

        Let $\mu$ be a probability measure defined on $\mathcal{A}$,
and let $f$ be a real or complex-valued function on $X$ which is
measurable with respect to $\mathcal{A}$ and integrable with respect
to $\mu$.  If $\nu = \nu_f$ is as in (\ref{nu_f(A) = int_A f d mu}),
then it is well known that
\begin{equation}
\label{|nu|(A) = int_A |f| d mu}
        |\nu|(A) = \int_A |f| \, d\mu
\end{equation}
for every $A \in \mathcal{B}$.  Let $\mathcal{B}$ be a
$\sigma$-subalgebra of $\mathcal{A}$, and let $\nu_\mathcal{B}$ be the
restriction of $\nu$ to $\mathcal{B}$, as in the previous paragraph.
If $f_\mathcal{B}$ is the conditional expectation of $f$ with respect
to $\mathcal{B}$, then $\nu_\mathcal{B}(A)$ is equal to the integral
of $f_\mathcal{B}$ over $A$ with respect to $\mu$ for every $A \in
\mathcal{B}$, as in (\ref{int_A f_mathcal{B} d mu = nu_f(A) = int_A f
  d mu}), and hence
\begin{equation}
\label{|nu_mathcal{B}|(A) = int_A |f_mathcal{B}| d mu}
        |\nu_\mathcal{B}|(A) = \int_A |f_\mathcal{B}| \, d\mu
\end{equation}
for every $A \in \mathcal{B}$, as in (\ref{|nu|(A) = int_A |f| d mu}).
It follows that
\begin{equation}
\label{int_A |f_mathcal{B}| d mu le int_A |f| d mu}
        \int_A |f_\mathcal{B}| \, d\mu \le \int_A |f| \, d\mu
\end{equation}
for every $A \in \mathcal{B}$, because of (\ref{|nu_mathcal{B}|(A) le
  |nu|(A)}).

        Even if a real or complex measure $\nu$ defined on $\mathcal{A}$
is not absolutely continuous with respect to $\mu$, it may be that the
restiction $\nu_\mathcal{B}$ of $\nu$ to a $\sigma$-subalgebra
$\mathcal{B}$ of $\mathcal{A}$ is absolutely continuous with respect
to the restriction of $\mu$ to $\mathcal{A}$.  Under these conditions,
the Radon--Nikodym theorem again implies that $\nu_\mathcal{B}$ can be
represented on $\mathcal{B}$ by integration of a function
$f_\mathcal{B}$ on $X$ that is measurable with respect to
$\mathcal{B}$ and integrable with respect to $\mu$.  As before, one
can combine (\ref{|nu_mathcal{B}|(A) le |nu|(A)}) and
(\ref{|nu_mathcal{B}|(A) = int_A |f_mathcal{B}| d mu}) to get that
\begin{equation}
\label{int_A |f_mathcal{B}| d mu le |nu|(A)}
        \int_A |f_\mathcal{B}| \, d\mu \le |\nu|(A)
\end{equation}
for every $A \in \mathcal{B}$.  If $\mathcal{B} =
\mathcal{B}(\mathcal{P})$ is the $\sigma$-algebra generated by a
partition $\mathcal{P}$ of $X$ into finitely or countably many
elements of $\mathcal{A}$, each of which has positive measure with
respect to $\mu$, then every measure defined on $\mathcal{B}$ is
absolutely continuous with respect to the restriction of $\mu$ to
$\mathcal{B}$.  In this case, we have that
\begin{equation}
\label{f_mathcal{B}(x) = frac{nu(A)}{mu(A)}}
        f_\mathcal{B}(x) = \frac{\nu(A)}{\mu(A)}
\end{equation}
for every $A \in \mathcal{B}$ and $x \in A$, instead of
(\ref{f_mathcal{B}(x) = frac{1}{mu(A)} int_A f d mu}).

        Suppose now that $\mathcal{B}$ and $\mathcal{C}$ are 
$\sigma$-subalgebras of $\mathcal{A}$, with $\mathcal{B} \subseteq 
\mathcal{C}$, and let $\nu$ be a real or complex-valued measure
defined on $\mathcal{A}$.  If $\nu_\mathcal{B}$ and $\nu_\mathcal{C}$
are the restrictions of $\nu$ to $\mathcal{B}$, $\mathcal{C}$,
respectively, then $\nu_\mathcal{B}$ is also the same as the
restriction of $\nu_\mathcal{C}$ as a measure defined on $\mathcal{C}$
to a measure on $\mathcal{B}$.  This is basically trivial, but it has
the following nice interpretation for conditional expectation.  Let
$f$ be a real or complex-valued function on $X$ that is measurable
with respect to $\mathcal{A}$ and integrable with respect to $\mu$,
and let $f_\mathcal{B} = E_\mathcal{A}(f \mid \mathcal{B})$,
$f_\mathcal{C} = E_\mathcal{A}(f \mid \mathcal{C})$ be the conditional
expectations of $f$ with respect to $\mathcal{B}$, $\mathcal{C}$,
respectively.  Also let
\begin{equation}
\label{(f_mathcal{C})_mathcal{B} = E_mathcal{C}(f_mathcal{C} mid mathcal{B})}
 (f_\mathcal{C})_\mathcal{B} = E_\mathcal{C}(f_\mathcal{C} \mid \mathcal{B})
\end{equation}
be the conditional expectation of $f_\mathcal{C}$ with respect to
$\mathcal{B}$, where $f_\mathcal{C}$ is considered as a measurable
function with respect to $\mathcal{C}$ instead of $\mathcal{A}$.
Under these conditions, it is easy to see that
\begin{equation}
\label{(f_mathcal{C})_mathcal{B} = f_mathcal{B}}
        (f_\mathcal{C})_\mathcal{B} = f_\mathcal{B},
\end{equation}
This follows from the previous statement about measures, applied to
$\nu = \nu_f$ as in (\ref{nu_f(A) = int_A f d mu}).  In particular, if
$f$ is already measurable with respect to $\mathcal{C}$, then
\begin{equation}
\label{E_mathcal{A}(f mid mathcal{B}) = E_mathcal{C}(f mid mathcal{B})}
        E_\mathcal{A}(f \mid \mathcal{B}) = E_\mathcal{C}(f \mid \mathcal{B}),
\end{equation}
as in the case of $f_\mathcal{C}$ in (\ref{(f_mathcal{C})_mathcal{B} =
  E_mathcal{C}(f_mathcal{C} mid mathcal{B})}).

\section{Additional properties}
\label{additional properties}

        Let $(X, \mathcal{A}, \mu)$ be a probability space again, and let
$\mathcal{B}$ be a $\sigma$-subalgebra of $\mathcal{A}$.  Also let
$f$ be a real or complex-valued function on $X$ that is measurable with
respect to $\mathcal{A}$ and integrable with respect to $\mu$, and let
$f_\mathcal{B}$ be the conditional expectation of $f$ with respect to
$\mathcal{B}$.  Of course, $|f|$ is a nonnegative real-valued integrable
function on $X$, and so the conditional expectation $(|f|)_\mathcal{B}$
of $|f|$ with respect to $\mathcal{B}$ is real-valued and nonnegative as
well.  Observe that
\begin{equation}
\label{int_A |f_mathcal{B}| d mu le ... = int_A (|f|)_mathcal{B} d mu}
        \int_A |f_\mathcal{B}| \, d\mu \le \int_A |f| \, d\mu
                              = \int_A (|f|)_\mathcal{B} \, d\mu
\end{equation}
for every $A \in \mathcal{B}$, by (\ref{int_A |f_mathcal{B}| d mu le
  int_A |f| d mu}) and the definition of $(|f|)_\mathcal{B}$.  It
follows that
\begin{equation}
\label{|f_mathcal{B}| le (|f|)_mathcal{B}}
        |f_\mathcal{B}| \le (|f|)_\mathcal{B}
\end{equation}
almost everywhere on $X$ with respect to $\mu$, since both sides of the
inequality are measurable with respect to $\mathcal{B}$.

        Now let $p \in (1, +\infty)$ be given, and suppose that $|f|^p$
is integrable on $X$ with respect to $\mu$.  Thus the conditional
expectation $(|f|^p)_\mathcal{B}$ of $|f|^p$ with respect to
$\mathcal{B}$ can be defined as before, and is real-valued and
nonnegative.  If $A \in \mathcal{B}$ and $\mu(A) > 0$, then
\begin{eqnarray}
\label{(frac{1}{mu(A)} int_A |f_mathcal{B}| d mu)^p le ...}
 \Big(\frac{1}{\mu(A)} \, \int_A |f_\mathcal{B}| \, d\mu\Big)^p
     & \le & \Big(\frac{1}{\mu(A)} \, \int_A |f| \, d\mu\Big)^p \\
     & \le & \frac{1}{\mu(A)} \, \int_A |f|^p \, d\mu \nonumber \\
      & = &  \frac{1}{\mu(A)} \, \int_A (|f|^p)_\mathcal{B} \, d\mu. \nonumber
\end{eqnarray}
This uses (\ref{int_A |f_mathcal{B}| d mu le int_A |f| d mu}) in the
first step, Jensen's or H\"older's inequality in the second step, and
the definition of $(|f|^p)_\mathcal{B}$ in the third step.  One can
check that this implies that
\begin{equation}
\label{(|f_mathcal{B}|)^p le (|f|^p)_mathcal{B}}
        (|f_\mathcal{B}|)^p \le (|f|^p)_\mathcal{B}
\end{equation}
almost everywhere on $X$, because $|f_\mathcal{B}|$ and
$(|f|^p)_\mathcal{B}$ are both measurable with respect to
$\mathcal{B}$.  It follows that
\begin{equation}
\label{int_X (|f_mathcal{B}|)^p d mu le ... = int_X |f|^p d mu}
 \int_X (|f_\mathcal{B}|)^p \, d\mu \le \int_X (|f|^p)_\mathcal{B} \, d\mu
                                     = \int_X |f|^p \, d\mu,
\end{equation}
using the definition of $(|f|^p)_\mathcal{B}$ in the second step.  In
particular, $|f_\mathcal{B}|^p$ is integrable with respect to $\mu$ as
well.

        Let $L^p(X, \mathcal{A}, \mu)$ be the usual space of real
or complex-valued functions $f$ on $X$ such that $f$ is measurable
with respect to $\mathcal{A}$ and $|f|^p$ is integrable with respect
to $\mu$, for $1 \le p < \infty$.  More precisely, $L^p(X,
\mathcal{A}, \mu)$ consists of equivalence classes of such functions,
which are equal to each other almost everywhere with respect to $\mu$
on $X$.  It is well known that $L^p(X, \mathcal{A}, \mu)$ is complete
with respect to the metric associated to the the $L^p$ norm
\begin{equation}
\label{||f||_p = (int_X |f|^p d mu)^{1/p}}
        \|f\|_p = \Big(\int_X |f|^p \, d\mu\Big)^{1/p}.
\end{equation}
Similarly, $L^\infty(X, \mathcal{A}, \mu)$ consists of equivalence
classes of functions on $X$ that are measurable with respect to
$\mathcal{A}$ and essentially bounded on $X$.  The $L^\infty$ norm
$\|f\|_\infty$ is defined to be the essential supremum of $|f|$ on
$X$, and $L^\infty(X, \mathcal{A}, \mu)$ is complete with respect to
the metric associated to this norm.

        Of course, if $f$ is measurable with respect to $\mathcal{B}$, 
then $f$ is measurable with respect to $\mathcal{A}$ too.  This leads
to a natural linear mapping from $L^p(X, \mathcal{B}, \mu)$ into
$L^p(X, \mathcal{A}, \mu)$ for each $p$, $1 \le p \le \infty$, which
is an isometry with respect to the $L^p$ norm.  Note that the image of
$L^p(X, \mathcal{B}, \mu)$ under this mapping is a closed linear
subspace of $L^p(X, \mathcal{A}, \mu)$, because $L^p(X, \mathcal{B},
\mu)$ is complete.

        It is easy to see that $f \mapsto f_\mathcal{B}$ is a linear
mapping from $L^1(X, \mathcal{A}, \mu)$ into $L^1(X, \mathcal{B},
\mu)$, using the uniqueness of $f_\mathcal{B}$.  Moreover,
\begin{equation}
\label{int_X |f_mathcal{B}| d mu le int_X |f| d mu}
        \int_X |f_\mathcal{B}| \, d\mu \le \int_X |f| \, d\mu
\end{equation}
for every $f \in L^1(X, \mathcal{A}, \mu)$, by (\ref{int_A
  |f_mathcal{B}| d mu le int_A |f| d mu}) with $A = X$.  If $f \in
L^p(X, \mathcal{A}, \mu)$ and $1 < p < \infty$, then $f_\mathcal{B}
\in L^p(X, \mathcal{B}, \mu)$ and
\begin{equation}
\label{||f_mathcal{B}||_p le ||f||_p}
        \|f_\mathcal{B}\|_p \le \|f\|_p,
\end{equation}
by (\ref{int_X (|f_mathcal{B}|)^p d mu le ... = int_X |f|^p d mu}).
Similarly, if $f \in L^\infty(X, \mathcal{A}, \mu)$, then
\begin{equation}
\label{int_A |f_mathcal{B}| d mu le int_A |f| d mu le ||f||_infty mu(A)}
        \int_A |f_\mathcal{B}| \, d\mu \le \int_A |f| \, d\mu 
                                        \le \|f\|_\infty \, \mu(A)
\end{equation}
for every $A \in \mathcal{B}$, using (\ref{int_A |f_mathcal{B}| d mu
  le int_A |f| d mu}) in the first step.  This implies that
$f_\mathcal{B}$ is essentially bounded on $X$ as well, and that
(\ref{||f_mathcal{B}||_p le ||f||_p}) holds when $p = \infty$, because
$f_\mathcal{B}$ is measurable with respect to $\mathcal{B}$.

        Let $f \in L^1(X, \mathcal{A}, \mu)$ and $B \in \mathcal{B}$
be given, and let ${\bf 1}_B(x)$ be the characteristic or indicator
function on $X$ associated to $B$, which is equal to $1$ when $x \in B$
and to $0$ otherwise.  If $A \in \mathcal{B}$, then $A \cap B \in \mathcal{B}$
too, and hence
\begin{equation}
\label{int_A f_mathcal{B} {bf 1}_B d mu = ... = int_A f {bf 1}_B d mu}
  \int_A f_\mathcal{B} \, {\bf 1}_B \, d\mu
     = \int_{A \cap B} f_\mathcal{B} \, d\mu = \int_{A \cap B} f \, d\mu
                                         = \int_A f \, {\bf 1}_B \, d\mu,
\end{equation}
using the definition of $f_\mathcal{B}$ in the second step.  Of
course, ${\bf 1}_B$ is measurable with respect to $\mathcal{B}$ on
$X$, because $B \in \mathcal{B}$, and $f_\mathcal{B}$ is measurable
with respect to $\mathcal{B}$ by construction, so that $f_\mathcal{B}
\, {\bf 1}_B$ is measurable with respect to $\mathcal{B}$ as well.
This shows that $f_\mathcal{B} \, {\bf 1}_B$ is equal to the
conditional expectation of $f \, {\bf 1}_B$ with respect to
$\mathcal{B}$, since it satisfies the requirements of the conditional
expectation.

        Similarly, if $g \in L^\infty(X, \mathcal{B}, \mu)$, then
\begin{equation}
\label{(f g)_mathcal{B} = f_mathcal{B} g}
        (f \, g)_\mathcal{B} = f_\mathcal{B} \, g.
\end{equation}
This reduces to the discussion in the previous paragraph when $g =
{\bf 1}_B$ for some $B \in \mathcal{B}$, which implies that (\ref{(f
  g)_mathcal{B} = f_mathcal{B} g}) holds when $g$ is a simple function
on $X$ that is measurable with respect to $\mathcal{B}$, by linearity.
One can use this to get that (\ref{(f g)_mathcal{B} = f_mathcal{B} g})
holds for every $g \in L^\infty(X, \mathcal{B}, \mu)$, because simple
functions on $X$ that are measurable with respect to $\mathcal{B}$ are
dense in $L^\infty(X, \mathcal{B}, \mu)$ with respect to the
$L^\infty$ norm.  If $f \in L^p(X, \mathcal{A}, \mu)$ for some $p$, $1
\le p \le \infty$, and if $q$ is the exponent conjugate to $p$, in the
sense that $1 \le q \le \infty$ and $1/p + 1/q = 1$, then one can
check that (\ref{(f g)_mathcal{B} = f_mathcal{B} g}) holds for every
$g \in L^q(X, \mathcal{B}, \mu)$.  This also uses H\"older's
inequality, and the fact that $f_\mathcal{B} \in L^p(X, \mathcal{B},
\mu)$ when $f \in L^p(X, \mathcal{A}, \mu)$ and $1 \le p \le \infty$,
as before.

\section{Another maximal function}
\label{another maximal function}

        Let $(X, \mathcal{A}, \mu)$ be a probability space, and let
$\mathcal{B}_1, \ldots, \mathcal{B}_n$ be finitely many $\sigma$-subalgebras
of $\mathcal{A}$, with $\mathcal{B}_j \subseteq \mathcal{B}_{j + 1}$
for $j = 1, \ldots, n - 1$.  Also let $f \in L^1(X, \mathcal{A}, \mu)$
be given, and let
\begin{equation}
\label{f_j = f_{mathcal{B}_j} = E(f mid mathcal{B}_j)}
        f_j = f_{\mathcal{B}_j} = E(f \mid \mathcal{B}_j)
\end{equation}
be the conditional expectation of $f$ with respect to $\mathcal{B}_j$
for each $j$.  Put
\begin{equation}
\label{f_l^*(x) = max_{1 le j le l} |f_j(x)|}
        f_l^*(x) = \max_{1 \le j \le l} |f_j(x)|
\end{equation}
for each $l = 1, \ldots, n$, and observe that $f_l^*$ is measurable
with respect to $\mathcal{B}_l$ for each $l$.  If $g \in L^1(X,
\mathcal{A}, \mu)$ too, then it is easy to see that
\begin{equation}
\label{(f + g)_l^*(x) le f_l^*(x) + g_l^*(x)}
        (f + g)_l^*(x) \le f_l^*(x) + g_l^*(x)
\end{equation}
for each $l$.  Similarly,
\begin{equation}
\label{(t f)_l^*(x) = |t| f_l^*(x)}
        (t \, f)_l^*(x) = |t| \, f_l^*(x)
\end{equation}
for each $l$ and real or complex number $t$, as appropriate, so that
the mapping from $f$ to $f_l^*$ is sublinear.

        Let $t > 0$ be given, and put
\begin{equation}
\label{A_l(t) = {x in X : f_l^*(x) > t}}
        A_l(t) = \{x \in X : f_l^*(x) > t\}
\end{equation}
for each $l = 1, \ldots, n$, which is an element of $\mathcal{B}_l$,
because $f_l^*$ is measurable with respect to $\mathcal{B}_l$.  Note
that $f_1^* = |f_1|$, so that
\begin{equation}
\label{A_1(t) = {x in X : |f_1(x)| > t}}
        A_1(t) = \{x \in X : |f_1(x)| > t\},
\end{equation}
and hence
\begin{equation}
\label{mu(A_1(t)) le ... le t^{-1} int_{A_1(t)} |f| d mu}
        \mu(A_1(t)) \le t^{-1} \, \int_{A_1(t)} |f_1| \, d\mu
                      \le t^{-1} \, \int_{A_1(t)} |f| \, d\mu,
\end{equation}
using (\ref{int_A |f_mathcal{B}| d mu le |nu|(A)}) in the second step.
If $l > 1$, then
\begin{eqnarray}
\label{A_l(t) setminus A_{l - 1}(t) = ...}
 A_l(t) \setminus A_{l - 1}(t) & = & \{x \in X : f_{l - 1}^*(x) \le t, 
                                                     \, f_l^*(x) > t\} \\
 & = & \{x \in X : f_{l - 1}^*(x) \le t, \, |f_l(x)| > t\}, \nonumber
\end{eqnarray}
by the definition of $f_l^*$, and in particular $|f_l(x)| > t$ on
(\ref{A_l(t) setminus A_{l - 1}(t) = ...}).  Thus
\begin{eqnarray}
\label{mu(A_l(t) setminus A_{l - 1}(t)) le ...}
        \mu(A_l(t) \setminus A_{l - 1}(t))
 & \le & t^{-1} \, \int_{A_l(t) \setminus A_{l - 1}(t)} |f_l| \, d\mu \\
 & \le & t^{-1} \, \int_{A_l(t) \setminus A_{l - 1}(t)} |f| \, d\mu, \nonumber
\end{eqnarray}
again using (\ref{int_A |f_mathcal{B}| d mu le int_A |f| d mu}) in the
second step, and the fact that $A_l(t) \setminus A_{l - 1}(t) \in
\mathcal{B}_l$, since $A_l(t) \in \mathcal{B}_l$ and $A_{l - 1}(t) \in
\mathcal{B}_{l - 1} \subseteq \mathcal{B}_l$.

        By construction, $f_l^*$ is monotone increasing in $l$, which
implies that
\begin{equation}
\label{A_l(t) subseteq A_{l + 1}(t)}
        A_l(t) \subseteq A_{l + 1}(t)
\end{equation}
for $l = 1, \ldots, n - 1$.  It follows that the sets $A_j(t)
\setminus A_{j - 1}(t)$ are pairwise disjoint for $j \ge 2$, and
disjoint from $A_1(t)$.  Using (\ref{mu(A_1(t)) le ... le t^{-1}
  int_{A_1(t)} |f| d mu}) and (\ref{mu(A_l(t) setminus A_{l - 1}(t))
  le ...}), we get that
\begin{eqnarray}
\label{mu(A_l(t)) = ... le t^{-1} int_X |f| d mu}
 \mu(A_l(t)) & = & \mu(A_1(t)) +
                     \sum_{j = 2}^l \mu(A_j(t) \setminus A_{j - 1}(t)) \\
 & \le & t^{-1} \, \int_{A_1(t)} |f| \, d\mu
  + \sum_{j = 2}^l t^{-1} \, \int_{A_j(t) \setminus A_{j - 1}(t)} |f| \, d\mu
                                                    \nonumber \\
  & = & t^{-1} \, \int_{A_l(t)} |f| \, d\mu \le t^{-1} \, \int_X |f| \, d\mu
                                                         \nonumber
\end{eqnarray}
for each $l$, with the obvious simplifications when $l = 1$.

        If $f \in L^\infty(X, \mathcal{A}, \mu)$, then 
$f_l^* \in L^\infty(X, \mathcal{B}_l, \mu)$ and
\begin{equation}
\label{||f_l^*||_infty le ||f||_infty}
        \|f_l^*\|_\infty \le \|f\|_\infty
\end{equation}
for each $l$, by the $p = \infty$ version of (\ref{||f_mathcal{B}||_p
  le ||f||_p}).  Using this and the weak-type estimate on $L^1$ in
(\ref{mu(A_l(t)) = ... le t^{-1} int_X |f| d mu}), one can get $L^p$
estimates for $f_l^*$ when $1 < p < \infty$, as in Section \ref{L^p
  estimates}, with $C_1 = 1$.

\backmatter

\newpage

\addcontentsline{toc}{chapter}{Index}

\printindex


\begin{thebibliography}{79}
%\begin{thebibliography}{ABC}

\addcontentsline{toc}{chapter}{Bibliography}


\bibitem {a} P.~Assouad, {\it Plongements lipschitziens dans ${\bf
    R}^n$}, Bulletin de la Soci\'et\'e Math\'ematique de France {\bf
  111} (1983), 429--448.

\bibitem {a-m} M.~Atiyah and I.~Macdonald, {\it Introduction to
  Commutative Algebra}, Addison-Wesley, 1969.

\bibitem {rb} R.~Bass, {\it Real Analysis for Graduate Students:
  Measure and Integration Theory}, 2011.
http://homepages.uconn.edu/\~{}rib02005/real.html.

\bibitem {b-c} J.~Benedetto and W.~Czaja, {\it Integration and Modern
  Analysis}, Birkh\"auser, 2009.

\bibitem {b-mac} G.~Birkhoff and S.~Mac Lane, {\it A Survey of Modern
  Algebra}, 4th edition, Macmillan, 1977.

\bibitem {cas} J.~Cassels, {\it Local Fields}, Cambridge University
  Press, 1986.

\bibitem {c-w-1} R.~Coifman and G.~Weiss, {\it Analyse Harmonique
  Non-Commutative sur certains Espaces Homog\`enes}, Lecture Notes in
  Mathematics {\bf 242}, Springer-Verlag, 1971.

\bibitem {c-w-2} R.~Coifman and G.~Weiss, {\it Extensions of Hardy
  spaces and their use in analysis}, Bulletin of the American
  Mathematical Society {\bf 83} (1977), 569--645.

\bibitem {d-s} G.~David and S.~Semmes, {\it Fractured Fractals and
  Broken Dreams: Self-Similar Geometry through Metric and Measure},
  Oxford University Press, 1997.

\bibitem {e-g} Evans and Gariepy, {\it Measure Theory and Fine
  Properties of Functions}, CRC Press, 1991.

\bibitem {fa1} K.~Falconer, {\it The Geometry of Fractal Sets},
  Cambridge University Press, 1986.

\bibitem {fa2} K.~Falconer, {\it Fractal Geometry: Mathematical
  Foundations and Applications}, 2nd edition, Wiley, 2003.

\bibitem {fe} H.~Federer, {\it Geometric Measure Theory},
  Springer-Verlag, 1969.

\bibitem {gf1} G.~Folland, {\it A Course in Abstract Harmonic
  Analysis}, CRC Press, 1995.

\bibitem {gf2} G.~Folland, {\it Real Analysis}, 2nd edition, Wiley, 1999.

\bibitem {fg} F.~Gouv\^ea, {\it $p$-Adic Numbers: An Introduction},
  2nd edition, Springer-Verlag, 2007.

\bibitem {jh1} J.~Heinonen, {\it Lectures on Analysis on Metric
  Spaces}, Springer-Verlag, 2001.

\bibitem {jh2} J.~Heinonen, {\it Geometric embeddings of metric
  spaces}, Reports of the Department of Mathematics and Statistics
  {\bf 90}, University of Jyv\"askyl\"a, 2003.

\bibitem {h-r} E.~Hewitt and K.~Ross, {\it Abstract Harmonic
  Analysis}, Volumes I and II, 1970, 1979.

\bibitem {h-s} E.~Hewitt and K.~Stromberg, {\it Real and Abstract
  Analysis}, Springer-Verlag, 1975.

\bibitem {h-w} W.~Hurewicz and H.~Wallman, {\it Dimension Theory},
  Princeton University Press, 1969.

\bibitem {fj} F.~Jones, {\it Lebesgue Integration on Euclidean Space},
  Jones and Bartlett, 1993.

\bibitem {ktz} Y.~Katznelson, {\it An Introduction to Harmonic
  Analysis}, 3rd edition, Cambridge University Press, 2004.

\bibitem {jk} J.~Kelley, {\it General Topology}, Springer-Verlag, 1975.

\bibitem {k-s} J.~Kelley and T.~Srinivasan, {\it Measure and
  Integral}, Springer-Verlag, 1988.

\bibitem {sk} S.~Krantz, {\it A Panorama of Harmonic Analysis},
  Mathematical Association of America, 1999.

\bibitem {k-p} S.~Krantz and H.~Parks, {\it The Geometry of Domains
  in Space}, Birkh\"auser, 1999.

\bibitem {l-m} J.~Luukkainen and H.~Movahedi-Lankarani, {\it Minimal
  bi-Lipschitz embedding dimension of ultrametric spaces}, Fundamenta
  Mathematicae {\bf 144} (1994), 181--193.

\bibitem {m-s-1} R.~Mac\'{\i}as and C.~Segovia, {\it Lipschitz
  functions on spaces of homogeneous type}, Advances in Mathematics
  {\bf 33} (1979), 257--270.

\bibitem {m-s-2} R.~Mac\'{\i}as and C.~Segovia, {\it A decomposition
  into atoms of distributions on spaces of homogeneous type}, Advances
  in Mathematics {\bf 33} (1979), 271--309.

\bibitem {mac-b} S.~Mac Lane and G.~Birkhoff, {\it Algebra}, 3rd
  edition, Chelsea, 1988.

\bibitem {mat} P.~Mattila, {\it Geometry of Sets and Measures in
  Euclidean Spaces}, Cambridge University Press, 1995.

\bibitem {r} H.~Royden, {\it Real Analysis}, 3rd edition, Macmillan,
  1988.

\bibitem {r1} W.~Rudin, {\it Principles of Mathematical Analysis}, 3rd
  edition, McGraw-Hill, 1976.

\bibitem {r2} W.~Rudin, {\it Real and Complex Analysis}, 3rd edition,
  McGraw-Hill, 1987.

\bibitem {r3} W.~Rudin, {\it Fourier Analysis on Groups}, Wiley, 1990.

\bibitem {r4} W.~Rudin, {\it Functional Analysis}, 2nd edition,
  McGraw-Hill, 1991.

\bibitem {cs} C.~Sadosky, {\it Interpolation of Operators and Singular
  Integrals: An Introduction to Harmonic Analysis}, Dekker, 1979.

\bibitem {jps} J.-P.~Serre, {\it Lie Algebras and Lie Groups}, Lecture
  Notes in Mathematics {\bf 1500}, Springer-Verlag, 2006.

\bibitem {st1} E.~Stein, {\it Singular Integrals and Differentiability
  Properties of Functions}, Princeton University Press, 1970.

\bibitem {st2} E.~Stein, {\it Topics in Harmonic Analysis Related to
  the Littlewood--Paley Theory}, Annals of Mathematics Studies {\bf
  63}, Princeton University Press, 1970.

\bibitem {st3} E.~Stein, {\it Harmonic Analysis: Real-Variable
  Methods, Orthogonality, and Oscillatory Integrals}, with the
  assistance of T.~Murphy, Princeton University Press, 1993.

\bibitem {st-sh-1} E.~Stein and R.~Shakarchi, {\it Fourier Analysis:
  An Introduction}, Princeton University Press, 2003.

\bibitem {st-sh-2} E.~Stein and R.~Shakarchi, {\it Real Analysis:
  Measure Theory, Integration, and Hilbert Spaces}, Princeton
  University Press, 2005.

\bibitem {st-sh-3} E.~Stein and R.~Shakarchi, {\it Functional
  Analysis: Introduction to Further Topics in Analysis}, Princeton
  University Press, 2011.

\bibitem {s-w} E.~Stein and G.~Weiss, {\it Introduction to Fourier
  Analysis on Euclidean Spaces}, Princeton University Press, 1971.

\bibitem {ks1} K.~Stromberg, {\it Introduction to Classical Real
  Analysis}, Wadsworth, 1981.

\bibitem {ks2} K.~Stromberg, {\it Probability for Analysts}, Lecture
  notes prepared by K.~Ravindran, Chapman \& Hall, 1994.

\bibitem {ds1} D.~Stroock, {\it Probability Theory: An Analytic View},
  2nd edition, Cambridge University Press, 2010.

\bibitem {ds2} D.~Stroock, {\it Essentials of Integration Theory for
  Analysis}, Springer-Verlag, 2011.

\bibitem {ds} D.~Sullivan, {\it Linking the universalities of
  Milnor--Thurston, Feigenbaum, and Ahlfors--Bers}, in {\it
  Topological Methods in Modern Mathematics}, 543--564, Publish or
  Perish, 1993.

\bibitem {t} M.~Taibleson, {\it Fourier Analysis on Local Fields},
  Princeton University Press, 1975.

\bibitem {at} A.~Torchinsky, {\it Real-Variable Methods in Harmonic
  Analysis}, Dover, 2004.

\bibitem {w-z} R.~Wheeden and A.~Zygmund, {\it Measure and Integral:
  An Introduction to Real Analysis}, Dekker, 1977.

\bibitem {z} A.~Zygmund, {\it Trigonometric Series}, Volumes I and II,
  3rd edition, with a foreword by R.Fefferman, Cambridge University
  Press, 2003.



\end{thebibliography}
\end{document}